\documentclass[a4paper]{article}

\usepackage{MyLaTeX,MyArticle}

\RequirePackage{pgf,pgffor}
\usepackage{tikz}
\newcommand{\Proj}[1]{\mathcal{P}({#1})}
\numberwithin{equation}{section}
\newcommand{\CSp}{\operatorname{CSp}}
\newcommand{\CSO}{\operatorname{CSO}}
\newcommand{\CO}{\operatorname{CO}}
\newcommand{\PCO}{\operatorname{PCO}}
\newcommand{\Irr}{\mathrm{Irr}}
\newcommand{\MOD}[1]{#1\mathrm{-mod}}
\title{On the Cohomology of Deligne--Lusztig Varieties}
\author{David A. Craven}
\date{July 2011}
\begin{document}
\maketitle

\begin{abstract}In this paper, we present a conjecture on the degree of unipotent characters in the cohomology of particular Deligne--Lusztig varieties for groups of Lie type, and derive consequences of it. These degrees are a necessary piece of data in the geometric version of Brou\'e's abelian defect group conjecture, and can be used to verify this geometric conjecture in new cases. The geometric version of Brou\'e's conjecture should produce a more combinatorially defined derived equivalence, called a perverse equivalence. We prove that our conjectural degree is an integer (which is not obvious) and has the correct parity for a perfect isometry, and verify that it induces a perverse equivalence for all unipotent blocks of groups of Lie type with cyclic defect groups, whenever the shape of the Brauer tree is known (i.e., not $E_7$ and $E_8$). It has also been used to find perverse equivalences for some non-cyclic cases. This paper is a contribution to the conjectural description of the exact form of a derived equivalence proving Brou\'e's conjecture for groups of Lie type.
\end{abstract}

\tableofcontents



\section{Introduction}

The cohomology of Deligne--Lusztig varieties has been an object of intense study in recent years (see for example \cite{brouemichel1997}, \cite{cravenrouquier2010un}, \cite{dmr2007} and \cite{dudas2010un}), and there are various conjectures and results concerning its structure; one aspect of particular interest is understanding the degree in the cohomology of the Deligne--Lusztig varieties in which a given unipotent character appears. This problem is of interest not only intrinsically, but also because of its application to Brou\'e's abelian defect group conjecture \cite{cravenrouquier2010un} \cite{rouquier2006}. One of the parameters in a potential perverse equivalence between a unipotent block and its Brauer correspondent is given by this degree, and so a conjecture for it would be of great use in constructing derived equivalences for groups of Lie type. (We define perverse equivalences in Section \ref{sec:pervequiv}, and describe an algorithm that should produce perverse equivalences for groups of Lie type.)

In this article we present a general conjecture on this degree, and give considerable evidence to support our conjecture. Previously, only the cases where the prime $\ell$ divides $q\pm 1$ were conjectured \cite{dmr2007}, and so this offers a considerable extension to the previous understanding of this structure. We begin by defining a modified version of Euler's totient function. Let $\mc F$ denote the set of all polynomials in $q$ that are products of cyclotomic polynomials and powers of $q$.

\begin{defn} For $r$ and $d$ integers with $r\geq 2$, write $\phi_d(r)$ for the number of natural numbers prime to $r$ that are at most $r/d$, and write $\phi_d(1)=1/2$. Set $B_d(q)=2$, $B_d(\Phi_r(q))=\phi(r)+d\phi_d(r)=\deg(\Phi_r(q))+d\phi_d(r)$ for $r\geq 1$, and extend the function $B_d:\mc F\to \Q$ linearly, so that $B_d(fg)=B_d(f)+B_d(g)$.
\end{defn}

Let $G=G(q)$ be a group of Lie type, and let $\ell$ be a prime dividing $|G|$, writing $d$ for the multiplicative order of $q$ modulo $\ell$. (The precise groups considered will be defined in Section \ref{sec:DLvars}.) If $\chi$ is a unipotent character of $G$ ($G$ not of type $^2\!B_2$, $^2\!G_2$ or $^2\!F_4$, which are considered in Section \ref{sec:suzree}), then the degree $f$ of $\chi$ is a polynomial in the set $\mc F$. Hence $B_d(f)$ is defined, and we abuse notation slightly and write $B_d(\chi)$ for this rational number. If $(L,\psi)$ is a unipotent $d$-cuspidal pair for the $\ell$-block containing $\chi$, we define $\pi(\chi)=(B_d(\chi)-B_d(\psi))/d$, another rational number.

Our first conjecture is on the cohomology of particular Deligne--Lusztig varieties; here we are deliberately vague, and provide a much more specific conjecture on this is given in Section \ref{sec:DLvars}.

\begin{conj}\label{conj:DLcohom} Let $\ell$ be a prime dividing $\Phi_d(q)$. If $\chi$ is a unipotent character of $\bar{\Q}_\ell G$, then $\pi(\chi_i)$ is the degree of the cohomology of a Deligne--Lusztig variety in which $\chi$ appears.
\end{conj}

The cohomology of the varieties should provide the perverse equivalences for Brou\'e's conjecture, and so the geometric version of Brou\'e's conjecture becomes the following.

\begin{conj}\label{conj:perverse} If $\chi_1,\dots,\chi_s$ are the unipotent ordinary characters in a unipotent $\ell$-block $B$ of $kG$ with abelian defect group, then there is a perverse equivalence from $D^b(\MOD B)$ to $D^b(\MOD b)$ with perversity function given by $\pi(\chi_i)$, where $b$ is the Brauer correspondent of $B$.
\end{conj}

Again, we are more specific about this conjecture should hold in Section \ref{sec:DLvars}. The first test that Conjectures \ref{conj:DLcohom} and \ref{conj:perverse} might hold is to prove that $\pi(\chi)$ is always an integer, which is the content of our first theorem. 

\begin{thm}\label{thm:pifnisinteger} Let $d\geq 1$ be such that $\Phi_d(q)$ divides $|G(q)|$. If $\chi$ is a unipotent character of $KG$ then $\pi(\chi)$ is an integer.
\end{thm}

The next theorem checks that in a bijection with signs predicted between a unipotent block and its Brauer correspondent, the sign attached to $\chi$ is $(-1)^{\pi(\chi)}$.

\begin{thm}\label{thm:bijectionwithsigns} Let $\ell$ be a large prime dividing $\Phi_d(q)$, and let $B$ be a unipotent $\ell$-block of $kG$, with Brauer correspondent $b$. In a bijection with signs $\Irr(B)\to \Irr(b)$, the sign attached to a unipotent character $\chi$ is $(-1)^{\pi(\chi)}$.
\end{thm}

(The definition of a large prime is, in the split case, that it does not divide the order of the Weyl group. In general, see \cite{bmm1993} or Section \ref{sec:DLvars} for a definition.) We prove Theorems \ref{thm:pifnisinteger} and \ref{thm:bijectionwithsigns} simultaneously in Section \ref{sec:integrality}; the proof is not case-by-case, and is remarkably short, needing no facts about groups of Lie type beyond the statement that $\chi(1)/\psi(1)$ is a constant modulo $\Phi_d(q)$, which is known \cite[\S5]{bmm1993}. In particular, we get a geometric interpretation of $B_d(\Phi_r)$; the quantity $B_d(\Phi_r)\pi/d$ is (modulo $2\pi$) the argument of the complex number $\Phi_r(\zeta)$, where $\zeta=\e^{2\pi\I/d}$. This proof exposes the meaning behind the somewhat obscure function $B_d(\Phi_r)$.

Having established this much we prove that, for unipotent $\ell$-blocks with cyclic defect groups, wherever the Brauer tree is known the $\pi$-function induces a perverse equivalence.

\begin{thm}\label{thm:brauertreesperverse} Let $\ell$ be a prime dividing $\Phi_d(q)$, and let $G$ be not of type $E_7$ or $E_8$. If $B$ is a unipotent $\ell$-block of $G$ with cyclic defect group, then the function $\pi(-)$ induces a perverse equivalence between $B$ and its Brauer correspondent, consistent with the geometric version of Brou\'e's conjecture.
\end{thm}

In the cases of $E_7$ and $E_8$, the reason we do not prove the result is that the Brauer trees are not known in these cases. Using the $\pi$-function however, it is possible to produce guesses for the Brauer trees in these cases, and if these guesses are true then in these cases as well the theorem holds. For primes dividing $\Phi_d$ for certain $d$, one can prove this result for $E_7$ and $E_8$ even without explicit knowledge of the Brauer tree (e.g., $d=3$ and $G=E_7$).

As well as always finding a perverse equivalence for blocks with cyclic defect groups, we also consider a few previous conjectures and results about the degrees in Deligne--Lusztig cohomology, and show that our results are consistent with these (see Section \ref{sec:previouswork}). Other evidence includes the fact that, with the ordering on the simple modules in an $\ell$-block induced by the $\pi$-function, the unipotent part of the decomposition matrix should be lower triangular; the author has tested some of the groups where this decomposition matrix is known, and found that the decomposition matrix is lower triangular in all cases. (This includes making certain new conjectures about currently unknown parameters in the decomposition matrices of $F_4(q)$: see Section \ref{sec:F4}).

\bigskip

We make the following notation: if $M$ is a module, write $\Proj{M}$ for its projective cover and $\Omega(M)$ for the kernel of the natural map $\Proj M\to M$. Similarly, we write $\Omega^{-1}(M)$ for the kernel of the natural map from $M$ to its injective hull. If $M$ is a module with $i$th radical layer $R_i$, we write $M=R_1/R_2/\cdots/R_d$. (Of course, in general this does not determine a module up to isomorphism, but does in the cases we use it.)

\section{Deligne--Lusztig Varieties}
\label{sec:DLvars}

In this section we give information on the groups and varieties that we deal with in this paper. After this section, since we are only proving properties of the function $B_d(-)$ for various groups of Lie type, we will have no need for the Deligne--Lusztig varieties, and this will be the only section that will concern them.

Let $p$ be a prime, and let $\b G$ be a connected, reductive algebraic group over the field $\bar\F_p$. Let $F$ be an endomorphism of $\b G$, with $F^\delta$ a Frobenius map for some $\delta\geq 1$ relative to an $\F_{q^\delta}$-structure on $\b G$. Let $W$ denote the Weyl group of $\b G$, $B^+$ the braid monoid of $W$, and let $\phi$ denote the automorphism of $W$ (and hence $B^+$) induced by $F$. Let $\ell\neq p$ be a prime number dividing $G=\b G^F$, and let $P$ be a Sylow $\ell$-subgroup of $G$. Finally, we let $\mc O$, $K$ and $k$ be, as usual, a complete discrete valuation ring, its field of fractions, and its residue field; we assume that $\mc O$ is an extension of the $p$-adic integers $\Z_\ell$, so that $K$ is an extension of $\Q_\ell$ and $k$ is an extension of $\F_\ell$; we assume, again as usual, that these extensions are sufficiently large, for example the algebraic closures. (The assumption that $\Q_\ell\subset K$ makes it easier for the theory of Deligne--Lusztig varieties.) 
We first assume that $\ell$ does not divide the order of $|W\gen\phi|$; in this case $P$ is abelian \cite{brouemalle1992}. 

Let $B$ be a unipotent $\ell$-block of $G$, which has abelian defect group $D$ since $P$ is abelian. By Brou\'e's abelian defect group conjecture $B$ should be derived equivalent to its Brauer correspondent $b$. However, in the case of groups of Lie type one expects that this derived equivalence can be chosen to have a geometric origin. More precisely, it is expected that there is a Deligne--Lusztig variety $Y$ associated to $G$, whose cohomology initially carries an action of $G$ on the one side and an action of a $\Phi_d$-torus containing $D$ on the other, and which can be extended to an action of $\Norm_G(D)$ inducing a derived equivalence between $B$ and $b$.

More structure can be placed on this expected derived equivalence: it should be \emph{perverse} (see Section \ref{sec:pervequiv}), and here we conjecture that the associated perversity function should be the function $\pi(-)$ given in the introduction. This statement has two consequences: the first is that there exists a perverse equivalence with certain properties, and the second is that the cohomology of $Y$ has particular properties.

\begin{conj}\label{conj:perverseequiv} If $\chi_1,\dots,\chi_s$ are the unipotent ordinary characters of $B$, yielding the non-negative integers $\pi(\chi_1),\dots,\pi(\chi_s)$, then there is a bijection between the simple $B$-modules and the simple $b$-modules yielding a perverse equivalence between $\MOD B$ and $\MOD b$ with $\pi(\chi_i)$ as perversity function.
\end{conj}

(The nature of the bijection itself is delayed until a later paper, and is related to the eigenvalues of the Frobenius on the $\chi_i$.) Examining the cohomology of $Y$ yields the next conjecture.

\begin{conj}\label{conj:cohomfound} In the cohomology $H^\bullet(Y,K)$, the character $\chi_i$ appears in degree $\pi(\chi_i)$ and no other degree.
\end{conj}

Both of these conjectures are related to the most general form of Brou\'e's conjecture, the so-called \emph{geometric form}, which is that the complex of $Y$ over $\mc O$ induces a perverse equivalence between $B$ and $b$.

\medskip

We now describe the cases in which $Y$ has been identified. If $B$ is a unipotent $\ell$-block of $kG$ then associated to $B$ is a $d$-cuspidal pair $(L,\psi)$, where $L$ is a $d$-split Levi subgroup and $\psi$ is a $d$-cuspidal unipotent character of $L$. If $L$ is a torus then the variety $Y$ was identified in \cite{brouemichel1997}, and we briefly describe this case (see also \cite[\S3.4]{cravenrouquier2010un}).

Let $w\mapsto \b w$ be the length-preserving lift $W\to B^+$ of the canonical map $B^+\to W$, and let $\b w_0$ be the lift of the longest element of $W$ in $B^+$. Choose $b_d\in B^+$ such that $(b_d\phi)^d=(\b w_0)^2\phi^d$; the variety $Y$ should be the Deligne--Lusztig variety $Y(b_d)$.

Recently \cite{dignemichel2011un} a generalization of this construction of $Y(b_d)$ was given, producing so-called \emph{parabolic Deligne--Lusztig varieties}. The construction of these is considerably more complicated, and we do not attempt it here. In \cite{dignemichel2011un} a candidate variety $Y$ is identified in the case where $L$ is minimal (i.e., the trivial character of $L$ is $d$-cuspidal). Thus in these cases the variety $Y$ has been found, but in general the case remains open.

\medskip

We now relax the condition that $\ell$ does not divide $|W\gen\phi|$. The geometric version of Brou\'e's conjecture no longer applies, but one may still search for a perverse equivalence with perversity function given above, and indeed in \cite{cravenrouquier2010un} Rouquier and the author constructed perverse equivalences with this perversity function for, among others, $\PSL_5(q)$ and $\ell=3$, whenever the Sylow $\ell$-subgroup has order $9$. One may expect Conjecture \ref{conj:perverseequiv} to extend to all cases where $\ell$ is good.

\section{Integrality of $\pi$ and a Bijection with Signs}
\label{sec:integrality}
In this section we prove that the $\pi$-function is always an integer, and demonstrate that, in a bijection with signs $\Irr(B)\to \Irr(b)$, that the sign attached to $\chi$ is $(-1)^{\pi(\chi)}$. Throughout this section, let $d$ be a positive integer at least $2$, and let $\zeta=\e^{2\pi\I/d}$, a primitive $d$th root of unity. (The case where $d=1$ is easy, and omitted.) Let $\mc F$ denote the set of all polynomials in $q$ that are products of cyclotomic polynomials and powers of $q$. 

We begin by determining $B_d(q^r-1)$, and then show that the argument of $\zeta^r-1$, as a complex number, is $(\pi/d)\cdot B_d(q^r-1)$. This implies that, for any $f(q)\in \mc F$, the argument of $f(\zeta)$ is $(\pi/d)\cdot B_d(f)$. We then apply this to groups of Lie type.


\begin{lem}\label{lem:Bdonqrminus1} For $r$ an integer at least $1$,
\[ B_d(q^r-1)=r+d\left\lfloor \frac rd \right\rfloor+\frac d2.\]
The function $B_d(-)$ is the unique homomorphism $f$ from $\mc F$ to $\Q$ that satisfies this property with $f(q)=2$.
\end{lem}
\begin{pf}Since $q^r-1$ is the product of $\Phi_i(q)$ for $i\mid r$, and $B_d(\Phi_1(q))=1+d/2$ and $B_d(\Phi_r(q))=\phi(r)+d\phi_d(r)$, it suffices to show that
\[ \sum_{i\mid r} \phi_d(r)=\left\lfloor\frac{r}{d}\right\rfloor.\]
The proof of this is a simple generalization of the proof that the sum of $\phi(i)$ for $i\mid r$ is $r$: the map $n\mapsto (r/\gcd(n,r),n/\gcd(n,r))$ is a bijection between the set of integers $\{1,\dots,\lfloor r/d\rfloor\}$ and the set of pairs of coprime positive integers $(i,j)$, for $i\mid r$ and $1\leq j\leq i/d$. This bijection proves that the equality is correct, and completes the proof.

The second statement is clear by proving it for $\Phi_r(q)$, either by M\"obius inversion or by induction on $r$.
\end{pf}

Having determined $B_d(q^r-1)$, we can now compute the complex argument of $\zeta^r-1$.

\begin{prop}\label{prop:argumentfirstcase} Let $r$ be an integer not divisible by $d$. The argument of the complex number $\zeta^r-1$ is $\pi(r/d+\lfloor r/d\rfloor+1/2)$.
\end{prop}
\begin{pf} Suppose firstly that $r<d$: then $\zeta^r=\e^{2r\pi i/d}$, and it is easy to see that if $z$ is a complex number with argument $\alpha$ and norm $1$ then $z-1$ has argument $(\alpha+\pi)/2$, so that $\zeta^r-1$ has argument $\pi(r/d+1/2)$, as claimed.

If $r>d$, then write $r=ad+b$: we see that $\zeta^r-1$ has argument $\pi(b/d+1/2)$, which is (modulo $2\pi$), $\pi(2a+b/d+1/2)=\pi(r/d+\lfloor r/d\rfloor+1/2)$.
\end{pf}

This proves that, for $\Phi_r(q)$ with $d\nmid r$, $B_d(\Phi_r)\pi/d$ is the argument of $\Phi_r(\zeta)$ (modulo $2\pi$). We must now deal with $\Phi_r(q)$ for $d\mid r$. We will not have to consider $B_d(\Phi_d)$ in what follows.

\begin{prop} Let $r>1$ be an integer. If $f(q)=(q^{rd}-1)/(q^d-1)$, then the argument of $f(\zeta)$ is zero. Consequently, the argument of $\Phi_r(\zeta)$ is $B_d(\Phi_r)\pi/d$, for any $r\neq d$.\end{prop}
\begin{pf} This is clear since
\[ \frac{q^{rd}-1}{q^d-1}=1+q^d+q^{2d}+\cdots+q^{(r-1)d},\]
and if $q^d=1$ then $(q^{rd}-1)/(q^d-1)=r$.
Hence the argument of $f(\zeta)$ is zero, as needed. To see the consequence, $B_d(f)$ is a multiple of $2d$, so that $B_d(f)\pi/d$ is the argument of $f(\zeta)$ modulo $2\pi$; the statement that $B_d(\Phi_r)\pi/d$ is the argument of $\Phi_r(\zeta)$ now follows as before, by M\"obius inversion for example.
\end{pf}

Clearly the argument of $q$, evaluated at $\zeta$, is $2\pi/d=B_d(q)\pi/d$, and hence we have the following general theorem.

\begin{thm}\label{thm:argumentgeneral} Let $f$ be a polynomial in $\mc F$, and suppose that $\Phi_d(q)\nmid f$. Modulo $2\pi$, the argument of $f(\zeta)$ is $B_d(f)\pi/d$.
\end{thm}

Let $\chi$ be a unipotent ordinary character in a block $B$, with associated $d$-cuspidal pair $(L,\psi)$. It is known \cite[\S5]{bmm1993} that $\psi(1)$ divides $\chi(1)$ (as polynomials in $q$) and $\chi(1)/\psi(1)\equiv (-1)^\ep\alpha\bmod \Phi_d(q)$ (again, as polynomials), for some positive $\alpha\in \Q$ and $\ep\in\Z$. Hence $\chi(1)/\psi(1)$ is a polynomial which, when $q$ is a primitive $d$th root of unity, evaluates to $\pm\alpha$, a real number. Thus $(B_d(\chi)-B_d(\psi))/d$, which modulo $2$ is the argument of $\pm\alpha$ divided by $\pi$, must be $\ep$ modulo $2$; in particular, $\pi(\chi)$ is always an integer, proving Theorem \ref{thm:pifnisinteger}.

If $\ell$ is large (i.e., does not divide the order of the Weyl group in the split case, in general see \cite{bmm1993}) then it is also proved in \cite[\S5]{bmm1993} that $(-1)^\ep=(-1)^{\pi(\chi)}$ is the sign in a perfect isometry between $B$ and $b$, so this proves Theorem \ref{thm:bijectionwithsigns} as well.


\section{Perverse Equivalences}
\label{sec:pervequiv}
In this section we briefly recap the theory of perverse equivalences, at least those parts that affect our results here. We begin with an effective definition of a perverse equivalence, a special type of derived equivalence.

\subsection{Definition and Algorithm}
\label{sec:defnalg}
\begin{defn} Let $A$ and $B$ be $R$-algebras, and let $f:D^b(\MOD A)\to D^b(\MOD B)$ be a derived equivalence. Then $f$ is \emph{perverse} if there exist orderings $S_1,\dots,S_r$ and $T_1,\dots,T_r$ of the simple $A$- and $B$-modules, and a function $\pi:\{1,\dots,r\}\to \Z_{\geq 0}$, such that, in the cohomology of $f(S_i)$, the only composition factors of $H^{-j}(f(S_i))$ are $T_\alpha$ for those $\alpha$ such that $\pi(\alpha)<j$, and a single copy of $T_i$ in $H^{-\pi(i)}(f(S_i))$.
\end{defn}

Often it is assumed that the ordering on the simples for $A$ is such that the $\pi$-function is (weakly) increasing, but if one removes this requirement it means that we can have two different orderings simultaneously; this makes what is going on more transparent for blocks of groups of Lie type with cyclic defect groups, where we really do have orderings on the simple modules for both the group (coming from the Brauer tree) and for the normalizer of the defect group (coming from ``jumps" in the perversity function -- see Section \ref{sec:perverseBrauer}). We also can envisage this as a bijection between the simple modules for the two algebras, which we will have occasion to do.

Along with the concept of a perverse equivalence, and why it is so useful in practice, is an algorithm to compute it. The orderings on the simple modules for $A$, together with the perversity function $\pi$, determine $\MOD B$ up to Morita equivalence, and the algorithm produces the `unique' perverse equivalence with these data. To discuss this algorithm, we let $G$ be a finite group and $B$ be a block of $kG$. Let $D$ denote its defect group, $N=\Norm_G(D)$, and $b$ the Brauer correspondent of $B$. Let $S_1,\dots,S_r$ denote the simple $B$-modules and $T_1,\dots,T_r$ denote the simple $b$-modules. Let $\pi:\{1,\dots,r\}\to\Z_{\geq 0}$ be a perversity function. We describe the image of the simple module $S_i$ as a complex $X_i$ in $D^b(\MOD b)$, describing first the case where induction and restriction is a stable equivalence.

The first term of the complex is the projective cover of $T_i$, denoted $\Proj{T_i}$, in degree $-\pi(i)$. The cohomology $H^{-\pi(i)}(X_i)$ consists of $T_i$ in the socle, and the largest submodule of $\Proj{T_i}/T_i$ consisting of those $T_\alpha$ such that $\pi(\alpha)<\pi(i)$. This module $M_{\pi(i)}$ will be the kernel of the map from degree $-\pi(i)$ to $-\pi(i)+1$; let $L_{\pi(i)}=\Omega^{-1}(M_{\pi(i)})$, i.e., $\Proj{T_i}/M_{\pi(i)}$.

For $0<j<\pi(i)$, the $-j$th term of the complex $X_i$ is the injective hull $P_j$ of $L_{j+1}$; the submodule $L_{j+1}$ is the image of the previous map, and define $M_j$ to be the largest submodule of $P_j$, containing $L_{j+1}$, such that $M_j/L_{j+1}$ has composition factors only those $T_\alpha$ such that $\pi(\alpha)<j$. The module $M_j/L_{j+1}$ is $H^{-j}(X_i)$, and $M_j$ is the kernel of the map from degree $-j$ to degree $-j+1$. Again, write $L_j=\Omega^{-1}(M_j)$.

Finally, the $0$th term of $X_j$ is the module $L_1$, which should be the Green correspondent of $S_i$.

If induction and restriction is not a stable equivalence (or even if it is, but it is not the desired stable equivalence) then we need to modify the Green correspondent to take account of this, by inserting relative $Q$-projective modules in all degrees from $-1$ to $-\eta(Q)$, for the various proper subgroups $Q$ of $D$ and integers $\eta(Q)$. The precise structure of this stable equivalence is not known in general; see \cite{cravenrouquier2010un} for more details, and examples of this for various groups. In the cases here where we find perverse equivalences directly, we will deal with cyclic defect groups, so that induction and restriction will work for us.

An important remark is that, if the injective module $P_j$ in degree $-j$ has a simple module $T_\alpha$ in its socle, then $\pi(\alpha)>j$, since otherwise in degree $j-1$ the module $T_\alpha$, which lies in the socle of $L_{j+1}$, would have been subsumed into $M_{j+1}$.

We finally discuss the cohomology of the complexes $X_i$, and how this may be used to reconstruct the (unipotent part of) the decomposition matrix of the block $B$. Let $\pi$ and the $S_i$ and $T_i$ be as above, and let $X_i$ be the complex in $D^b(\MOD b)$ obtained by running the algorithm. The \emph{alternating sum of the cohomology} $H(S_i)$ of $X_i$ is the virtual $b$-character
\[ \bigoplus_{j=0}^{\pi(S_i)} \;\mathop{\bigoplus_{S\in \mathrm{cf}(H^{-j}(X_i))}} (-1)^{j-\pi(S)}S,\]
where $\mathrm{cf}(M)$ is the set of composition factors of $M$. These virtual $b$-characters determine $r$ rows of the decomposition matrix in an easy way, and can determine the rest of the decomposition matrix if the corresponding rows of $b$ are known. (We will assume that $b$ is the only $\ell$-block of a group $P\rtimes H$, where $P$ is an $\ell$-group and $H$ is an $\ell'$-group in our description, but this condition can be relaxed somewhat, to the statement that all simple $b$-modules come from irreducible representations in characteristic $0$.)

We will explain this description via example.

\begin{example} Let $G=G_2(3)$ and $\ell=13\mid\Phi_3(3)$. Let $P$ denote a (cyclic) Sylow $\ell$-subgroup, $N=\Norm_G(P)=\Z_\ell\rtimes \Z_6$, and order the simple $kN$-modules so that the $i$th radical layer of $\Proj{k}$ is $T_i$ for $1\leq i\leq 6$. The $\pi$-function for $G$ is given in Section \ref{sec:G2}, and the ordering on the simples for the principal block $B$ of $kG$ is $\phi_{1,0}=k$, $G_2[\theta^2]$, $\phi_{2,2}$, $G_2[\theta]$, $\phi_{1,6}$, $G_2[1]$. The $\pi$-function is, with this ordering, $0,3,3,3,4,4$. (It is a coincidence that, in this case, the obvious ordering on the simple $kN$-modules makes the $\pi$-function weakly increasing, and in general this does not happen.)

The Green correspondents of the simple $B$-modules have dimensions $1$, $12$, $11$, $12$, $5$ and $1$, and have radical layers (writing $i$ for $T_i$)
\[ C_1=1,\;\; C_2=6/\cdots/5,\;\; C_3=2/\cdots/6,\;\; C_4=3/\cdots/2,\;\; C_5=5/\cdots/3,\;\; C_6=4.\]
(We can delete the inner radical layers since a $kN$-module is determined by its dimension and socle (or top).)

Running the algorithm with the $\pi$-function above on the simple $kN$-modules, we get six complexes, of the form:
\[\begin{array}{lr}
   X_2:&\;\; \Proj2\to\Proj6\to\Proj6\to C_2.
\\ X_3:&\;\; \Proj3\to\Proj2\to\Proj2\to C_3.
\\ X_4:&\;\; \Proj4\to\Proj3\to\Proj3\to C_4.
\\ X_5:&\;\; \Proj5\to\Proj6\to\Proj4\to\Proj5\to C_5.
\\ X_6:&\;\; \Proj6\to\Proj5\to\Proj5\to\Proj4\to C_6.
\end{array}\]
The cohomology of the complexes above is displayed in the following table.
\smallskip
\begin{center}\begin{tabular}{cccccc}
\hline $X_i$ & $H^{-4}$ & $H^{-3}$ & $H^{-2}$ & $H^{-1}$ & Total
\\ \hline $2$&& $1/2$& $1$& & $2$
\\ $3$& & $3$& & $1$& $3-1$
\\ $4$& & $4$ & & & $4$
\\ $5$& $1/2/3/4/5$ & & & & $5-4-3-2+1$
\\ $6$& $6$ & & & & $6$
\\ \hline \end{tabular}\end{center}

The column `Total' gives the alternating sum in cohomology. To construct the first six rows of the decomposition matrix for $B$, we stipulate that the vector consisting of $0$ everywhere except a $1$ in the $i$th position should be the sum of the rows (with signs) given in the Total column. Hence the third row, minus the first row, should be $(0,0,1,0,0,0)$, and hence the third row is $(1,0,1,0,0,0)$. Continuing this, we get the matrix below.
\[ \begin{tabular}{ccccccccc}
\hline Name & Degree & $\pi$ & $S_1$ & $S_2$ & $S_3$ & $S_4$ & $S_5$ & $S_6$
\\ \hline $\phi_{1,0}$ & $1$ & $0$ & $1$ &&&&&
\\ $G_2[\theta^2]$ & $q\Phi_1^2\Phi_2^2/3$ & $3$&&1&&&&
\\ $\phi_{2,2}$ & $q\Phi_2^2\Phi_6/2$  & $3$&1&&1&&&
\\ $G_2[\theta]$ & $q\Phi_1^2\Phi_2^2/3$  & $3$&&&&1&&
\\ $\phi_{1,6}$ & $q^6$ & $4$ & &1&1&1&1&
\\ $G_2[1]$ & $q\Phi_1^2\Phi_6/6$ & $4$ &&&&&&1
\\ \hline\end{tabular}\]
\end{example}

A crucial remark is that, when the simple $B$-modules are ordered in terms of increasing perversity, the rows of the decomposition matrix are naturally lower triangular, as we see in this example. This yields a bijection between (some of) the irreducible ordinary $B$-characters and the simple $B$-modules; this allows us to transfer the $\pi$-function from the ordinary unipotent characters for unipotent blocks (and the subset of irreducible ordinary characters are the unipotent ones) to the simple modules in characteristic $\ell$. In the case of a Brauer tree this bijection (between the non-exceptional vertices and edges) is easy to describe: if a non-exceptional vertex has degree $1$, associate it to its unique adjacent edge, and remove both, repeating this process until all simple modules are exhausted.

\subsection{Genericity}
\label{sec:genericity}

Let $\ell$ be an integer (not necessarily prime, nor even a prime power), and let $H$ be an $\ell'$-group. Let $\rho:H\to \GL_n(\C)$ be a faithful representation of $H$. It is well known that there exists an algebraic number field $K$, with ring of integers $\mc O$, such that $\rho$ may be written as $\rho:H\to\GL_n(\mc O)$. Suppose that $\ell$ is chosen so that there is a surjective homomorphism $\mc O\to \Z_\ell$ (the ring $\Z/\ell\Z$), inducing the map $\alpha:\GL_n(\mc O)\to \GL_n(\Z_\ell)$ with $\ker\alpha\cap\ker\rho=1$. This yields a map $H\to\Aut(\Z_\ell^n)$ (where here $\Z_\ell$ is considered simply as a group), so we may form the group $G_\ell=(\Z_\ell)^n\rtimes H$; this group is in some sense generic in the integer $\ell$. These groups can be found as the normalizers of $\Phi_d$-tori in groups of Lie type, where $|\Phi_d|=\ell$.

Now specify $\ell$ to be a prime power, and let $k$ be a field whose characteristic divides $\ell$. In the situation of Brou\'e's conjecture, we can let $b$ be the group algebra $kG_\ell$ (since, at least if $\ell$ is a prime power, this group algebra has only one block). The simple $b$-modules are in natural one-to-one correspondence with the simple $kH$-modules, so the simple $b$-modules are `independent' of $\ell$, in the sense that there is a natural identification of the simple $b$-modules for all suitable $\ell$. If $\ell$ and $\ell'$ are two suitable prime powers, we say that the simple $kG_\ell$- and $kG_{\ell'}$-modules are \emph{identified}.

With this identification of the simple $kG_\ell$-modules, we may run the algorithm `generically', without necessarily specifying $\ell$, with a fixed $\pi$-function. In general, the results of the algorithm do depend on $\ell$, but for sufficiently large $\ell$ this should not be the case. While this has not yet been proved (this is ongoing work of Rapha\"el Rouquier and the author) the case where $n=1$ (i.e., the $\ell$-group is cyclic) can easily be proved, as we demonstrate below.

Before we start, we want to extend our definition of identified modules: let $\ell$ be a power of a prime $p$, and suppose that $d\mid(p-1)$. We can construct the group $G_\ell=\Z_\ell\rtimes \Z_d$, and consider the indecomposable $kG_\ell$-modules, where $k$ is an algebraically closed field of characteristic $p$. The group algebra $kG_\ell$ has a single block, with cyclic defect group, and the Brauer tree of $kG_\ell$ is a star, with $d$ vertices on the boundary. The projective cover of any simple module is uniserial: label the simple $kG_\ell$-modules so that $T_1$ is the trivial module, and the first $d$ radical layers of $\Proj{T_i}$ are the simple modules $T_1$, $T_2$, \dots, $T_d$. For any $1\leq i,j\leq d$, there exists a unique uniserial module with $j$ layers and socle $T_i$: write $U_{i,j}$ for this indecomposable module. If $\ell'$ is power of another prime $p'$ with $d\mid(p'-1)$, then we can perform the same construction, and produce uniserial modules $U_{i,j}'$; we identify $U_{i,j}$ and $U_{i,j}'$.

\begin{prop}\label{prop:genericity} Let $H$ be a cyclic group, of order $d$, represented as a complex reflection group. Let $\ell$ and $\ell'$ be powers of primes $p$ and $p'$ such that $d\mid (p-1),(p'-1)$, and write $G_1=G_\ell$ and $G_2=G_{\ell'}$, using the construction above. If $\pi:\{1,\dots,d\}\to\Z_{\geq 0}$ is a perversity function then, if $X_i$ and $X_i'$ ($1\leq i\leq d$) are the complexes describing the results of the algorithm applied to $G_1$ and $G_2$ respectively, we have:
\begin{enumerate}
\item for $1\leq j\leq \pi(i)$, the projective module in degree $-j$ for both $X_i$ and $X_i'$ is the projective cover $\Proj{T_\alpha}$ for some $1\leq \alpha\leq d$ (where we identify the simple modules $T_\alpha$);
\item the module $H^{-j}(X_i)$ is a uniserial module $U_{\alpha,\beta}$, and this is identified with $H^{-j}(X_i')$;
\item writing $A_i$ for the term in degree $0$ of $X_i$, and $A_i'$ for the term in degree $0$ of $X_i'$, if $\pi(i)$ is even then $A_i$ and $A_i'$ are identified uniserial modules, and if $\pi(i)$ is odd then $\Omega(A_i)$ and $\Omega(A_i')$ are identified uniserial modules.
\end{enumerate}
\end{prop}
\begin{pf} Label the uniserial $kG_1$-modules of length at most $d$ (and hence also the $kG_2$-modules via identification) $U_{\alpha,\beta}$, as above. Fix $1\leq i\leq d$, and for $kG_1$ and $1\leq j\leq \pi(i)$, we construct the modules $P_j$, $M_j$ and $L_j$, as in the algorithm, so that $P_j$ is the injective hull of $L_{j+1}$, and $M_j$ is the largest submodule of $P_j$, containing $L_{j-1}$, such that $M_j/L_{j-1}$ contains as composition factors only modules $T_\alpha$ where $\pi(\alpha)<j$. For $kG_2$ we construct the modules $P_j'$, $M_j'$ and $L_j'$ similarly.

We proceed by reverse induction on $j$, starting with the case $j=\pi(i)$. Here, $P_j=\Proj{T_i}$ and $P_j'=\Proj{T_i}$, so (i) of the proposition is true for $j=\pi(i)$. Additionally, $H^{-\pi(i)}(X_i)$ is uniserial of length $r+1$ for some $r\geq 0$, so is the module $U_{i,r+1}$, with radical layers $T_{i-r},T_{i-r+1},\dots,T_i$ (with indices read modulo $d$); this is the largest $r\geq 0$ such that all of $T_{i-r},T_{i-r+1},\dots,T_{i-1}$ have $\pi$-value less than $\pi(i)$. Clearly $r<d$, as the $(d+1)$th socle layer of $\Proj{T_i}$ is $T_i$, which cannot be part of $H^{-\pi(i)}(X_i)$; hence $r$ is independent of the particular exceptionality of the vertex, and so $H^{-\pi(i)}(X_i)$ and $H^{-\pi(i)}(X_i')$ are both $U_{i,r+1}$, proving (ii) for $j=\pi(i)$.

Now let $j$ be less than $\pi(i)$. We notice that, if the top of $H^{-(j+1)}(X_i)$ -- which is the top of $M_{j+1}$ -- is $T_\alpha$ for some $\alpha$, then the projective module in degree $-j$ is $\Proj{T_{\alpha-1}}$; since $T_{\alpha-1}$ was not included in $M_{j+1}$, we must have that $\pi(\alpha-1)\geq j+1$. Since $H^{-(j+1)}(X_i)$ is identified with $H^{-(j+1)}(X_i')$, we see that both $P_j$ and $P_j'$ are $\Proj{T_{\alpha-1}}$, and so (i) is true for $j$. Also, if $P_{j+1}=\Proj{T_\beta}$, then the top of $P_{j+1}$, and hence the top of $L_j$, is $T_\beta$: by the remark just before the start of this section, $\pi(\beta)>j$. 

The module $M_j/L_{j-1}$ is uniserial, with radical layers $T_{\beta-s},T_{\beta-s+1},\dots,T_{\beta-1}$ (with indices read modulo $d$), and some $s$, possibly zero; this is the largest $s\geq 0$ such that all of $T_{\beta-s},T_{\beta-s+1},\dots,T_{\beta-1}$ have $\pi$-value less than $j$. Clearly $s<d$, as the $T_{\beta-d}=T_\beta$, and $\pi(\beta)>j$. Hence $H^{-j}(X_i)=U_{\beta-1,s}$; as the top of $L_j'$ is also $T_\beta$, we must also have that $H^{-j}(X_i')=U_{\beta-1,s}$, proving (ii) for this $j$. Hence, by reverse induction, (i) and (ii) hold.

It remains to deal with (iii). We note that $A_i=\Omega^{-1}(M_1)$ and $A_i'=\Omega^{-1}(M_1')$; since all projective modules have dimension $\ell$ and $\ell'$ respectively, $\dim(A_i)+\dim(M_1)=\ell$, and $\dim(A_i')+\dim(M_1')=\ell'$. Since the tops of $M_1$ and $M_1'$ are identified simple modules, the socles of $A_i$ and $A_i'$ are identified simple modules; as $A_i$ and $A_i'$ are determined by their dimension and their socle, we need to show that if $\pi(i)$ is even then $\dim A_i=\dim A_i'$, and if $\pi(i)$ is odd then $\dim(\Omega(A_i))=\dim(\Omega(A_i'))$, or equivalently $\dim(M_1)=\dim(M_1')$.

Firstly, $\dim(L_j)+\dim(M_j)=\ell$, and $\dim(M_j)=\dim(L_{j+1})+\dim(H^{-j}(X_i))$; by repeating this calculation, we see that if $\pi(i)-j$ is even, we have
\[ \dim(M_j)=\sum_{\alpha=j}^{\pi(i)} (-1)^{\alpha-j}\dim(H^{-\alpha}(X_i)).\]
If $\pi(i)-1$ is even, so $\pi(i)$ is odd, then $\dim(M_1)=\dim(M_1')$, as the cohomology of $X_i$ and $X_i'$ is the same, yielding (iii) in this case. If $\pi(i)$ is even,
\[ \dim(A_i)=\sum_{\alpha=1}^{\pi(i)} (-1)^{\alpha-j}\dim(H^{-\alpha}(X_i)),\]
and so we get $\dim(A_i)=\dim(A_i')$, as needed for (iii).
\end{pf}

\section{Perverse Equivalences and Brauer Trees}
\label{sec:perverseBrauer}
In this section, $B$ is a block of a finite group with cyclic defect group and $b$ is its Brauer correspondent. The simple $B$-modules are labelled $S_1,\dots,S_d$ with some ordering to be given, and the simple $b$-modules are labelled $T_1,\dots,T_d$, again with some ordering to be determined later.

In \cite{rickard1989}, Rickard proved (although not couched in these terms) that there is always a perverse equivalence for blocks with cyclic defect groups; the proof that this equivalence is perverse is in \cite{chuangrouquierun}. While the perversity function itself is easy to describe, the bijection between the simple modules for the block and its Brauer correspondent is less easy to describe, and we omit it here. (It relates to, but is not the same as, Green's walk on the Brauer tree \cite{green1974}. This ordering will be described explicitly in a later paper, but it is not of importance here.)

\begin{thm}\label{thm:standardperversecyclic} Let $B$ be a block of $kG$ with a cyclic defect group $D$, and let $b$ be its Brauer correspondent in $k\Norm_G(D)$. For $S$ a simple $B$-module, let $f(S)$ denote the length of the path from the exceptional vertex of the Brauer tree of $B$ to the vertex incident to $S$ that is closest to the exceptional vertex; let $r$ be the maximum of the $f(S)$. Depending on the perfect isometry between $B$ and $b$, either $\pi_0(S)=r-f(S)$ or $\pi_0(S)=r-f(S)+1$ is the perversity function for a perverse equivalence between $B$ and $b$, with some ordering on the simple $B$- and $b$-modules.
%
\end{thm}

We will describe a family of perverse equivalences for blocks with cyclic defect groups: by varying the perversity function in a natural way, we get infinitely many different perverse equivalences, for \emph{some} bijection between the simple modules. The perversity function given in Theorem \ref{thm:standardperversecyclic} will be referred to as the \emph{canonical} perversity function, and the ordering on the simple $B$-modules alluded to in this theorem will be referred to as the \emph{canonical ordering}. For the simple $b$-modules, we choose some module to be $T_1$ (the trivial module if $b$ is a principal block), and the canonical ordering is the ordering where $T_i$ is the $i$th radical layer of the projective cover of $T_1$, for all $1\leq i\leq d$. Therefore, if the exceptionality of the vertex of the Brauer tree is $1$, then the projective cover of $T_1$ has radical layers
\[ 1/2/3/4/\cdots/d/1.\] An example of this canonical perversity function is the case $d=6$ in Section \ref{sec:G2}.

We can extend the perversity function given in Theorem \ref{thm:standardperversecyclic} to an arbitrary Brauer tree algebra, as since we will be proceeding by induction on the number of vertices, we need to consider Brauer trees that do not necessarily come from groups.

\medskip

Let $B$ be a Brauer tree algebra and let $b$ be the star with the same number $d$ of non-exceptional vertices. Order the simple $b$-modules $T_1,\dots,T_d$ as above, and let $S_1,\dots,S_d$ and $S_1',\dots,S_d'$ be two orderings on the simple $B$-modules. Let $\pi(-)$ and $\pi'(-)$ be perversity functions on the simple $B$-modules. We say that the pairs $(\pi,\{S_i\})$ and $(\pi',\{S_i'\})$ are \emph{algorithmically equivalent} if, for all $1\leq i\leq d$, if $S_i\cong S_j'$, then when one applies the algorithm to yield complexes $X_i$ and $X_j'$,
\begin{enumerate}
\item The terms in degree 0 of $X_i$ and $X_j'$ are isomorphic as $b$-modules, and
\item The alternating sum of cohomologies of $X_i$ and $X_j'$ are identical as virtual $b$-characters.
\end{enumerate}

The theorem we wish to prove is the following.

\begin{thm}\label{thm:algorithmicallyequivalent} Let $B$ be a Brauer tree algebra, and let $b$ denote the star, with the same number $d$ of vertices as $B$. Let $\pi_0$ denote the canonical perversity function on the simple $B$-modules $S_1,\dots,S_d$ with some ordering, and let $\pi'(-)$ be a $\Z_{\geq 0}$-valued function on the set $\{S_1,\dots,S_d\}$, such that:
\begin{enumerate}
\item if $S_i$ and $S_j$ share a non-exceptional vertex in the Brauer tree of $B$, with $S_j$ closer to the exceptional vertex than $S_i$, then $\pi'(S_j)-\pi'(S_i)$ is positive;
\item $\pi'(S_i)-\pi_0(S_i)$ is a non-negative, even integer.
\end{enumerate}
There is an ordering $S_1',S_2',\dots,S_d'$ of the simple $B$-modules such that the pairs $(\pi_0,\{S_i\})$ and $(\pi',\{S_i'\})$ are algorithmically equivalent.
\end{thm}
\begin{pf} We proceed by induction on the number $d$ of simple $B$-modules, the case where $d=1$ being trivial. For each $d$, we also proceed by induction on the sum of $\pi_0(S_i)$ for all simple $B$-modules.

Suppose firstly that $\pi'(S_i)\geq 2$ for all $1\leq i\leq d$, and let $\pi''(S_i)=\pi'(S_i)-2$. By induction there is an ordering $S_1'',\dots,S_d''$ such that $(\pi_0,\{S_i\})$ and $(\pi'',\{S_i''\})$ are algorithmically equivalent. Apply the algorithm to the pairs $(\pi'',\{S_i''\})$ and $(\pi',\{S_i''\})$, to yield complexes $X_1'',\dots,X_d''$ and $\bar X_1',\dots,\bar X_d'$. Notice that the alternating sums of the cohomology of $X_i''$ and $\bar X_i'$ are identical, and the effect on the degree $0$ term is to apply $\Omega^{-2}(-)$ to it. It is easily seen that $\Omega^{-2}(T_i)=T_{i-1}$ (with indices taken modulo $d$), and so if the $j$th radical layer of the degree 0 term of $X_i''$ is $T_\alpha$ then the $j$th radical layer of the degree 0 term of $\bar X_i'$ is $T_{\alpha-1}$.

Finally, let $S_i'=S_{i-1}''$,  and let $X_1',\dots,X_d'$ be the complexes obtained applying the algorithm to the pair $(\pi',\{S_i'\})$. Note that replacing $S_i''$ with $S_{i-1}''$ is the same as fixing the $S_i''$ and replacing $T_i$ with $T_{i+1}$, so the effect is that the degree 0 term of $X_{i+1}'$ (the complex corresponding to $S_{i+1}'=S_i''$) has the same Green correspondent as the degree 0 term of $X_i''$. Since replacing $T_i$ by $T_{i+1}$ induces a rotation on the star, there is no effect on the alternating sums of the cohomology, and so $(\pi',\{S_i'\})$ and $(\pi'',\{S_i''\})$ are algorithmically equivalent, completing the proof in this case.

\medskip

Hence we may assume that $\pi'(S_a)=\pi_0(S_a)$ is either $0$ or $1$ for some simple $B$-module $S_a$, and choose $S_a$ so that is lies incident to the boundary of the Brauer tree of $B$. Assume firstly that $\pi'(S_a)=0$.

We may remove the edge corresponding to $S_a$ (and the now-isolated vertex) to get the Brauer tree for some algebra $B'$, and we do the same with $T_a$ to get a Brauer tree algebra $b'$. Notice that the functions $\pi'$ and $\pi_0$ on the remaining edges of $B$ still satisfy the hypotheses of the theorem, so there is an ordering $S_1',\dots,S_{a-1}',S_{a+1}',\dots,S_d'$ on the simple $B'$-modules such that $(\pi_0,\{S_i\}\setminus\{S_a\})$ and $(\pi',\{S_i'\})$ are algorithmically equivalent. This induces an ordering $S_1',\dots,S_d'$ by setting $S_a'=S_a$. We claim that $(\pi_0,\{S_i\})$ and $(\pi',\{S_i'\})$ are algorithmically equivalent.

To see this, since $\pi_0(S_a)=0$, at any stage in the algorithm where $T_a$ could be taken as cohomology it is. This means that whether we are running the algorithm on a simple $S_i\not\cong S_a$ in either the algebra $B$ or $B'$, we get the same projective modules in the complexes, and (ignoring $S_a$) the same alternating sum of the cohomology. In particular, the projective $\Proj{S_\alpha}$ in degree $-1$ is the same, and $S_a$ is therefore the head of the term in degree 0. Thus the heads of the terms in degree 0 of the complexes do not depend on whether the algorithm is run in $B$ or $B'$. Similarly, the socle of the terms in degree 0 of the complexes do not depend on whether the algorithm is run in $B$ or $B'$.

However, if the head and the socle of a simple $b$-module $M$ are fixed, and $\dim M\leq d$, then $M$ is fixed, and similarly if $\dim(\Omega(M))\leq d$. One of these cases holds for the degree 0 terms, depending on whether $\pi_0(S_i)$ is even or odd respectively, as we saw in Proposition \ref{prop:genericity}. This proves that the degree 0 terms in the complexes for the pairs $(\pi_0,\{S_i\})$ and $(\pi',\{S_i'\})$ are the same; in particular, their dimensions are the same. As we know that the contribution to the alternating sum of the cohomology of all modules but $S_a$ is the same in the complexes for $(\pi_0,\{S_i\})$ and $(\pi',\{S_i'\})$, the fact that the dimensions of the degree 0 terms are the same means that the contribution of $S_a$ is also the same. Hence both the degree 0 terms and the alternating sums of the cohomology are the same for $S_i\not\cong S_a$, for both $(\pi_0,\{S_i\})$ and $(\pi',\{S_i'\})$.

Lastly, as $\pi_0(S_a)=\pi'(S_a)=0$, we see that these complexes are identical. This proves that $(\pi_0,\{S_i\})$ and $(\pi',\{S_i'\})$ are algorithmically equivalent, as claimed.

\medskip

We are left with the case where there is no $S_a$ with $\pi_0(S_a)=0$. To deal with the case where $\pi_0(S_a)=1$, we simply note that if $\bar\pi_0(S_i)=\pi_0(S_i)-1$, then this also a valid canonical perversity function, and $\bar\pi_0$ and $\bar\pi'(S_i)=\pi'(S_i)-1$ also satisfy the hypotheses of the theorem, so there is an ordering $S_1',\dots,S_d'$ on the simple $B$-modules such that $(\bar\pi_0,\{S_i\})$ and $(\bar\pi',\{S_i'\})$ are algorithmically equivalent. Notice that the alternating sum of the cohomology of the complexes for $(\bar\pi_0,\{S_i\})$ and $(\pi_0,\{S_i\})$ are identical, and the degree 0 terms of the one are simply $\Omega^{-1}$ applied to the degree 0 terms of the other. Since the same statements hold for $(\bar\pi',\{S_i'\})$ and $(\pi',\{S_i'\})$, we get that $(\pi_0,\{S_i\})$ and $(\pi',\{S_i'\})$ are algorithmically equivalent, as required.
\end{pf}

Notice that it is relatively easy to understand the modification to the ordering $S_1,\dots,S_d$ needed to produce the ordering $S_1',\dots,S_d'$: essentially, if $\pi_0(S_a)=\pi'(S_a)$ then set $S_a'=S_a$, and remove these modules from contention. We then subtract $2$ from $\pi'(-)$ and cycle the remaining $S_i$, and repeat.

As a corollary to Theorem \ref{thm:algorithmicallyequivalent}, we get the following result on blocks with cyclic defect groups.

\begin{cor}\label{cor:allperversecyclic} Let $B$ be a block of $kG$ with a cyclic defect group $D$, and let $b$ be its Brauer correspondent in $k\Norm_G(D)$. Let $\pi_0$ denote the canonical perversity function on the simple $B$-modules $S_1,\dots,S_d$ with some ordering, and let $\pi'(-)$ be a $\Z_{\geq 0}$-valued function on the set $\{S_1,\dots,S_d\}$, such that:
\begin{enumerate}
\item if $S_i$ and $S_j$ share a non-exceptional vertex in the Brauer tree of $B$, with $S_j$ closer to the exceptional vertex than $S_i$, then $\pi'(S_j)-\pi'(S_i)$ is positive;
\item $\pi'(S_i)-\pi_0(S_i)$ is a non-negative, even integer.
\end{enumerate}
There is a perverse equivalence from $B$ to $b$ with $\pi'$ as perversity function.
\end{cor}

Notice that actually, from the fact that $\pi_0(-)$ is the slowest-increasing function that takes the correct signs on the simple modules, (ii) can be relaxed to the statement that $\pi'(S_i)-\pi(S_i)$ is even, and non-negative for simple modules incident to a boundary vertex of the Brauer tree. For classical groups, where the boundary has either one or two non-exceptional characters, this is a considerable saving in effort.

We will show in later sections that the perversity function on blocks with cyclic defect group, for groups of Lie type, do satisfy the hypotheses of this corollary in the cases where the Brauer tree is known.

\section{Relationship to Previous Work}
\label{sec:previouswork}

In this section we will summarize some of the previous work on this problem, and how it interacts with Conjecture \ref{conj:DLcohom}.

In the cases of $d=1$ and $d=2$, there is already a conjecture from \cite{dmr2007}, which states that the degree should be $2\deg(\chi(1))/d$.

\begin{prop}\label{prop:conjd12} If $d=1$ or $d=2$ then for $\chi$ in the principal $\ell$-block, $\pi(\chi)=2\deg(\chi(1))/d$.
\end{prop}
\begin{pf} If $d=1$ and $r>1$ then $\phi_d(r)=\phi(r)$, so that $B_d(\Phi_r)=2\deg(\Phi_r)$; the result follows since $\Phi_1$ cannot divide $\chi(1)$ if $\chi$ lies in the principal $\ell$-block.

If $d=2$ then $B_2(\Phi_1)=2=2\deg(\Phi_1)$, and for $r>2$ we see that $\phi_2(r)=\phi(r)/2$, so that $B_d(\Phi_r)=2\deg(\Phi_r)$ again. Hence $B_d(\chi)=2\deg(\chi(1))$ in both cases, as claimed.
\end{pf}

The other case where much is known about the structure of the Deligne--Lusztig variety is when $d$ is the Coxeter number, which for the groups considered here (not the Ree and Suzuki groups) is simply the largest integer $d$ such that $\Phi_d(q)\mid |G(q)|$. In this case, both the structure of the cohomology of the Deligne--Lusztig variety and the geometric version of Brou\'e's conjecture are known.

\begin{thm}[Lusztig \cite{lusztig1976}] Conjecture \ref{conj:cohomfound} on the cohomology of Deligne--Lusztig varieties holds whenever $d$ is the Coxeter number.
\end{thm}

If $d$ is the Coxeter number then the Sylow $\Phi_d$-subgroups are cyclic, so Rickard's theorem holds and there is a perverse equivalence (see Theorem \ref{thm:standardperversecyclic}). In this case, it is actually seen that the perversity function for $d$ the Coxeter number is the canonical perversity function in Theorem \ref{thm:standardperversecyclic}. It is easy to see that, if $A(\chi)$ is the degree of $\chi(1)$ as a polynomial in $q$, $a(\chi)$ is the power of $q$ in a factorization of $\chi(1)$ and $c(\chi)$ is the power of $(q-1)$ in a factorization of $\chi(1)$, then $\pi(\chi)=(A(\chi)+a(\chi))/d+c(\chi)/2$. (In the Coxeter case $\pi(\chi)=B_d(\chi)/d$.) Hence we need to show that $(A(\chi)+a(\chi))/d+c(\chi)/2=\pi_0(\chi)$. For $\GL_n(q)$ and the exceptional types we do this explicitly, but for the other classical groups, and the exceptional groups, we leave it as a simple exercise. Hence we get the following result.

\begin{thm}[Rickard] Conjecture \ref{conj:DLcohom} holds whenever $d$ is the largest integer such that $\Phi_d(q)$ is a factor of $|G(q)|$ (as polynomials in $q$).
\end{thm}

By work of Olivier Dudas \cite[Theorem B]{dudas2010un}, for $G$ not of type $E_7$, $E_8$ and $^2\!F_4$, it is known that the complex of the Deligne--Lusztig variety, over $\mc O$, does indeed induce a perverse equivalence, and so even the geometric version of Brou'e's conjecture holds in this case (see the remark after Conjecture \ref{conj:cohomfound}).

\medskip

In addition to this result, Dudas and Jean Michel have calculated the cohomology of various Deligne--Lusztig varieties, and the results are consistent with the conjecture here. A non-exhaustive list is the following:
\begin{enumerate}
\item $G=\GU_4(q)$, $d=4$;
\item $G=\GU_6(q)$, $d=6$;
\item $G=E_6(q)$, $d=9$;
\item $G={}^2\!E_6(q)$, $d=12$;
\item $G=E_7(q)$, $d=14$.
\end{enumerate}

\section{The Combinatorial Objects}
\label{sec:combobjs}
In this section we introduce partitions and symbols. Much of this is well known and we summarize it briefly here, both to fix notation and for the reader's convenience.

If $\lambda=(\lambda_1,\lambda_2,\dots,\lambda_s)$ is a partition of $n$ (with, for now, $\lambda_i\neq 0$), the \emph{first-column hook lengths} of $\lambda$ is the set $X=\{x_1,\dots,x_s\}$, where $x_i=\lambda_i+s-i$. It is easy to see that the set of all partitions (including the empty partition) is in bijection with the set of all finite subsets of $\Z_{>0}$, via sending a partition to its set of first-column hook lengths.

A \emph{$\beta$-set} is a finite subset of $\Z_{\geq 0}$. We introduce an equivalence relation on all such sets, generated by $X\sim X'$ if $X'=\{0\}\cup\{x+1\,:\,x\in X\}$. The \emph{rank} of $X$ is the quantity $\sum_{x\in X} x-s(s-1)/2$, where $s=|X|$. Notice that the rank is independent of the representative of the equivalence class of $\beta$-set; indeed, if we take the unique representative $X$ with $0\notin X$, then the rank of $X$ is the size of the partition $\lambda$ whose first-column hook lengths are $X$. We tend to order the elements of a $\beta$-set $X=\{x_1,\dots,x_s\}$ so that $x_i>x_{i+1}$.

If $X=\{x_1,\dots,x_s\}$ is a $\beta$-set, then adding a \emph{$d$-hook} to $X$ involves replacing some $x_i$ by $x_i+d$ (of course, this assumes that $x_i+d$ is not an element of $X$), and similarly removing a $d$-hook involves replacing some $x_i$ by $x_i-d$. The \emph{$d$-core} of $X$ is the $\beta$-set obtained by removing all possible $d$-hooks.

The $\beta$-sets of partitions can be more easily understood on the abacus. If $d$ is a positive integer, the \emph{$d$-abacus} is a diagram consisting of $d$ columns, or \emph{runners}, labelled $0,\dots,d-1$ from left to right. Starting with $0$ at the top of the left-most runner, we place all non-negative integers on the runners of the abacus, first by moving across the runners left to right, then moving down the runners, as below.
\begin{center}\begin{tabular}{|c|c|c|c|c|}
\hline 0&1&2&3&4
\\ 5&6&7&8&9
\end{tabular}\end{center}
Each number occupies a \emph{position} of the abacus. A \emph{row} of the abacus is a set of $d$ positions corresponding to $di+j$ for $j=0,\dots,d-1$, for some $i\geq 0$. If $X$ is a $\beta$-set, it can be represented on the $d$-abacus by placing a \emph{bead} at position $i$ whenever $i\in X$, and a \emph{gap} at position $i$ whenever $i\notin X$. For ease of description, we often replace $X$ with an equivalent $\beta$-set so that $x_1$ lies on the far-right runner, or equivalently $x_1\equiv -1\bmod d$. The \emph{first row} of the abacus is the row containing $x_1$, and subsequent rows are numbered upwards.

The act of adding or removing a $d$-hook is very easy to describe on the abacus: it consists of moving a bead one place on its runner, down or up respectively. The $d$-core of $X$ is obtained by moving all beads on the $d$-abacus as far upwards as possible.

\medskip

A \emph{symbol} is an unordered pair $\lambda=\{X,Y\}$ of subsets of $\Z_{\geq 0}$. We will write $X=\{x_1,\dots,x_s\}$ with $x_i>x_{i+1}$, and $Y=\{y_1,\dots,y_t\}$ with $y_i>y_{i+1}$. We introduce an equivalence relation on the set of symbols, which is generated by the relation that $\{X,Y\}\sim \{X',Y'\}$ if $X'=\{0\}\cup\{x+1\,:\,x\in X\}$ and $Y'=\{0\}\cup\{y+1\,:\,y\in Y\}$. If $X=Y$ then the symbol is \emph{degenerate}, and otherwise is \emph{non-degenerate}.

The \emph{defect} of $\lambda=\{X,Y\}$ is the quantity $|\,|X|-|Y|\,|$, and the \emph{rank} of $\lambda$ is the quantity $\sum_{x\in X} x+\sum_{y\in Y} y-\lfloor(|X|+|Y|-1)^2/4\rfloor$. Notice that equivalent symbols have the same defect and rank.

Let $\lambda=\{X,Y\}$ be a symbol. Adding a \emph{$d$-hook} to $\lambda$ involves adding $d$ to one of the elements of either $X$ or $Y$ to get another symbol $\mu$. Adding a \emph{$d$-cohook} to $\lambda$ involves adding $d$ to one of the elements of $X$ and transferring it to $Y$, or vice versa, to get another symbol $\mu$. By removing all $d$-hooks we get the \emph{$d$-core}, and by removing all $d$-cohooks we get the \emph{$d$-cocore}. Adding a $d$-hook does not change the defect of a symbol, but adding a $d$-cohook adds or subtracts $2$.

(If one envisages a symbol as a pair of $\beta$-sets, adding a $d$-hook is simply adding a $d$-hook on the abacus one of the partitions; a $d$-cohook is less easy to visualize.)

\section{Combinatorics for Linear and Unitary Groups}

In this section we describe the unipotent characters for $\GL_n(q)$ and $\GU_n(q)$ and their distribution into blocks, and calculate the function $B_d(-)$ on various polynomials that appear when calculating the $\pi$-function for these groups. Let $G=\GL_n(q)$ for some $n$ and $q$, let $\ell\mid |G|$ be a prime, and write $d$ for the multiplicative order of $q$ modulo $\ell$. We describe briefly the unipotent characters and blocks of $\GL_n(q)$, as discussed in \cite{fongsrin1982}.

\medskip

The unipotent characters of $\GL_n(q)$ are labelled by partitions $\lambda$ of $n$, or equivalently $\beta$-sets of rank $n$. Let $X=\{x_1,\dots,x_s\}$ (with $x_i>x_{i+1}$) be a $\beta$-set of rank $n$, and let $\lambda$ be its corresponding partition. If $\chi_\lambda$ is the unipotent character of $\GL_n(q)$ corresponding to $\lambda$, then
\begin{equation} \chi_\lambda(1)=\frac{\D\(\prod_{i=1}^n (q^i-1)\)\(\prod_{1\leq i<j\leq s}(q^{x_i}-q^{x_j})\)}{\D \(\vphantom{\prod_{i=1}^s}q^{\binom{s-1}{2}+\binom{s-2}{2}+\cdots}\)\(\prod_{i=1}^s\prod_{j=1}^{x_i}(q^j-1)\)}.\label{eq:lineardegrees}\end{equation}
(Later we will refer to the `first' and `second' terms of the numerator and denominator of this equation: these have the obvious meanings.)

It is easy to see that $\chi_\lambda(1)$ does not depend on the choice of $\beta$-set $X$ representing $\lambda$. Two $\beta$-sets $X$ and $Y$, with partitions $\lambda$ and $\mu$, have the same $d$-core if and only if the corresponding unipotent characters, $\chi_\lambda$ and $\chi_\mu$, lie in the same $\ell$-block of $G$: the $d$-cuspidal pair for that block has character labelled by the $d$-core of $\lambda$.

\medskip

Let $G=\GU_n(q)$ for some $n$ and $q$, let $\ell\mid |G|$ be a prime, and write $d$ and $e$ for the multiplicative orders of $q$ and $-q$ respectively modulo $\ell$; then $e=d$ if $4\mid d$, $e=2d$ if $d$ is odd and $e=d/2$ otherwise. As with the linear groups, we summarize briefly the facts about unipotent characters and blocks that we need, taken from \cite{fongsrin1982}.

The unipotent characters of $\GU_n(q)$ are similar to those of $\GL_n(q)$, in that they are again associated to partitions of $n$. If $\chi_\lambda$ is the unipotent character of $\GL_n(q)$ associated to $\lambda$ and $\phi_\lambda$ is the unipotent character of $\GU_n(q)$ associated to $\lambda$, then the degree of $\phi_\lambda$ is obtained from that of $\chi_\lambda$ by replacing $q$ with $(-q)$ (with possibly a sign change if this makes the character degree negative). In the expansion of $\phi_\lambda(1)$ into powers of $q$ and cyclotomic polynomials, this has the effect of replacing $\Phi_r$ with $\Phi_{2r}$ and vice versa, whenever $r$ is odd.

The structure of the $\ell$-blocks of $G$ is similar as well: these are parametrized by $e$-cores, and two unipotent characters $\phi_\lambda$ and $\phi_\mu$ lie in the same $\ell$-block of $G$ if and only if $\lambda$ and $\mu$ have the same $e$-core: the $d$-cuspidal pair for that block has character labelled by the $e$-core of $\lambda$.

\medskip

To describe the $\pi$-function for Brauer trees of $\GL_n(q)$ and $\GU_n(q)$ in the next section, we need to evaluate $(B_d(\chi_\mu)-B_d(\chi_\lambda))/d$, where $\lambda$ is a $d$-core (or $e$-core for unitary groups) and $\mu$ is obtained from $\lambda$ by adding a single $d$-hook (or $e$-hook). Since $B_d(-)$ is a homomorphism, we consider $B_d(\chi_\mu(1)/\chi_\lambda(1))$, so only need to consider the difference between $\chi_\lambda(1)$ and $\chi_\mu(1)$ in Equation \ref{eq:lineardegrees}. Correspondingly we need to evaluate quantities such as $B_d(q^{x_i+d}-q^{x_j})-B_d(q^{x_i}-q^{x_j})$. Firstly, we describe $B_d(q^i\pm q^j)$.

\begin{prop} Let $i$ and $j$ be integers with $i>j$. We have
\[ B_d(q^i-q^j)=i+j+d\cdot\left\lfloor \frac{i-j}{d}\right\rfloor+\frac{d}{2},\]
and
\[ B_d(q^i+q^j)=i+j+d\cdot\(\left\lfloor\frac{2(i-j)}{d}\right\rfloor-\left\lfloor\frac{i-j}{d}\right\rfloor\).\]
\end{prop}
\begin{pf} The case where $j=0$ is Lemma \ref{lem:Bdonqrminus1}. The general case easily follows since $q^i-q^j=q^j(q^{i-j}-1)$.

For the second equality, we have
\[ q^i+q^j=\frac{q^{2i}-q^{2j}}{q^i-q^j},\]
so that 
\begin{align*} B_d(q^i+q^j)&=\(2i+2j+d\cdot\left\lfloor \frac{2i-2j}{d}\right\rfloor+\frac{d}{2}\)-\(i+j+d\cdot\left\lfloor \frac{i-j}{d}\right\rfloor+\frac{d}{2}\)
\\ &=i+j+d\cdot\(\left\lfloor\frac{2(i-j)}{d}\right\rfloor-\left\lfloor\frac{i-j}{d}\right\rfloor\).\end{align*}
\end{pf}

This yields the following proposition in an obvious way, which deals with the effect on the second term in the numerator for the character degree, when going from $\chi_\lambda(1)$ to $\chi_\mu(1)$, for $\GL_n(q)$. (We also include a case that will be needed for symplectic and orthogonal groups.)

\begin{prop}\label{prop:diffsinBdslinear} Let $i$ and $j$ be integers, and let $d\geq 1$ be an integer. We have that
\begin{align*} B_d(q^{i+d}-q^j)-B_d(q^i-q^j)&=\begin{cases}2d & i-j>0\\ d & -d<i-j<0 \\ 0 & i-j<-d\end{cases}
\\ B_d(q^{i+d}+q^j)-B_d(q^i+q^j)&=\begin{cases} 2d& i-j>-d/2\\d&i-j=-d/2\\0&i-j<-d/2\end{cases}
\end{align*}
\end{prop}

As we said before, when working with unitary groups we need two integers, $d$ and $e$, where $e=d$ if $4\mid d$, $e=2d$ if $d$ is odd, and $e=d/2$ otherwise. We will be evaluating $B_d(f)$, with the terms in $f$ generally depending on $e$: we want the analogue of Proposition \ref{prop:diffsinBdslinear} in this case.

\begin{prop}\label{prop:diffsinBdsunitary}
Suppose that $d=e$. Then
\[ B_d((-q)^{i+e}-(-q)^j)-B_d((-q)^i-(-q)^j)=\begin{cases}2d & i-j>0\text{ or }i-j>-d/2\text{ and odd}\\ d&0>i-j>-d\text{ and even}\\0& i-j<-d\text{ or }i-j<-d/2\text{ and odd}\end{cases}\]
Suppose that $d=e/2$. Then
\[ B_d((-q)^{i+e}-(-q)^j)-B_d((-q)^i-(-q)^j)=\begin{cases}4d & i-j>0\text{ or }i-j>-d/2\text{ and odd}\\ 3d & 0>i-j>-d\text{ and even}\\ 2d& -d/2>i-j>-3d/2\text{ and odd}\\d& -d>i-j>-2d\text{ and even}\\0& i-j<-2d\text{ or }i-j<3d/2\text{ and odd}\end{cases}\]
Suppose that $d=2e$. Then
\[ B_d((-q)^{i+e}-(-q)^j)-B_d((-q)^i-(-q)^j)=\begin{cases}d & i-j>0\text{ or }i-j>-d/2\text{ and odd}\\ 0& i-j<-d/2\text{ or }i-j<0\text{ and even}\end{cases}\]
\end{prop}
\begin{pf} The first two statements are easy, and the last one follows once one notices that, if $e=d/2$ and $a\in \N$, then
\[ \left\lfloor\frac{a+e}{d}\right\rfloor+\left\lfloor\vphantom{\frac{a+e}{d}}\frac{a}{d}\right\rfloor=\left\lfloor\frac{2a}{d}\right\rfloor.\]
\end{pf}

We will need to take products of these polynomials when dealing with the first term of the numerator in $\chi_\lambda(1)$ (and $\phi_\lambda(1)$) and the second term of the denominator of $\chi_\lambda(1)$ (and $\phi_\lambda(1)$).

\begin{prop}\label{prop:prodsofpolys} Let $i$ and $j$ be integers, and let $d\geq 1$ be an integer. Write $e=d$ if $4\mid d$, $e=2d$ if $d$ is odd, and $e=d/2$ otherwise. We have
\[ B_d\(\prod_{i=n+1}^{n+d} (q^i-1)\)=2nd+d^2+\frac{3d}{2},\qquad B_d\(\prod_{i=n+1}^{n+e} ((-q)^i-1)\)=2ne+e(e+1)+\frac{d}{2}.\]
\end{prop}
\begin{pf} In both cases, we proceed by induction on $n$. The inductive step, from $n-1$ to $n$, is clear, since in the first case we replace $(q^n-1)$ with $(q^{n+d}-1)$, which contributes $2d$ by Proposition \ref{prop:diffsinBdslinear}, and in the second case we replace $((-q)^n-1)$ with $((-q)^{n+e}-1)$, which contributes $2e$ by Proposition \ref{prop:diffsinBdsunitary}. Hence we only need to prove the formula for $n=0$; this case is trivial.
\end{pf}

\section{Brauer Trees for Linear and Unitary Groups}

The first Brauer trees for which we prove Theorem \ref{thm:brauertreesperverse} are the linear and unitary groups, using the results of the previous section.
\subsection{$\GL_n(q)$}

Let $n$ be a positive integer, let $q$ be a prime power, let $\ell$ be a prime, and write $d$ for the multiplicative order of $q$ modulo $\ell$. Let $B$ be an $\ell$-block of $G=\GL_{n+d}(q)$ with a cyclic defect group, with $d$-core a partition $\lambda$ of $n$; let $X=\{x_1,\dots,x_s\}$ (with $x_i>x_{i+1}$) be a $\beta$-set corresponding to $\lambda$. We will compute the function $\pi(-)$ for the unipotent characters in $B$. There are $d$ unipotent characters $\chi_\mu$, each with $\lambda$ as $d$-core and $|\mu|-|\lambda|=d$; by choosing $X$ sufficiently large, we have the subset $X'=\{x_{i_1},\dots,x_{i_d}\}$ of $X$ consisting of those $d$ integers such that $x_{i_j}+d\notin X$ (i.e., they represent the possible $d$-hooks that may be added), and order them so that $x_{i_j}>x_{i_{j+1}}$. Notice that if one adds $d$ to $x_{i_j}$, then $j$ is the leg length of the corresponding $d$-hook added to $\lambda$.

Label the unipotent characters $\chi_1,\dots,\chi_d$ in $B$ by $\chi_j$ having partition with $x_{i_j}$ incremented by $d$. By \cite{fongsrin1984}, the Brauer tree of a block $B$, with $d$-core $\lambda$, is a line, with the exceptional vertex at the right end, $\chi_d$ adjacent to it, and $\chi_i$ adjacent to $\chi_{i+1}$, as in the following diagram.
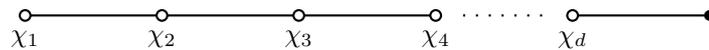
\tikzstyle{every node}=[circle, fill=black!0,
                        inner sep=0pt, minimum width=4pt]
\begin{center}\begin{tikzpicture}[thick,scale=1.8]
\draw \foreach \x in {0,1,2,4}{
(\x,0) -- (\x+1,0)};
\draw[loosely dotted] (3.2,0) -- (3.8,0);
\draw (0,0) node [draw,label=below:$\chi_{1}$] (l0) {};
\draw (1,0) node [draw,label=below:$\chi_{2}$] (l1) {};
\draw (2,0) node [draw,label=below:$\chi_{3}$] (l2) {};
\draw (3,0) node [draw,label=below:$\chi_{4}$] (l3) {};
\draw (4,0) node [draw,label=below:$\chi_{d}$] (l5) {};
\draw (5,0) node [fill=black!100] (ld) {};
\end{tikzpicture}\end{center}

\begin{prop}\label{prop:GLnperversity} With the setup above, we have that
\[ \pi(\chi_j)=2(n-x_{i_j}+s-i_j)+(j-1).\]
\end{prop}
\begin{pf}Let $\mu$ be obtained from $\lambda$ by replacing $x_{i_j}$ by $x_{i_j}+d$, and note that $\pi(\chi_\mu)=(B_d(\chi_\mu)-B_d(\chi_\lambda))/d$. To calculate $B_d(\chi_\mu)-B_d(\chi_\lambda)$, it suffices to evaluate the function $B_d(-)$ on the difference between the formulae for $\chi_\mu(1)$ and $\chi_\lambda(1)$, using Equation \ref{eq:lineardegrees}. This formula has two terms in both the numerator and denominator: the function $B_d(-)$, applied to the difference for the first term of the numerator is $2dn+d^2+3d/2$ by Proposition \ref{prop:prodsofpolys}, and similarly the difference for the second term of the denominator is $2dx_{i_j}+d^2+3d/2$; the difference for the rest of the denominator is zero, since $s$ does not change. Hence so far we have a contribution to $\pi(\chi_\mu)$ of
\[ \frac{(2dn+d^2+3d/2)-(2dx_{i_j}+d^2+3d/2)}{d}=2(n-x_{i_j}).\]
We must now evaluate the difference in the second term of the numerator, which consists of adding terms of the form $(q^{x_{i_j}+d}-q^{x_\alpha})$ and removing terms of the form $(q^{x_{i_j}}-q^{x_\alpha})$, of course using Proposition \ref{prop:diffsinBdslinear}.

The set $X\setminus\{x_{i_j}\}$ is split into three subsets: $X_1=\{x\in X\,:\,x<x_{i_j}\}$, $X_2=\{x\in X\,:\,x_{i_j}<x<x_{i_j}+d\}$, and $X_3=\{x\in X\,:\,x_{i_j}+d<x\}$. Using Proposition \ref{prop:diffsinBdslinear}, we see that for $x\in X_k$, $B_d(q^{x_{j_i}+d}-q^x)-B_d(q^{x_{j_i}}-q^x)$ is one of $2d$, $d$ or $0$, depending on whether $k=1$, $k=2$ or $k=3$. Hence, when computing $B_d(\chi_i)-B_d(\chi_\lambda)=d\pi(\chi_i)$, the contribution from the second term of the numerator is $2d\cdot |X_1|+d\cdot |X_2|$. There are $s-i_j$ elements in $X_1$, and $j-1$ elements in $X_2$, so the contribution from this term is $2d(s-i_j)+d(j-1)$. Dividing by $d$ and adding to the previous contributions gives the claimed formula.
\end{pf}

Notice that $\pi(\chi_1)$ is always even, and $\pi(\chi_j)<\pi(\chi_{j+1})$; since $\pi(\chi_1)\geq 0=\pi_0(\chi_1)$, this proves that there is a perverse equivalence with this as perversity function, by Corollary \ref{cor:allperversecyclic} (and the remark thereafter).

\medskip

We finally note that, for principal blocks, the $\pi$-function is particularly easy to determine, and if $B$ is the principal $\ell$-block of $\GL_{d+r}(q)$ ($0\leq r<d$) then the $\pi$-function is as below.
\tikzstyle{every node}=[circle, fill=black!0,
                        inner sep=0pt, minimum width=4pt]
\begin{center}\begin{tikzpicture}[thick,scale=1.8]
\draw (0,0.15) node{$0$};
\draw (1,0.15) node{$2r+1$};
\draw (2,0.15) node{$2r+2$};
\draw (3,0.15) node{$2r+3$};
\draw (4,0.15) node{$2r+d-1$};
\draw \foreach \x in {0,1,2,4}{
(\x,0) -- (\x+1,0)};
\draw[loosely dotted] (3.2,0) -- (3.8,0);
\draw (0,0) node [draw,label=below:$\chi_1$] (l0) {};
\draw (1,0) node [draw,label=below:$\chi_2$] (l1) {};
\draw (2,0) node [draw,label=below:$\chi_3$] (l2) {};
\draw (3,0) node [draw,label=below:$\chi_4$] (l3) {};
\draw (4,0) node [draw,label=below:$\chi_d$] (l5) {};
\draw (5,0) node [fill=black!100] (ld) {};
\end{tikzpicture}\end{center}

\subsection{$\GU_n(q)$}

We now prove Theorem \ref{thm:brauertreesperverse} for the unitary groups $\GU_n(q)$. Let $n$ be a positive integer, let $q$ be a prime power, let $\ell\mid |G|$ be a prime, and write $d$ and $e$ for the multiplicative orders of $q$ and $-q$ respectively modulo $\ell$; then $e=d$ if $4\mid d$, $e=2d$ if $d$ is odd and $e=d/2$ otherwise. Let $G=\GU_{n+e}(q)$, and let $B$ be an $\ell$-block of $G$ with cyclic defect group.

We use the description of the Brauer trees from \cite{fongsrin1990}. Let $\lambda$ be an $e$-core of size $n$ and let $X$ be a $\beta$-set corresponding to $\lambda$. Let $X'$ denote the subset of $X$ consisting of all $x\in X$ such that $x+e\notin X$, as in the case of $\GL_n(q)$. By replacing $X$ with an equivalent $\beta$-set, we have $|X'|=e$. Divide $X'$ into $Y$ and $Z$, where $Y$ consists of all even elements of $X'$, and $Z$ consists of all odd elements of $X'$, with the ordering on $Y=\{y_1,\dots,y_a\}$ and $Z=\{z_1,\dots,z_b\}$ given by $y_i>y_{i+1}$ and $z_i>z_{i+1}$, as with $X$. Let $\sigma_i$ be the character of $\GU_{n+e}(q)$ obtained by replacing $y_i$ with $y_i+e$, and similarly let $\tau_i$ be the character obtained by replacing $z_i$ with $z_i+e$. The Brauer tree is as follows.

\tikzstyle{every node}=[circle, fill=black!0,
                        inner sep=0pt, minimum width=4pt]
\begin{center}\begin{tikzpicture}[thick,scale=1.8]
\draw \foreach \x in {0,1,3,4,6,7}{
(\x,0) -- (\x+1,0)};
\draw[loosely dotted] (2.2,0) -- (2.8,0);
\draw[loosely dotted] (5.2,0) -- (5.8,0);
\draw (0,0) node [draw,label=below:$\sigma_{1}$] (l0) {};
\draw (1,0) node [draw,label=below:$\sigma_{2}$] (l1) {};
\draw (2,0) node [draw,label=below:$\sigma_{3}$] (l2) {};
\draw (3,0) node [draw,label=below:$\sigma_{a}$] (l5) {};
\draw (4,0) node [fill=black!100] (ld) {};
\draw (8,0) node [draw,label=below:$\tau_{1}$] (l0) {};
\draw (7,0) node [draw,label=below:$\tau_{2}$] (l1) {};
\draw (6,0) node [draw,label=below:$\tau_{3}$] (l2) {};
\draw (5,0) node [draw,label=below:$\tau_{b}$] (l5) {};
\end{tikzpicture}\end{center}
If $e$ is even then the two branches of the tree have the same length, and so it is obvious that the $\pi$-function (since it is non-negative and has the correct parity) satisfies the second requirement of Corollary \ref{cor:allperversecyclic}, so we only need to check the first condition. (When $e$ is odd, both conditions need to be checked.)

\medskip

\begin{prop} Let $T$ denote the Brauer tree of a unipotent block of a group $\GU_{n+e}(q)$. If $\chi$ and $\psi$ are adjacent non-exceptional vertices on $T$ such that $\chi$ is closer to the exceptional node than $\psi$, then $\pi(\chi)>\pi(\psi)$.
\end{prop}
\begin{pf}
We need to prove that $\pi(\sigma_{i+1})-\pi(\sigma_i)$ is positive, and similarly for $\pi(\tau_{i+1})-\pi(\tau_i)$. Since it is odd by Theorem \ref{thm:bijectionwithsigns}, we actually only need to show that it is non-negative; also, replacing a $\beta$-set for $\lambda$ with one with one more member swaps the two branches, so we only need to prove that $\pi(\sigma_{i+1})-\pi(\sigma_i)$ is non-negative, as this will automatically prove the result for the $\tau_i$.

Consider $y_i$ and $y_{i+1}$: we have that $y_i=x_\alpha$ and $y_{i+1}=x_\beta$ for some $\alpha<\beta$. One quantity that will appear in our analysis is $x_\alpha-x_\beta-\beta+\alpha$, which is equal to the number of gaps between the positions $x_\alpha$ and $x_\beta$.

\medskip

\noindent\textbf{Case 1}: $e=d$. Write $f_1(x_i)=s-i$, $f_2(x_i)$ for the number of odd $x\in X$ such that $0<x-x_i<d/2$, and $f_3(x_i)$ for the number of even $x\in X$ such that $0<x-x_i<d$. We have, by Propositions \ref{prop:diffsinBdsunitary} and \ref{prop:prodsofpolys}, together with the same argument as for $\GL_n(q)$, that
\[ \pi(\sigma_i)=2(n-x_\alpha)+2(f_1(x_\alpha)+f_2(x_\alpha))+f_3(x_\alpha).\]
It is clear since $d$ is even that $f_3(x_\alpha)=i-1$ (recall that $x_\alpha=y_i$), since it measures the number of even beads on a given row of the abacus. Hence we have
\[ \pi(\sigma_i)=2(n-x_\alpha+s-\alpha)+(i-1)+2f_2(x_\alpha),\]
and if $x_\beta=y_{i+1}$ (i.e., $x_\alpha$ and $x_\beta$ label adjacent characters on the Brauer tree), then
\[ \pi(\sigma_{i+1})-\pi(\sigma_i)=2(x_\alpha-x_\beta-\beta+\alpha)+1+2(f_2(x_\beta)-f_2(x_\alpha)).\]

If $x_\alpha-x_\beta>d/2$, then $f_2(x_\beta)-f_2(x_\alpha)$ measures the difference between the number of odd beads in $(x_\beta,x_\beta+d/2)$ and $(x_\alpha,x_\alpha+d/2)$, and if $x_\alpha-x_\beta<d/2$ then these intervals overlap, so that $f_2(x_\beta)-f_2(x_\alpha)$ measures the difference between the number of odd beads in $(x_\beta,x_\alpha)$ and $(x_\beta+d/2,x_\alpha+d/2)$. In both cases, if we include the term $x_\alpha-x_\beta-\beta+\alpha$, which counts the gaps between $x_\beta$ and $x_\alpha$, we see that this term is non-negative. This proves the result for the case where $d=e$.

%

\medskip

\noindent\textbf{Case 2}: $e=d/2$. We keep $f_1(x_i)$ and $f_2(x_i)$ from the previous argument, and we have that
\[ \pi(\sigma_i)=n-x_\alpha+f_1(x_\alpha)+f_2(x_\alpha),\]
so that
\[ \pi(\sigma_{i+1})-\pi(\sigma_i)=x_\alpha-x_\beta-\beta+\alpha+f_2(x_\beta)-f_2(x_\alpha).\]
This is almost identical to the previous formula, except that 1 is subtracted and it is halved. Since the previous expression was positive, this one must be non-negative, and this is all that is required, so the proposition holds in this case as well.

\medskip

\noindent\textbf{Case 3}: $e=2d$, and here there are more contributions to consider, as we see from Proposition \ref{prop:diffsinBdsunitary}. We use $f_2(x_i)$ and $f_3(x_i)$ from the previous cases, and also introduce $f_4(x_i)$, the number of odd $x\in X$ such that $d/2<x-x_i<3d/2$, and $f_5(x_i)$, the number of even $x\in X$ such that $d<x-x_i<2d$. We see, from Propositions \ref{prop:diffsinBdsunitary} and \ref{prop:prodsofpolys}, that
\begin{align*} \pi(\sigma_i)&=4(n-x_i)+4f_1(x_i)+4f_2(x_i)+3f_3(x_i)+2f_4(x_i)+f_5(x_i)
\\ &=4(n-x_i+s-i)+2f_2(x_i)+2(f_2(x_i)+f_4(x_i))+2f_3(x_i)+(f_3(x_i)+f_5(x_i)).\end{align*}
We have organized these terms in this way because $f_2(x_i)+f_4(x_i)$ counts the number of odd $x\in X$ in the interval $(x_i,x_i+3d/2)$, and similarly for $f_3(x_i)+f_5(x_i)$.
We get a very similar expression for the difference between consecutive $\pi(\sigma_i)$:
\begin{align*} \pi(\sigma_{i+1})-\pi(\sigma_i)=4(&x_\alpha-x_\beta-\beta+\alpha)+4(f_2(x_\beta)-f_2(x_\alpha))+3(f_3(x_\beta)-f_3(x_\alpha))
\\ &+2(f_4(x_\beta)-f_4(x_\alpha))+(f_5(x_\beta)-f_5(x_\alpha)).\end{align*}
Using the grouping of the $f_k(x_i)$ above, we may apply the same argument as for the first case to see that $4(x_\alpha-x_\beta-\beta+\alpha)$ cancels out all of these differences, so that $\pi(\sigma_{i+1})-\pi(\sigma_i)$ is non-negative, as needed.
\end{pf}

%
%

We must now prove the second condition for being a perverse equivalence, namely that this perversity function is always at least the canonical one. This is only necessary in the case where $e$ is odd, since otherwise there are the same number of vertices either side of the exceptional vertex in the Brauer tree, and the result is clear.

We keep the notation for the Brauer tree from above, and assume that $b>a$. Since the canonical perversity function $\pi_0(-)$ on $\tau_1$ is either $0$ or $1$, and $\pi(\tau_1)$ has the same parity as this, $\pi_0(\tau_1)\leq \pi(\tau_1)$. It suffices therefore to check that $\pi(\sigma_1)\geq \pi_0(\sigma_1)$. Notice that, writing $c=b-a$, an odd integer, we have that $\pi_0(\sigma_1)$ is either $c$ or $c+1$, depending on whether $\pi_0(\tau_1)=0$ or $\pi_0(\tau_1)=1$. As the parities of $\pi_0(\sigma_1)$ and $\pi(\sigma_1)$ are the same, we must prove that $\pi(\sigma_1)\geq c$.

We give a lemma which contains most of the details in the proof that we need.

\begin{lem}\label{lem:unitaryhelper} Let $\lambda$ be a partition of $n$, let $e>1$ be an odd integer, and let $X=\{x_1,\dots,x_s\}$ (with $x_i>x_{i+1}$) denote a $\beta$-set of $\lambda$ on the $e$-abacus. Assume that $X$ is large enough so that there is a subset $X'$ of $e$ elements $x$ such that $x+e\notin X$, and write $X'=Y\cup Z$, where the elements of $Y$ are even and those of $Z$ are odd. Suppose that $|\,|Z|-|Y|\,|=m$.
\begin{enumerate}
\item If $X$ cannot be chosen so that $x_1=2e-1$ (i.e., represented on a two-row abacus with the second row having no gaps) then $n\geq e+1\geq m+1$, and if $\alpha\neq 1$ then $n-\lambda_\alpha\geq e$.
\item Suppose that $X$ is chosen so that $x_1=2e-1$. If $|Y|>|Z|$ then the sum of all parts whose corresponding elements of $X$ are even is at least $m-2$, and similarly if $|Z|>|Y|$ then the sum of all parts whose corresponding elements of $X$ are odd is at least $m-2$.
\end{enumerate}
\end{lem}
\begin{pf} If $X$ cannot be represented on a two-row abacus with the second row having no gaps then this means that there is an element $x\in \{1,\dots,x_1\}\setminus X$ such that $x_1-x\geq e+1$. This implies that the largest of the first-column hook lengths of $\lambda$ is at least $e+1$, and so in particular $n\geq e+1$ and $n-\lambda_\alpha\geq e$. Obviously $e\geq m$, and this completes the proof of (i).

For (ii), we have that $x_1=2e-1$ is odd. First assume that $|Z|>|Y|$, so we are counting the contribution from parts corresponding to odd elements of $X$. Partition the numbers $e,\dots,2e-1$ into subsets $A_1,\dots,A_v$ and $B_1,\dots,B_w$, where the $A_i$ consist of consecutive gaps in the abacus of $X$ and the $B_i$ consist of consecutive beads in the abacus of $X$. (Perform this partitioning so that $v$ and $w$ are minimized.) Write $A_i'=\{a-e:a\in A_i\}$. Notice that $X'$ is the union of the $B_i$ and the $A_i'$. It is easy to see that $w-v$ is either $0$ or $1$, depending on whether $e$ is a gap or a bead respectively.

Each $A_i'$ or $B_i$ can contribute at most one more odd element to $X'$ than even; each of the $B_i$ also contains at least one odd element of $X'$ unless $B_i=\{x\}$ for some even $x$, in which case this $B_i$ actually contributes an extra even, not an extra odd. Suppose that there are $u$ such $B_i$, so that there are at least $w-u$ different odd elements of $X'$ in sets $B_i$: by counting the $A_i'$ and $B_i$ that contribute an extra odd, taking away those that contribute an extra even, we have that $(w-u)-u+v\geq m$. If $w=v$ then $w-u\geq m/2$, but since $m$ is odd we actually have $w-u\geq (m+1)/2$. If $w=v+1$ then we again get $w-u\geq (m+1)/2$. If $w=v$, then the $\beta$-set elements in each $B_i$ correspond to a particular \emph{non-zero} size of part of $\lambda$ (and different $B_i$ correspond to different sizes of part), whereas if $w=v+1$ then $B_1$ consists $\beta$-set elements corresponding to parts of $\lambda$ of size $0$. Thus,
\begin{itemize}
\item if $w=v$, the parts corresponding to odd $\beta$-numbers in $X$ have size at least $1+\cdots+(m+1)/2=(m+1)(m+3)/8$, and
\item if $w=v+1$, the parts corresponding to odd $\beta$-numbers in $X$ have size at least $0+1+\cdots+(m-1)/2=(m^2-1)/8$.
\end{itemize}
In either case, we get that the contributions from parts whose corresponding elements of $X$ are odd is at least $(m^2-1)/8\geq m-2$, since $m$ is a positive odd integer.

If $|Y|>|Z|$ then an identical argument works with `odd' and `even' swapped, and this yields the result.
%
%
%
%
%
%
\end{pf}

Using this lemma we will prove that $\pi(\sigma_1)\geq \pi_0(\sigma_1)$.

\begin{prop} Let $B$ denote a unipotent block of a group $\GU_{n+e}(q)$ with cyclic defect group, and let $\pi(-)$ and $\pi_0(-)$ be as defined above. For any non-exceptional character $\chi$ in $B$, we have that $\pi(\chi)\geq \pi_0(\chi)$.
\end{prop}
\begin{pf} As we have mentioned, it suffices to check the case where $e$ is odd. As before, let $X=\{x_1,\dots,x_s\}$ be a $\beta$-set associated to $\lambda=(\lambda_1,\dots,\lambda_s)$ (where some of the $\lambda_i$ may be $0$), let $X'$, $Y$ and $Z$ be as above, assume that $|Y|<|Z|$, and let $y_1=x_\alpha$; it suffices from the arguments above to prove that $\pi(\sigma_1)\geq \pi_0(\sigma_1)$. We choose $X$ so that $x_1$ lies on the far-right runner of the abacus for ease of explanation. As above, we have
\[ \pi(\sigma_1)=n-x_\alpha+(s-\alpha)+f_2(x_\alpha)=n-\lambda_\alpha+f_2(x_\alpha).\]
There are two cases to consider: when $\alpha=1$ and when $\alpha\neq 1$.

\medskip

\noindent\textbf{Case 1}: $\alpha=1$. In this case, let $\bar\lambda=(\lambda_2,\lambda_3,\dots,\lambda_s)$. Notice that the effect of removing $\lambda_1$ from $\lambda$ is to add an extra element to $Z$ and remove one from $Y$, thus incrementing $c$ by $2$: thus by Lemma \ref{lem:unitaryhelper} we have that $n-\lambda_1=|\bar\lambda|\geq c$, so that $\pi(\sigma_1)\geq c$, as needed.
%
%
%
%
%
%

\medskip

\noindent\textbf{Case 2}: $\alpha\neq 1$. In this case $x_1=z_1$ is odd. By Lemma \ref{lem:unitaryhelper}, either $n-\lambda_\alpha\geq c$ or $\lambda$ is representable on a two-row abacus with second row having no gaps, and $n-\lambda_\alpha\geq c-2$. Since we must show that $n-\lambda_\alpha+f_2(x_\alpha)\geq c$, if we prove that $f_2(x_\alpha)\geq 2$ then we are done.

If $x_\alpha$ lies on the second row of the abacus then it is easy to see that $x_\alpha\leq e-3$, and if $x_\alpha<e-3$ then $f_2(x_\alpha)\geq 2$ as required. Hence $x_\alpha=e-3$ and there are no odd beads at all between $e-1$ and $2e-3$; hence $X=\{1,\dots,e-1,2e-1\}$ and $\lambda=(e-1)$, and $n-\lambda_\alpha+f_2(\lambda)=e\geq c$, as required. Therefore we may assume that $x_\alpha\geq e$ lies on the first row of the abacus.

If $x_\alpha\neq x_2$ then there are at least two odd beads to the right of $x-\alpha$, and so $f_2(x_\alpha)\geq 2$, as required; therefore $\alpha=2$, and $f_2(x_\alpha)=1$. Let $\bar\lambda=(\lambda_3,\dots,\lambda_s)$, noticing that removing the first two rows, which have one odd and one even hook lengths, does not change the sizes of $Y$ and $Z$. Therefore, as in Case 1, by Lemma \ref{lem:unitaryhelper} we have that $n-\lambda_1-\lambda_2=|\bar\lambda|\geq c-2$, so that 
\[ \pi(\sigma_1)\geq c-1+\lambda_1.\]
The only way this can be less than $c$ is if $\lambda_1=0$, so that $\lambda$ is the empty partition. However, in this case it is easy to see that $c=1$ and $\pi(\sigma_1)=1$, completing the proof.
\end{pf}

We have therefore completed the proof of both conditions of Corollary \ref{cor:allperversecyclic} for the unitary groups, and hence completed the proof of Theorem \ref{thm:brauertreesperverse}.

\section{Combinatorics for Symplectic and Orthogonal Groups}

For symplectic and orthogonal groups we must be slightly careful about our choice of groups, noting that at various points in the literature the `wrong' choice has been made. For the theory of unipotent characters developed here, $G$ should be a simple classical group, with \emph{diagonal} automorphisms and centre allowed; note that this \textbf{does not} allow groups such as $\SO^+_{2n}(q)$ and $\CSO_{2n}^+(q)$, as $\SO_{2n}^+(q)$ induces the graph automorphism on $\Omega_{2n}^+(q)$.

For definiteness, if $G$ is of type $B_n$ we choose $G$ to be $\SO_{2n+1}(q)$, if $G$ is of type $C_n$ we choose $\CSp_{2n}(q)$, if $G$ is of type $D_n$ or ${}^2\!D_n$ we choose $\CO_{2n}^+(q)^0$ and $\CO_{2n}^-(q)^0$ ($\PCO_{2n}^\pm(q)^0$ is the adjoint form of type ${}^\ep\!D$) as described in \cite{carterfinite} as a subgroup of index $2$ in $\CO_{2n}^\pm(q)$.

Let $\ell\nmid q$ be a prime dividing $|G|$, and write $d$ for the order of $q$ modulo $\ell$, so that $\ell\mid \Phi_d(q)$. Let $e$ be the order of $q^2$ modulo $\ell$, so that $e=d$ if $d$ is odd and $e=d/2$ if $d$ is even. The combinatorics behind the unipotent characters of $G$ and how they are distributed into blocks are very similar for all such $G$, and can be described simultaneously, using symbols (see Section \ref{sec:combobjs}).

The symbols of odd defect and a given rank $n$ parametrize the unipotent characters of the groups of type $B_n$ and $C_n$, whereas the symbols of defect divisible by $4$ correspond to unipotent characters of the groups of type $D_n$ (with two unipotent characters corresponding to each degenerate symbol), and symbols of defect congruent to $2$ modulo $4$ correspond to unipotent characters of the groups of type $^2\!D_n$.

Let $\Lambda=\{X,Y\}$ with defect $\delta$ and rank $n$, with $X=\{x_1,\dots,x_s\}$ and $Y=\{y_1,\dots,y_t\}$, ordered so that $x_i<x_{i+1}$ and $y_i<y_{i+1}$.

\medskip

In the case of $B_n$ and $C_n$, if $\chi_\Lambda$ is the unipotent character corresponding to the symbol $\Lambda$ (which has odd defect), then
\begin{equation} \chi_\lambda(1)=\frac{\D\(\prod_{i=1}^n (q^{2i}-1)\)\(\prod_{1\leq i<j\leq s}(q^{x_i}-q^{x_j})\)\(\prod_{1\leq i<j\leq t}(q^{y_i}-q^{y_j})\)\(\prod_{i,j}(q^{x_i}+q^{y_j})\)}{\D 2^{(s+t-1)/2} q^{\binom{s+t-2}{2}+\binom{s+t-4}{2}+\cdots}\(\prod_{i=1}^s\prod_{j=1}^{x_i}(q^{2j}-1)\)\(\prod_{i=1}^t\prod_{j=1}^{y_i}(q^{2j}-1)\)}.\label{eq:degreesymp}\end{equation}
As with the linear and unitary groups, this degree is invariant under the equivalence relation on symbols.

\medskip

In type $D_n$, so $G=\CO_{2n}^+(q)^0$, if $\chi_\Lambda$ is the (or `a' if $\Lambda$ is degenerate) unipotent character corresponding to the symbol $\Lambda$ (which has defect divisible by $4$), then
\begin{equation} \chi_\Lambda(1)=\frac{\D (q^n-1)\(\prod_{i=1}^{n-1}(q^{2i}-1)\)\(\prod_{1\leq i<j\leq s}(q^{x_i}-q^{x_j})\)\(\prod_{1\leq i<j\leq t}(q^{y_i}-q^{y_j})\)\(\prod_{i,j}(q^{x_i}+q^{y_j})\)}{\D 2^c q^{\binom{s+t-2}{2}+\binom{s+t-4}{2}+\cdots}\(\prod_{i=1}^s\prod_{j=1}^{x_i}(q^{2j}-1)\)\(\prod_{i=1}^t\prod_{j=1}^{y_i}(q^{2j}-1)\)},\label{eq:degreeorth1}\end{equation}
where $c=\lfloor (s+t-1)/2\rfloor$ if $X\neq Y$, and $s$ if $X=Y$. Again, this degree is invariant under the equivalence relation on symbols.

\medskip

In type $^2\!D_n$, so $G=\CSO_{2n}^-(q)$ ??????for $q$ odd and $G=\SO_{2n}^-(q)$ for $q$ even, if $\chi_\Lambda$ is the unipotent character corresponding to the symbol $\Lambda$ (which has even defect not divisible by $4$), then
\begin{equation} \chi_\Lambda(1)=\frac{\D (q^n+1)\(\prod_{i=1}^{n-1}(q^{2i}-1)\)\(\prod_{1\leq i<j\leq s}(q^{x_i}-q^{x_j})\)\(\prod_{1\leq i<j\leq t}(q^{y_i}-q^{y_j})\)\(\prod_{i,j}(q^{x_i}+q^{y_j})\)}{\D 2^c q^{\binom{s+t-2}{2}+\binom{s+t-4}{2}+\cdots}\(\prod_{i=1}^s\prod_{j=1}^{x_i}(q^{2j}-1)\)\(\prod_{i=1}^t\prod_{j=1}^{y_i}(q^{2j}-1)\)},\label{eq:degreeorth2}\end{equation}
where $c=(s+t-2)/2$. This degree is also invariant under the equivalence relation on symbols.

\medskip

In all of these groups, two unipotent characters lie in the same $\ell$-block of their respective group if and only if the corresponding symbols have the same $e$-core if $e=d$, and $e$-cocore if $e=d/2$.

\medskip

As with linear and unitary groups, we need certain evaluations of the $B_d$-function on polynomials to make our calculations easier in the next section.
\begin{prop}\label{prop:symporder} Let $d$ be a positive integer, and let $e=d$ if $d$ is odd, and $e=d/2$ if $d$ is even. Let $n\geq 0$ be an integer.
\[ B_d\(\prod_{i=n+1}^{n+e}(q^{2i}-1)\)=4ne+2e^2+2e+\frac d2\]
\end{prop}
\begin{pf}This is a simple calculation:
\begin{align*} B_d\(\prod_{i=n+1}^{n+e}(q^{2i}-1)\)&=\sum_{i=n+1}^{n+e}\(2i+d\left\lfloor\frac{2i}{d}\right\rfloor+\frac{d}{2}\)
\\ &=2ne+e(e+1)+\frac{ed}2+d\sum_{i=n+1}^{n+e}\left\lfloor\frac{2i}{d}\right\rfloor
\\ &=2ne+e(e+1)+\frac{ed}2+\begin{cases}2ne+d & e\neq d\\2ne+d(e+3)/2&e=d\end{cases}
\\ &=4ne+2e^2+2e+\frac d2.
\end{align*}
\end{pf}

For cohooks we can assume that $d$ is even. The next result has an easy proof, safely left to the reader.

\begin{prop}\label{prop:cohookdiff} Let $i$ and $j$ be integers, and let $d$ be an even integer. Write $e=d/2$, and if $j>i$ then suppose that $i-j$ is not divisible by $e$. We have
\[ B_d(q^{i+e}+q^j)-B_d(q^i-q^j)=\begin{cases}d &i-j>0
\\ 0& i-j<0\end{cases}.\]

\[ B_d(q^{i+e}-q^j)-B_d(q^i+q^j)=\begin{cases}d &i-j>-e
\\ 0& i-j<-e\end{cases}.\]
\end{prop}

\section{Brauer Trees for Classical Groups}

The Brauer trees for symplectic and orthogonal groups are very similar, and we will only give a complete treatment of the cases of $B_n$ and $C_n$, then describe the differences needed for the other orthogonal groups. The Brauer trees were described in \cite{fongsrin1990}, and we summarize their description.

\subsection{Symplectic and Odd-Dimensional Orthogonal Groups}

In this section, $d=e$ is an odd integer, $q$ is a prime power and $\ell\mid \Phi_d(q)$ is a prime. Let $G_n$ be one of the groups $\SO_{2n+1}(q)$ and $\CSp_{2n}(q)$. Let $\Lambda=\{X,Y\}$ be a symbol of rank $n$, with $X=\{x_1,\dots,x_s\}$ and $Y=\{y_1,\dots,y_t\}$, ordered so that $x_i<x_{i+1}$ and $y_i<y_{i+1}$. Assume that $\Lambda$ is an $e$-core. Recall that we view $\Lambda$ as a pair of partitions: let $X'$ denote the beads of $X$ on the end of their runners of the $e$-abacus, and let $Y'$ denote the beads of $Y$ on the end of their runners of the $e$-abacus. By choosing $\Lambda$ suitably, $|X'|=|Y'|=e$. Write $X'=\{x_1',\dots,x_e'\}$ and $Y'=\{y_1',\dots,y_e'\}$, with $x_i'>x_{i+1}'$ and $y_i'>y_{i+1}'$.

Let $\sigma_1,\dots,\sigma_e$ be the unipotent characters of $G=G_{n+e}$ corresponding to adding $e$ to the elements of $X'$, with $\sigma_i$ coming from $x_i'$; similarly, let $\tau_1,\dots,\tau_e$ be the unipotent characters of $G$ corresponding to adding $e$ to the elements of $Y'$, with $\tau_i$ coming from $y_i'$. In this case the Brauer tree is as follows.
\tikzstyle{every node}=[circle, fill=black!0,
                        inner sep=0pt, minimum width=4pt]
\begin{center}\begin{tikzpicture}[thick,scale=1.8]
\draw \foreach \x in {0,1,3,4,6,7}{
(\x,0) -- (\x+1,0)};
\draw[loosely dotted] (2.2,0) -- (2.8,0);
\draw[loosely dotted] (5.2,0) -- (5.8,0);
\draw (0,0) node [draw,label=below:$\sigma_{1}$] (l0) {};
\draw (1,0) node [draw,label=below:$\sigma_{2}$] (l1) {};
\draw (2,0) node [draw,label=below:$\sigma_{3}$] (l2) {};
\draw (3,0) node [draw,label=below:$\sigma_{e}$] (l5) {};
\draw (4,0) node [fill=black!100] (ld) {};
\draw (8,0) node [draw,label=below:$\tau_{1}$] (l0) {};
\draw (7,0) node [draw,label=below:$\tau_{2}$] (l1) {};
\draw (6,0) node [draw,label=below:$\tau_{3}$] (l2) {};
\draw (5,0) node [draw,label=below:$\tau_{e}$] (l5) {};
\end{tikzpicture}\end{center}
Since the two branches emanating from the exceptional node have the same length, the second requirement of Corollary \ref{cor:allperversecyclic} is automatically satisfied, so we have to show that $\pi(\sigma_i)\leq \pi(\sigma_{i+1})$, as with the linear and unitary groups before. (Since $\{X,Y\}=\{Y,X\}$, we do not need to show the same thing for the $\tau_i$.)

\begin{prop}\label{prop:proponesymplinear} Let $\Lambda=\{X,Y\}$ be as before, with $X'$, $Y'$, the $\sigma_i$ and $\tau_i$ as constructed. For $x\in X$, let $f_1(x)$ denote the size of the set $\{y\in Y:x-y>-d/2\}$. If $x_i'=x_\alpha$, then 
\[ \pi(\sigma_i)=4(n-x_\alpha)+2(s-\alpha+f_1(x_\alpha))+(i-1).\]
If $x_{i+1}'=x_\beta$, then 
\[ \pi(\sigma_{i+1})-\pi(\sigma_i)=4(x_\alpha-x_\beta)-2(\beta-\alpha)-2(f_1(x_\alpha)-f_1(x_\beta))+1,\]
and in particular is positive.
\end{prop}
\begin{pf} The proof of this statement follows the same pattern as that of the corresponding result for the linear and unitary groups, in other words tracking the change to $\chi_\Lambda(1)$ (Equation \ref{eq:degreesymp}) when replacing $x_\alpha$ by $x_\alpha+e$. Clearly the third term of the numerator, and first and last terms of the denominator, remain unchanged. Proposition \ref{prop:symporder} implies that the change to the first term of the numerator, and the second term of the denominator, is the $4(n-x_\alpha)$ term, and the proof of Proposition \ref{prop:GLnperversity} proves that the change to the second term of the numerator is $2(s-\alpha)+(i-1)$. In the same vein as for the unitary groups, Proposition \ref{prop:diffsinBdslinear} proves that the difference is $2f(x_\alpha)$ for the last term of the denominator, and this completes the determination of $\pi(\sigma_i)$.

That $\pi(\sigma_{i+1})-\pi(\sigma_i)$ is as claimed is obvious, so it suffices to show that it is positive. To see this, since the corresponding result holds for the $\pi$-function applied to $\GL_n(q)$, we see that $2(x_\alpha-x_\beta)-2(\beta-\alpha)\geq 0$, and so it remains to see that $2(x_\alpha-x_\beta)-2(f_1(x_\alpha)-f_1(x_\beta))\geq 0$, or equivalently that $f_1(x_\alpha)-f_1(x_\beta)\leq x_\alpha-x_\beta$; however, it is obvious that the difference between the sets $\{y\in Y:x_\beta-y>-d/2\}$ and $\{y\in Y:x_\alpha-y>-d/2\}$ is at most those numbers in $[x_\beta-d/2,x_\alpha-d/2]$, which has size $x_\alpha-x_\beta$ as needed.
\end{pf}

This proves Theorem \ref{thm:brauertreesperverse} for odd $d$, via Corollary \ref{cor:allperversecyclic}, so we now consider the more complicated case where $e=d/2$. The description of the Brauer tree is very similar to the previous case: let $\Lambda=\{X,Y\}$ be an $e$-cocore of odd defect $\delta$ and rank $n$, and let $X'$ and $Y'$ denote the subsets of $X$ and $Y$ given by
\[ X'=\{x\in X:x+e\notin Y\},\qquad Y'=\{y\in Y:y+e\notin X\}.\]
Assume that $|X|>|Y|$, so that $|X|-|Y|=\delta$. By \cite[(3E)]{fongsrin1990}, we have that $|X'|=e+\delta$ and $|Y'|=e-\delta$. Write $X'=\{x_1',\dots,x_{e+\delta}'\}$, ordered so that $x_i'>x_{i+1}'$, and similarly for $Y'$. If $\sigma_i$ is the unipotent character corresponding to the symbol obtained by adding an $e$-cohook to $x_i'$, and similarly for $\tau_i$ and $y_i'$, then the Brauer tree is as follows.
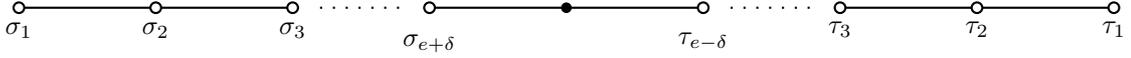
\tikzstyle{every node}=[circle, fill=black!0,
                        inner sep=0pt, minimum width=4pt]
\begin{center}\begin{tikzpicture}[thick,scale=1.8]
\draw \foreach \x in {0,1,3,4,6,7}{
(\x,0) -- (\x+1,0)};
\draw[loosely dotted] (2.2,0) -- (2.8,0);
\draw[loosely dotted] (5.2,0) -- (5.8,0);
\draw (0,0) node [draw,label=below:$\sigma_{1}$] (l0) {};
\draw (1,0) node [draw,label=below:$\sigma_{2}$] (l1) {};
\draw (2,0) node [draw,label=below:$\sigma_{3}$] (l2) {};
\draw (3,0) node [draw,label=below:$\sigma_{e+\delta}$] (l5) {};
\draw (4,0) node [fill=black!100] (ld) {};
\draw (8,0) node [draw,label=below:$\tau_{1}$] (l0) {};
\draw (7,0) node [draw,label=below:$\tau_{2}$] (l1) {};
\draw (6,0) node [draw,label=below:$\tau_{3}$] (l2) {};
\draw (5,0) node [draw,label=below:$\tau_{e-\delta}$] (l5) {};
\end{tikzpicture}\end{center}

With the picture above, it is fairly easy to prove the first of the two properties needed for $\pi(-)$ to induce a perverse equivalence.

\begin{prop}\label{prop:proponesympunit} Let $\Lambda=\{X,Y\}$ be as before, with $X'$, $Y'$, the $\sigma_i$ and $\tau_i$ as constructed. For $x\in X$, let $f_1(x)$ denote the size of the set $\{y\in Y:x-y>-d/2\}$. If $x_i'=x_\alpha$, then 
\[ \pi(\sigma_i)=2(n-x_\alpha)+(s-\alpha)+f_1(x_\alpha).\]
If $x_{i+1}'=x_\beta$, then 
\[ \pi(\sigma_{i+1})-\pi(\sigma_i)=2(x_\alpha-x_\beta)-(\beta-\alpha)-(f_1(x_\alpha)-f_1(x_\beta)),\]
and in particular is positive.
\end{prop}
\begin{pf} We follow the same strategy as with Proposition \ref{prop:proponesymplinear}, noting that removing an element of $X$ and adding it to $Y$ does not alter the quantity $s+t$. The first term, $2(n-x_\alpha)$, is produced exactly as in Proposition \ref{prop:proponesymplinear}. For the rest of the terms, we need to identify the effect on the last three terms of the numerator of Equation \ref{eq:degreesymp} simultaneously, using Proposition \ref{prop:cohookdiff}. We see that some terms (corresponding to $x\in X$) move from the second term in the numerator to the last term in the numerator, and clearly these contribute $s-\alpha$, as this is the size of the set $\{x\in X:x_\alpha>x\}$. Finally, some terms (corresponding to $y\in Y$) move from the last term in the numerator to the third term, and these contribute $f_1(x-\alpha)$, as this is the size of the set $\{y\in Y:x-y>-d/2\}$. (Both of these determinations use Proposition \ref{prop:cohookdiff}, of course.)

Clearly the difference $\pi(\sigma_{i+1})-\pi(\sigma_i)$ is as described, so it suffices to show that it is positive. However, this is simply the difference in Proposition \ref{prop:proponesymplinear}, minus $1$, and halved, thus is non-negative. However, since the $\pi$-function, applied to adjacent unipotent characters in the Brauer tree, must have different parities, means that $\pi(\sigma_{i+1})-\pi(\sigma_i)$ is positive, as needed.
\end{pf}

Thus it remains to deal with the second requirement of Corollary \ref{cor:allperversecyclic}, namely that $\pi(-)$ is at least the canonical perversity function. As with the unitary groups, this only needs to be checked on $\tau_1$.

\begin{prop}\label{prop:proptwosympunit} Let $\Lambda=\{X,Y\}$ be an $e$-cocore as above, and let $\lambda$ and $\mu$ be the partitions with $\beta$-sets $X$ and $Y$ respectively. Let $|X|=s$, $|Y|=t$, and $\delta=s-t$ is positive and odd. Let $X'$ and $Y'$, and the $\sigma_i$ and $\tau_i$ be as above. We have that $\pi(\tau_1)\geq 2\delta$, and consequently $\pi(\tau_1)\geq \pi_0(\tau_1)$, where $\pi_0(-)$ is the canonical perversity function.
\end{prop}
\begin{pf} We first notice that if $e=1$ then $\delta=1$, and so $Y'$ is empty; thus we may assume that $e\geq 2$. Let $y_1'=y_\alpha$. Choose $u$ and $v$ maximal such that $\{0,\dots,u-1\}\subs X$ and $\{0,\dots,v-1\}\subs Y$. (It could be that either $u$ or $v$, or even both, are $0$.) Notice that $\lambda$ has $s-u$ rows and $\mu$ has $t-v$ rows.

We begin by noting that $\sum_{x\in X} x=s(s-1)/2+|\lambda|$, and similarly $\sum_{y\in Y} y=t(t-1)/2+|\mu|$. Hence
\begin{align*} n&=\sum_{x\in X} x+\sum_{y\in Y} y-\frac{(s+t-1)^2}4
\\ &=|\lambda|+|\mu|+\frac{s^2}{4}+\frac{t^2}{4}-\frac{st}{2}-\frac{1}{4}
\\&=|\lambda|+|\mu|+\frac{\delta^2-1}{4}.\end{align*}
Notice that $\pi_0(\tau_1)=2\delta$ or $\pi_0(\tau_1)=2\delta+1$; since $\pi(-)$ and $\pi_0(-)$ have the same parity, it suffices to show that $\pi(\tau_1)-2\delta$ is non-negative. We get that
\begin{align*} \pi(\tau_1)-2\delta&=2(n-y_\alpha)+(t-\alpha)+f_1(y_\alpha)-2\delta
\\ &=\(\frac{(\delta^2-1)}2-2\delta\)-2y_\alpha+2|\lambda|+2|\mu|+(t-\alpha)+f_1(y_\alpha).\end{align*}
Since $\delta$ is a positive odd integer, $(\delta^2-1)/2-2\delta\geq -2$. In addition, $y_\alpha+(t-\alpha)=\mu_\alpha$, (which might be $0$), and therefore $|\mu|\geq y_\alpha+(t-\alpha)$.  Thus it remains to prove that
\begin{equation} 2|\lambda|+|\mu|-y_\alpha+f_1(y_\alpha)\geq 2.\label{eq:sympperv}\end{equation}

\medskip

\noindent\textbf{Case 1}: $u\leq v$. Since $s>t$, $\lambda$ has at least $v-u+1$ rows, so $|\lambda|\geq v-u+1$. If $y_\alpha<v$ then $f_1(y_\alpha)+|\lambda|\geq y_\alpha+1$, and since $|\lambda|\geq 1$,
\[ |\lambda|+(|\lambda|+f_1(y_\alpha)-y_\alpha)\geq 2,\]
as needed. If $y_\alpha>v$ then $f_1(y_\alpha)\geq u$, and so $|\lambda|+f_1(y_\alpha)\geq v+1$. However, the first-column hook lengths of $\mu$ are the $y_i-v$ (whenever this is positive) and so $y_\alpha-v\leq |\mu|$. Hence
\[ |\lambda|+(|\lambda|+f_1(y_\alpha)+|\mu|-y_\alpha)\geq 2,\]
proving this case.

\medskip

\noindent\textbf{Case 2}: $u>v$. If $\lambda$ is not the empty partition then $2|\lambda|\geq 2$, and for Equation \ref{eq:sympperv} it suffices to show that $|\mu|-y_\alpha+f_1(y_\alpha)\geq 0$.  If $y_\alpha<v$ then $f_1(y_\alpha)\geq y_\alpha$ clearly, and if $y_\alpha>v$ then, as in Case 1, $y_\alpha-v\leq |\mu|$, and $f_1(y_\alpha)\geq v$, so that $|\mu|+f_1(y_\alpha)\geq y_\alpha$, as claimed. Hence $|\lambda|=0$, and $X=\{0,\dots,s-1\}$. Since $f_1(y_\alpha)$ counts all elements of $X$ from $0$ up to (but not including) $y_\alpha+e$, and $e$ cannot be $1$, if $y_\alpha<v$ then $f_1(y_\alpha)\geq y_\alpha+2$, and so $f_1(y_\alpha)-y_\alpha\geq 2$, as needed. Finally, if $y_\alpha>v$ then $t>v$, so that $s=u\geq v+2$. As in the previous argument, $y_\alpha-v\leq |\mu|$, and $f_1(y_\alpha)\geq v+2$, so that $|\mu|+f_1(y_\alpha)-y_\alpha\geq 2$, completing the final case.
\end{pf}

With the proof of this proposition, Corollary \ref{cor:allperversecyclic} now implies Theorem \ref{thm:brauertreesperverse} for symplectic and odd-dimensional orthogonal groups.

\subsection{Even-Dimensional Orthogonal Groups}

In this case, the $d$-split Levi subgroups in cuspidal pairs involve other orthogonal groups, but the sign involved might change (see \cite[\S3]{bmm1993}). More precisely, removing an $e$-hook for a group of type $D_n^\ep$ still results in a group of type $D_{n-e}^\ep$, but removing an $e$-cohook for a group of type $D_n^\ep$ results in a group of type $D_{n-e}^{-\ep}$; as our blocks have weight $1$, the $e$-cocore involved will come from a group of different sign to the block itself.

It will be possible to treat both the plus-type and minus-type orthogonal groups simultaneously. The first point to notice when comparing Equations \ref{eq:degreesymp}, \ref{eq:degreeorth1} and \ref{eq:degreeorth2} is that the only difference is one term, which is $(q^{2n}-1)$, $(q^n-1)$ and $(q^n+1)$ respectively. Clearly, whichever of these terms is present does not affect the differences (using the notation of the previous section) $\pi(\sigma_{i+1})-\pi(\sigma_i)$, and so the second assertion in Proposition \ref{prop:proponesympunit}, that this difference is positive, still holds. However, the proof of Proposition \ref{prop:proptwosympunit} is slightly different, as the defect $\delta$ is now even, and whether the term is $(q^{2n}-1)$, $(q^n-1)$ or $(q^n+1)$ does affect (slightly) the quantity $\pi(\tau_1)$.

We now prove the analogue of Proposition \ref{prop:proptwosympunit}, noting that we may assume that the defect $\delta$ is positive, as else the number of nodes on either side of the exceptional node are the same, and there is nothing to prove. Hence let $d$ be a positive even integer, $e=d/2$, and let $\Lambda=\{X,Y\}$ denote an $e$-cocore. Writing $|X|=s$ and $|Y|=t$, we may assume that $s>t$, and $\delta=s-t$ is even, so that $\Lambda$ is non-degenerate. Write $n$ for the rank of $\Lambda$, let $G_n$ denote the orthogonal group $\PCO_{2n}^\ep(q)^0$, where $\ep$ is either $1$ or $-1$ depending on whether $4\mid \delta$ or $4\nmid\delta$ respectively, and let $G_{n+e}=\PCO_{2n}^{-\ep}(q)^0$. Denote by $X'$ the set of all $x\in X$ such that $x+e\notin Y$, and similarly for $Y'$, ordered so that $x_i'>x_{i+1}'$ and $y_i'>y_{i+1}'$. By an equivalent choice of $\Lambda$, we have that $|X'|=e+\delta$ and $|Y'|=e-\delta$. We denote by $\sigma_i$ the unipotent character of $G_{n+e}$ whose corresponding symbol is $\Lambda$ with the cohook corresponding to $x_i'$ added, and similarly for $\tau_i$ and $y_i'$. 

If $\Lambda$ is non-degenerate then the Brauer tree for this block is exactly the same as the one in the previous section. However, if $\Lambda$ is degenerate (so $\ep=-1$ in particular) then only one branch of the Brauer tree exists, and the tree looks like that for $\GL_n(q)$. Notice that (for $d$ even) this differs from the tree given in \cite{fongsrin1990}; whether there is a single branch or two branches for the group $\CSO_{2n}^-(q)$, the group actually treated in \cite{fongsrin1990}, depends on the situation. For example, for $\CSO_6^-(3)$, $\ell=7\mid\Phi_6(3)$, there is a single branch, but for $\CSO_8^-(3)$, $\ell=41\mid\Phi_8(3)$, there are two branches, one branch containing non-unipotent characters.


In the next result we use all of this notation.

\begin{prop}\label{prop:proptwoorthogunit} With the notation above, we have that $\pi(\tau_1)\geq 2\delta$, and therefore $\pi(\tau_1)\geq\pi_0(\tau_1)$, where $\pi_0(-)$ is the canonical perversity function.
\end{prop}
\begin{pf} We first notice that $e\geq 2$ since $e\geq \delta$, and if $e=2$ then $\delta=2$, and so $Y'$ is empty. Thus we may assume that $e\geq 3$. Let $y_1'=y_\alpha$. Choose $u$ and $v$ maximal such that $\{0,\dots,u-1\}\subs X$ and $\{0,\dots,v-1\}\subs Y$. (It could be that either $u$ or $v$, or even both, are $0$.) Notice that $\lambda$ has $s-u$ rows and $\mu$ has $t-v$ rows.

We begin by noting that $\sum_{x\in X} x=s(s-1)/2+|\lambda|$, and similarly $\sum_{y\in Y} y=t(t-1)/2+|\mu|$. Hence
\begin{align*} n&=\sum_{x\in X} x+\sum_{y\in Y} y-\left\lfloor\frac{(s+t-1)^2}4\right\rfloor
\\ &=|\lambda|+|\mu|+\frac{s^2}{4}+\frac{t^2}{4}-\frac{st}{2}
\\&=|\lambda|+|\mu|+\frac{\delta^2}{4}.\end{align*}
The determination of $\pi(\tau_1)$ is very similar to the odd-dimensional case, except that the $2n$ term must be replaced with $2(n-1)+1=2n-1$, which is
\[ \left.\left[ B_d\((q^{n+e}\pm 1)\prod_{i=1}^{n+e-1}(q^{2i}-1)\)-B_d\((q^n\mp 1)\prod_{i=1}^{n-1}(q^{2i}-1)\)\right]\right/d=2(n-1)+1,\]
via Propositions \ref{prop:symporder} and \ref{prop:cohookdiff}.

As before, $\pi_0(\tau_1)$ is either $2\delta$ or $2\delta+1$, and so by the parity argument it suffices to prove that $\pi(\tau_1)-2\delta$ is non-negative. We get that
\begin{align*} \pi(\tau_1)-2\delta&=2(n-y_\alpha)+(t-\alpha)+f_1(y_\alpha)-2\delta-1
\\ &=\(\frac{\delta^2}2-2\delta\)-2y_\alpha+2|\lambda|+2|\mu|+(t-\alpha)+f_1(y_\alpha)-1.\end{align*}
Since $\delta$ is a positive even integer, $\delta^2/2-2\delta\geq -2$. In addition, $y_\alpha+(t-\alpha)=\mu_\alpha$ (which might be $0$), and therefore $|\mu|\geq y_\alpha+(t-\alpha)$.  Thus it remains to prove that
\begin{equation} 2|\lambda|+|\mu|-y_\alpha+f_1(y_\alpha)\geq 3.\label{eq:orthogperv}\end{equation}

\medskip

\noindent\textbf{Case 1}: $u\leq v$. Since $s-t\geq 2$, $\lambda$ has at least $v-u+2$ rows, so $|\lambda|\geq v-u+2$. If $y_\alpha<v$ then $f_1(y_\alpha)+|\lambda|\geq y_\alpha+1$, and since $|\lambda|\geq 2$,
\[ |\lambda|+(|\lambda|+f_1(y_\alpha)-y_\alpha)\geq 3,\]
as needed. If $y_\alpha>v$ then $f_1(y_\alpha)\geq u$, and so $|\lambda|+f_1(y_\alpha)\geq v+2$. However, the first-column hook lengths of $\mu$ are the $y_i-v$ (whenever this is positive) and so $y_\alpha-v\leq |\mu|$. Hence, since $|\lambda|\geq 1$,
\[ |\lambda|+(|\lambda|+f_1(y_\alpha)+|\mu|-y_\alpha)\geq 3,\]
proving this case.

\medskip

\noindent\textbf{Case 2}: $u>v$. If $\lambda$ is not the empty partition then $2|\lambda|\geq 2$, and for Equation \ref{eq:orthogperv} it suffices to show that $|\mu|-y_\alpha+f_1(y_\alpha)\geq 1$.  If $y_\alpha<v$ then $f_1(y_\alpha)\geq y_\alpha+1$, and if $y_\alpha>v$ then, as in Case 1, $y_\alpha-v\leq |\mu|$, and $f_1(y_\alpha)\geq v+1$, so that $|\mu|+f_1(y_\alpha)\geq y_\alpha+1$, as claimed. Hence $|\lambda|=0$ and $X=\{0,\dots,s-1\}$. Since $f_1(y_\alpha)$ counts all elements of $X$ from $0$ up to (but not including) $y_\alpha+e$, and $e\geq 3$, if $y_\alpha<v$ then $f_1(y_\alpha)\geq y_\alpha+3$ (as $s-v\geq 2$), and so $f_1(y_\alpha)-y_\alpha\geq 3$, as needed. Finally, if $y_\alpha>v$ then $t>v$, so that $s=u\geq v+3$. As in the previous argument, $y_\alpha-v\leq |\mu|$, and $f_1(y_\alpha)\geq v+3$, so that $|\mu|+f_1(y_\alpha)-y_\alpha\geq 3$, completing the final case.
\end{pf}

This completes the proof of Theorem \ref{thm:brauertreesperverse} for the even-dimensional orthogonal groups, and hence for all classical groups.

\section{The Exceptional Cases}

%

For the exceptional groups, we need to prove Theorem \ref{thm:brauertreesperverse}, at least except for types $E_7$ and $E_8$, where the Brauer trees are not known. In addition, for the smaller groups ($G_2(q)$, $^3\!D_4(q)$, $F_4(q)$) information about the decomposition matrices is known, and we may check that our ordering on the simple modules produces a lower triangular shape. We do not make guesses for the Brauer trees of $E_7$ and $E_8$ here, and delay this to a later paper, where we also discuss the bijection between the simple modules for the block and its Brauer correspondent.

We will also give the $\pi$-function on the unipotent characters for those $d$ for which the Sylow $\ell$-subgroup is non-cyclic, except for $d=1,2$, where the $\pi$-function is easily calculable (see Proposition \ref{prop:conjd12}).

\subsection{$G_2(q)$}
\label{sec:G2}
We first consider the group $G_2(q)$, which has order $q^6\Phi_1^2\Phi_2^2\Phi_3\Phi_6$. The only Brauer trees for $G_2(q)$ are for the principal block and $d=3,6$. By \cite{shamash1989}, the Brauer trees for $G_2(q)$ are determined; labelled with the $\pi$-function, these are as follows, for $\ell\mid \Phi_3$ and $\ell\mid\Phi_6$ respectively.
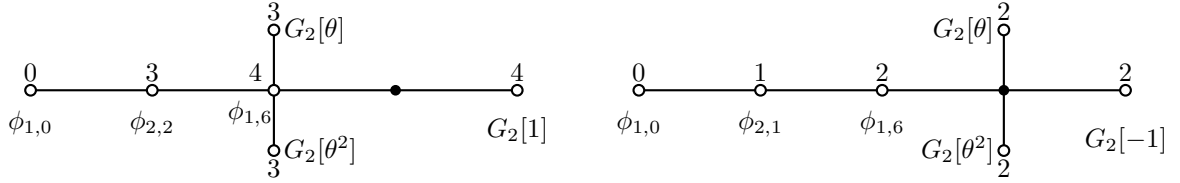
\begin{center}\begin{tikzpicture}[thick,scale=1.6]
\draw (0,0.5) -- (4,0.5);
\draw (2,0) -- (2,1);

\draw (0,0.5) node [draw,label=below:$\phi_{1,0}$] (l0) {};
\draw (1,0.5) node [draw,label=below:$\phi_{2,2}$] (l1) {};
\draw (4,0.5) node [draw,label=below:${G_2[1]}$] (l2) {};
\draw (3,0.5) node [fill=black!100] (ld) {};
\draw (2,0) node [draw,label=right:${G_2[\theta^2]}$] (l4) {};
\draw (2,0.5) node [draw,label=below left:$\phi_{1,6}$] (l4) {};
\draw (2,1) node [draw,label=right:${G_2[\theta]}$] (l4) {};
\draw (0,0.65) node{$0$};
\draw (1,0.65) node{$3$};
\draw (1.85,0.65) node{$4$};
\draw (4,0.65) node{$4$};
\draw (2,1.15) node{$3$};
\draw (2,-0.15) node{$3$};

\draw (5,0.5) -- (9,0.5);
\draw (8,0) -- (8,1);

\draw (5,0.5) node [draw,label=below:$\phi_{1,0}$] (l0) {};
\draw (6,0.5) node [draw,label=below:$\phi_{2,1}$] (l1) {};
\draw (7,0.5) node [draw,label=below:$\phi_{1,6}$] (l2) {};
\draw (8,0.5) node [fill=black!100] (ld) {};
\draw (8,0) node [draw,label=left:${G_2[\theta^2]}$] (l4) {};
\draw (9,0.5) node [draw,label=below:${G_2[-1]}$] (l4) {};
\draw (8,1) node [draw,label=left:${G_2[\theta]}$] (l4) {};

\draw (5,0.65) node{$0$};
\draw (6,0.65) node{$1$};
\draw (7,0.65) node{$2$};
\draw (9,0.65) node{$2$};
\draw (8,1.15) node{$2$};
\draw (8,-0.15) node{$2$};

\end{tikzpicture}\end{center}
This proves Theorem \ref{thm:brauertreesperverse} for $G_2(q)$. Finally, if $\ell\mid (q+1)$, we can check that the decomposition matrix is triangular with respect to the $\pi$-function: the decomposition matrix for the principal $\ell$-block of $G_2(q)$ is given in \cite{hiss1989}, and it indeed satisfies the triangularity condition.

\subsection{$^3\!D_4(q)$}

The group $^3\!D_4(q)$ has order $q^{12}\Phi_1^2\Phi_2^2\Phi_3^2\Phi_6^2\Phi_{12}$; much information on the decomposition numbers, and the single Brauer tree for $d=12$, was determined in \cite{geck1991}. There are no non-principal unipotent $\ell$-blocks that have non-trivial defect group, and so in this case we produce a table of the unipotent characters and their $\pi$-functions for $d=3,6$, together with the labelled Brauer tree for $d=12$.

\begin{center}\begin{tabular}{l|cc||l|cc}
Name & $d=3$ & $d=6$ & Name & $d=3$ & $d=6$
\\ \hline $\phi_{1,0}$& 0 & 0 & $^3\!D_4[-1]$ & - & 3
\\ $\phi_{1,3}'$& 3 &  2 & $^3\!D_4[1]$& 6 & 4 
\\ $\phi_{2,2}$& 5 & 3  &$\phi_{1,3}''$ & 7 &4
\\$\phi_{2,1}$& 6 & - & $\phi_{1,6}$ & $8$ & $4$
\end{tabular}
\end{center}


Here is the Brauer tree for $d=12$.
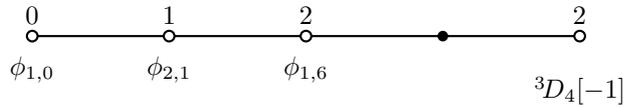
\begin{center}\begin{tikzpicture}[thick,scale=1.8]
\draw \foreach \x in {0,1,2,3}{
(\x,0) -- (\x+1,0)};
\draw (0,0) node [draw,label=below:$\phi_{1,0}$] (l0) {};
\draw (1,0) node [draw,label=below:$\phi_{2,1}$] (l1) {};
\draw (2,0) node [draw,label=below:$\phi_{1,6}$] (l2) {};
\draw (3,0) node [fill=black!100] (ld) {};
\draw (4,0) node [draw,label=below:${^3\!D_4[-1]}$] (l4) {};
\draw (0,0.15) node{$0$};
\draw (1,0.15) node{$1$};
\draw (2,0.15) node{$2$};
\draw (4,0.15) node{$2$};
\end{tikzpicture}\end{center}

We can also check that the decomposition matrices given in \cite{geck1991} are triangular with respect to the $\pi$-function: for $d=2,3$ this is true, but for $d=6$ it requires, in the notation of \cite{geck1991}, $c=d=0$, where these are currently unknown parameters.

\subsection{$F_4(q)$}
\label{sec:F4}
The group $F_4(q)$ has order $q^{24}\Phi_1^4\Phi_2^4\Phi_3^3\Phi_4^2\Phi_6^2\Phi_8\Phi_{12}$. Much of the structure of the decomposition matrix of $F_4(q)$ was discovered in \cite{kohler2006}, and the Brauer trees for $F_4(q)$ were determined in \cite{hisslubeck1998}, with the planar embedding for the case $d=12$ determined in \cite[Theorem 3.12(ii)]{dudas2010un2}. We give a table of the $\pi$-function for the cases where the Sylow $\ell$-subgroup is non-cyclic, with $d\geq 3$. (Here $B_i$ denotes a non-principal block of $F_4(q)$.)
\begin{center}
\begin{tabular}{l|ccc||l|ccc}
Degree & $d=3$ & $d=4$ & $d=6$ & Degree & $d=3$ & $d=4$ & $d=6$ 
\\ \hline
$\phi_{1,0}$&0&0&0&$\phi_{9,2}$&-&7&-
\\$\phi_{1,24}$&16&12&8&$\phi_{9,10}$&-&11&-
\\$\phi_{1,12}'$&14&-&-&$\phi_{9,6}'$&-&-&6
\\$\phi_{1,12}''$&14&-&-&$\phi_{9,6}''$&-&-&6
\\$\phi_{2,4}'$&7&$B_2$: 4&4&$\phi_{12,4}$&-&8&6
\\$\phi_{2,16}'$&15&$B_2$: 10&8&$\phi_{16,5}$&13&-&-
\\$\phi_{2,4}''$&7&$B_3$: 4&4&$B_2,1$&-&5&4
\\$\phi_{2,16}''$&15&$B_3$: 10&8&$B_2,r$&-&10&-
\\$\phi_{4,1}$&7&6&-&$B_2,\ep$&-&11&8
\\$\phi_{4,7}'$&13&$B_2$: 9&-&$B_2,\ep'$&-&$B_2$: 8&7
\\$\phi_{4,7}''$&13&$B_3$: 9&-&$B_2,\ep''$&-&$B_3$: 8&7
\\$\phi_{4,8}$&12&10&-&$F_4[-1]$&-&-&7
\\$\phi_{4,13}$&15&12&-&$F_4[\I]$&-&10&-
\\$\phi_{6,6}'$&-&10&-&$F_4[-\I]$&-&10&-
\\$\phi_{6,6}''$&-&-&-&$F_4[\theta]$&13&-&7
\\$\phi_{8,3}'$&10&-&5&$F_4[\theta^2]$&13&-&7
\\$\phi_{8,9}'$&14&-&7&$F_4^{\mathrm{I}}[1]$&-&10&8
\\$\phi_{8,3}''$&10&-&5&$F_4^{\mathrm{II}}[1]$&14&12&-
\\$\phi_{8,9}''$&14&-&7
\end{tabular}\end{center}
For $d=4$ there are two non-principal unipotent $\ell$-blocks with cyclic defect group, but for all other $d$, the only such blocks have Sylow $\ell$-subgroups as defect groups. Here are the trees for $d=4$.
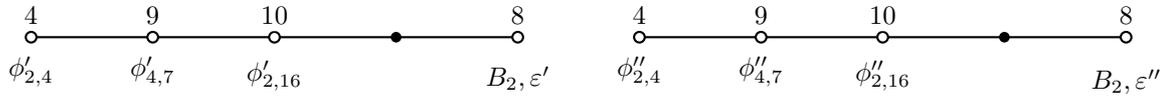
\begin{center}\begin{tikzpicture}[thick,scale=1.6]
\draw (0,0) -- (4,0);
\draw (0,0) node [draw,label=below:$\phi_{2,4}'$] (l0) {};
\draw (1,0) node [draw,label=below:$\phi_{4,7}'$] (l1) {};
\draw (2,0) node [draw,label=below:$\phi_{2,16}'$] (l2) {};
\draw (3,0) node [fill=black!100] (ld) {};
\draw (4,0) node [draw,label=below:${B_2,\ep'}$] (l4) {};
\draw (0,0.18) node{$4$};
\draw (1,0.18) node{$9$};
\draw (2,0.18) node{$10$};
\draw (4,0.18) node{$8$};

\draw (5,0) -- (9,0);
\draw (5,0) node [draw,label=below:$\phi_{2,4}''$] (l0) {};
\draw (6,0) node [draw,label=below:$\phi_{4,7}''$] (l1) {};
\draw (7,0) node [draw,label=below:$\phi_{2,16}''$] (l2) {};
\draw (8,0) node [fill=black!100] (ld) {};
\draw (9,0) node [draw,label=below:${B_2,\ep''}$] (l4) {};
\draw (5,0.18) node{$4$};
\draw (6,0.18) node{$9$};
\draw (7,0.18) node{$10$};
\draw (9,0.18) node{$8$};
\end{tikzpicture}\end{center}
For $d=8$, we have only the principal $\ell$-block.
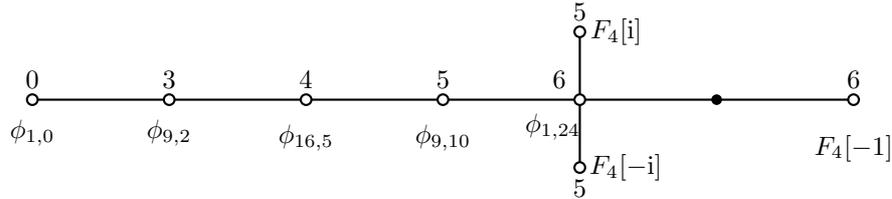
\begin{center}\begin{tikzpicture}[thick,scale=1.8]
\draw (0,0.5) -- (6,0.5);
\draw (4,0) -- (4,1);

\draw (0,0.5) node [draw,label=below:$\phi_{1,0}$] (l0) {};
\draw (1,0.5) node [draw,label=below:$\phi_{9,2}$] (l1) {};
\draw (2,0.5) node [draw,label=below:$\phi_{16,5}$] (l2) {};
\draw (3,0.5) node [draw,label=below:$\phi_{9,10}$] (l2) {};
\draw (4,0) node [draw,label=right:${F_4[-\I]}$] (l4) {};
\draw (4,0.5) node [draw,label=below left:$\phi_{1,24}$] (l4) {};
\draw (4,1) node [draw,label=right:${F_4[\I]}$] (l4) {};
\draw (5,0.5) node [fill=black!100] (ld) {};
\draw (6,0.5) node [draw,label=below:${F_4[-1]}$] (l2) {};
\draw (0,0.65) node{$0$};
\draw (1,0.65) node{$3$};
\draw (2,0.65) node{$4$};
\draw (3,0.65) node{$5$};
\draw (3.85,0.65) node{$6$};
\draw (4,1.15) node{$5$};
\draw (4,-0.15) node{$5$};
\draw (6,0.65) node{$6$};
\end{tikzpicture}\end{center}
Finally, for $d=12$, we again have only the principal $\ell$-block.
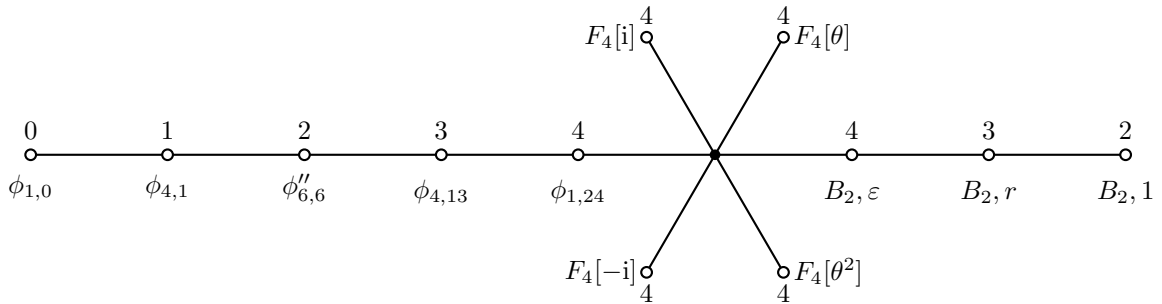
\begin{center}\begin{tikzpicture}[thick,scale=1.8]
\draw (0,1) -- (8,1);
\draw (4.5,0.1339) -- (5.5,1.866);
\draw (4.5,1.866) -- (5.5,0.1339);
\draw (0,1) node [draw,label=below:$\phi_{1,0}$] (l0) {};
\draw (1,1) node [draw,label=below:$\phi_{4,1}$] (l0) {};
\draw (2,1) node [draw,label=below:$\phi_{6,6}''$] (l0) {};
\draw (3,1) node [draw,label=below:$\phi_{4,13}$] (l0) {};
\draw (4,1) node [draw,label=below:$\phi_{1,24}$] (l0) {};
\draw (5,1) node [fill=black!100] (ld) {};
\draw (6,1) node [draw,label=below:${B_2,\ep}$] (l0) {};
\draw (7,1) node [draw,label=below:${B_2,r}$] (l0) {};
\draw (8,1) node [draw,label=below:${B_2,1}$] (l0) {};
\draw (0,1.18) node{$0$};
\draw (1,1.18) node{$1$};
\draw (2,1.18) node{$2$};
\draw (3,1.18) node{$3$};
\draw (4,1.18) node{$4$};
\draw (6,1.18) node{$4$};
\draw (7,1.18) node{$3$};
\draw (8,1.18) node{$2$};

\draw (4.5,1.866) node [draw,label=left:${F_4[\I]}$] (l0) {};
\draw (5.5,1.866) node [draw,label=right:${F_4[\theta]}$] (l0) {};
\draw (4.5,0.1339) node [draw,label=left:${F_4[-\I]}$] (l0) {};
\draw (5.5,0.1339) node [draw,label=right:${F_4[\theta^2]}$] (l0) {};

\draw (4.5,2.03) node{$4$};
\draw (5.5,2.03) node{$4$};
\draw (4.5,-0.02) node{$4$};
\draw (5.5,-0.02) node{$4$};
\end{tikzpicture}\end{center}

This completes the proof of Theorem \ref{thm:brauertreesperverse} for $F_4(q)$. For the principal $\ell$-block, one may check that the ordering on the unipotent characters given by the perversity function yields a decomposition matrix with lower triangular shape, using the tables in \cite{kohler2006}. The author has done this for $d=3,4,6$, the cases where our conjecture is new: for $d=3$, the decomposition matrix is definitely lower triangular this with ordering, but for $d=4$ and $d=6$ this requires certain currently unknown entries to be zero. (For the $d=4$ table in \cite{kohler2006}, we need $a$, $b$ and $f$ to be 0, and for the $d=6$ table, we need $b$, $c$, $d$, $f$, and the $\ast$ in the $\phi_{8,9}'$ and $\phi_{8,9}''$ rows to be 0. For each of these, 0 is an acceptable value.)

\subsection{$E_6(q)$}

The group $G=E_6(q)$ has order $q^{36}\Phi_1^6\Phi_2^4\Phi_3^3\Phi_4^2\Phi_5\Phi_6^2\Phi_8\Phi_9\Phi_{12}$; the block structure can be deduced from \cite{bmm1993}; for the blocks with cyclic defect groups we give the Brauer trees, which are given in \cite{hlm1995}, with the $\pi$-function attached.

For $d=3,4,6$, we give the table of the $\pi$-function.
\begin{center}\begin{tabular}{l|ccc||l|ccc}
Name & $d=3$ & $d=4$ & $d=6$ & Name & $d=3$ & $d=4$ & $d=6$
\\ \hline
$\phi_{1,0}$&0&0&0&$\phi_{30,15}$&22&-&10
\\$\phi_{1,36}$&24&18&12&$\phi_{60,8}$&19&-&9
\\$\phi_{10,9}$&19&16&-&$\phi_{80,7}$&17&14&9
\\$\phi_{6,1}$&7&6&4&$\phi_{90,8}$&-&14&-
\\$\phi_{6,25}$&23&18&12&$\phi_{60,5}$&16&$B_2$: 11&8
\\$\phi_{20,10}$&20&-&-&$\phi_{60,11}$&20&$B_2$: 14&10
\\$\phi_{15,5}$&14&12&-&$\phi_{64,4}$&15&-&-
\\$\phi_{15,17}$&22&18&-&$\phi_{64,13}$&21&-&-
\\$\phi_{15,4}$&13&11&7&$\phi_{81,6}$&-&13&-
\\$\phi_{15,36}$&21&17&11&$\phi_{81,10}$&-&15&-
\\$\phi_{20,2}$&10&$B_2$: 6&5&$D_4,1$&$B_2$: 8&11&7
\\$\phi_{20,20}$&22&$B_2$: 15&11&$D_4,\ep$&$B_2$: 16&17&11
\\$\phi_{24,6}$&16&-&8&$D_4,r$&$B_2$: 13&15&10
\\$\phi_{24,12}$&20&-&10&$E_6[\theta]$&19&-&10
\\$\phi_{30,3}$&14&-&6&$E_6[\theta^2]$&19&-&10
\end{tabular}
\end{center}
For $d=2$ there is a single unipotent Brauer tree, in a non-principal $\ell$-block.
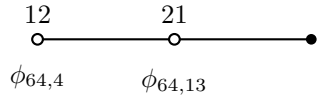
\begin{center}\begin{tikzpicture}[thick,scale=1.8]
\draw (0,0) -- (2,0);
\draw (0,0) node [draw,label=below:$\phi_{64,4}$] (l0) {};
\draw (1,0) node [draw,label=below:$\phi_{64,13}$] (l1) {};
\draw (2,0) node [fill=black!100] (ld) {};
\draw (0,0.18) node{$12$};
\draw (1,0.18) node{$21$};
\end{tikzpicture}\end{center}
For $d=3$ we get the following, non-principal Brauer tree.
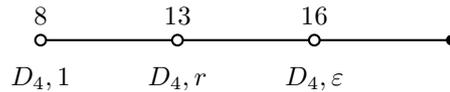
\begin{center}\begin{tikzpicture}[thick,scale=1.8]
\draw (0,0) -- (3,0);
\draw (0,0) node [draw,label=below:${D_4,1}$] (l0) {};
\draw (1,0) node [draw,label=below:${D_4,r}$] (l1) {};
\draw (2,0) node [draw,label=below:${D_4,\ep}$] (l2) {};
\draw (3,0) node [fill=black!100] (ld) {};
\draw (0,0.18) node{$8$};
\draw (1,0.18) node{$13$};
\draw (2,0.18) node{$16$};
\end{tikzpicture}\end{center}
For $d=4$ the picture is similar to $d=3$, with one small non-principal block.
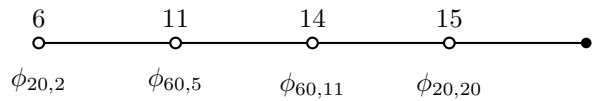
\begin{center}\begin{tikzpicture}[thick,scale=1.8]
\draw (0,0) -- (4,0);
\draw (0,0) node [draw,label=below:$\phi_{20,2}$] (l0) {};
\draw (1,0) node [draw,label=below:$\phi_{60,5}$] (l1) {};
\draw (2,0) node [draw,label=below:$\phi_{60,11}$] (l2) {};
\draw (3,0) node [draw,label=below:$\phi_{20,20}$] (l2) {};
\draw (4,0) node [fill=black!100] (ld) {};
\draw (0,0.18) node{$6$};
\draw (1,0.18) node{$11$};
\draw (2,0.18) node{$14$};
\draw (3,0.18) node{$15$};
\end{tikzpicture}\end{center}
For $d=5$ we have two unipotent $\ell$-blocks, both with cyclic defect since the Sylow $\Phi_d$-subgroup is cyclic, with the same Brauer tree but different $\pi$-functions.
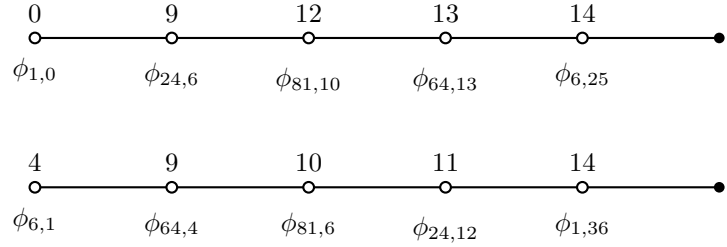
\begin{center}\begin{tikzpicture}[thick,scale=1.8]
\draw (0,1.1) -- (5,1.1);
\draw (0,1.1) node [draw,label=below:$\phi_{1,0}$] (l0) {};
\draw (1,1.1) node [draw,label=below:$\phi_{24,6}$] (l1) {};
\draw (2,1.1) node [draw,label=below:$\phi_{81,10}$] (l2) {};
\draw (3,1.1) node [draw,label=below:$\phi_{64,13}$] (l2) {};
\draw (4,1.1) node [draw,label=below:$\phi_{6,25}$] (l2) {};
\draw (5,1.1) node [fill=black!100] (ld) {};
\draw (0,1.28) node{$0$};
\draw (1,1.28) node{$9$};
\draw (2,1.28) node{$12$};
\draw (3,1.28) node{$13$};
\draw (4,1.28) node{$14$};

\draw (0,0) -- (5,0);
\draw (0,0) node [draw,label=below:$\phi_{6,1}$] (l0) {};
\draw (1,0) node [draw,label=below:$\phi_{64,4}$] (l1) {};
\draw (2,0) node [draw,label=below:$\phi_{81,6}$] (l2) {};
\draw (3,0) node [draw,label=below:$\phi_{24,12}$] (l2) {};
\draw (4,0) node [draw,label=below:$\phi_{1,36}$] (l2) {};
\draw (5,0) node [fill=black!100] (ld) {};
\draw (0,0.18) node{$4$};
\draw (1,0.18) node{$9$};
\draw (2,0.18) node{$10$};
\draw (3,0.18) node{$11$};
\draw (4,0.18) node{$14$};
\end{tikzpicture}\end{center}
When $d=8$, $9$ or $12$, we just get the Brauer trees of the principal $\ell$-blocks, which are given here, in ascending order.
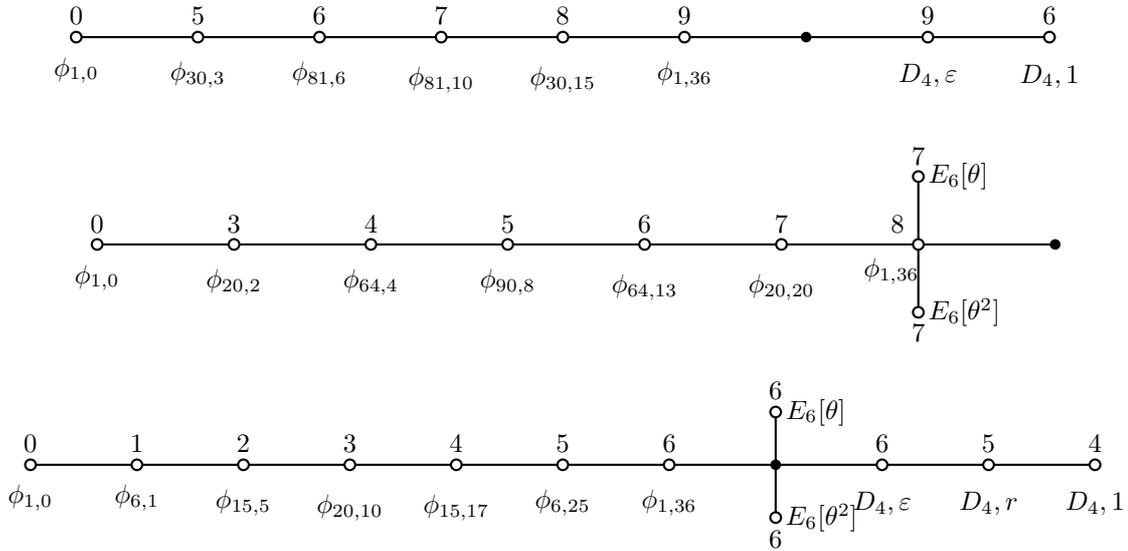
\begin{center}\begin{tikzpicture}[thick,scale=1.6]
\draw (0,0) -- (8,0);
\draw (0,0) node [draw,label=below:$\phi_{1,0}$] (l0) {};
\draw (1,0) node [draw,label=below:$\phi_{30,3}$] (l1) {};
\draw (2,0) node [draw,label=below:$\phi_{81,6}$] (l2) {};
\draw (3,0) node [draw,label=below:$\phi_{81,10}$] (l2) {};
\draw (4,0) node [draw,label=below:$\phi_{30,15}$] (l2) {};
\draw (5,0) node [draw,label=below:$\phi_{1,36}$] (l2) {};
\draw (6,0) node [fill=black!100] (ld) {};
\draw (7,0) node [draw,label=below:${D_4,\ep}$] (l2) {};
\draw (8,0) node [draw,label=below:${D_4,1}$] (l2) {};
\draw (0,0.18) node{$0$};
\draw (1,0.18) node{$5$};
\draw (2,0.18) node{$6$};
\draw (3,0.18) node{$7$};
\draw (4,0.18) node{$8$};
\draw (5,0.18) node{$9$};
\draw (7,0.18) node{$9$};
\draw (8,0.18) node{$6$};
\end{tikzpicture}\end{center}
\begin{center}\begin{tikzpicture}[thick,scale=1.8]
\draw (0,0.5) -- (7,0.5);
\draw (6,0) -- (6,1);

\draw (0,0.5) node [draw,label=below:$\phi_{1,0}$] (l0) {};
\draw (1,0.5) node [draw,label=below:$\phi_{20,2}$] (l1) {};
\draw (2,0.5) node [draw,label=below:$\phi_{64,4}$] (l2) {};
\draw (3,0.5) node [draw,label=below:$\phi_{90,8}$] (l2) {};
\draw (4,0.5) node [draw,label=below:$\phi_{64,13}$] (l2) {};
\draw (5,0.5) node [draw,label=below:$\phi_{20,20}$] (l2) {};
\draw (6,0) node [draw,label=right:${E_6[\theta^2]}$] (l4) {};
\draw (6,0.5) node [draw,label=below left:$\phi_{1,36}$] (l4) {};
\draw (6,1) node [draw,label=right:${E_6[\theta]}$] (l4) {};
\draw (7,0.5) node [fill=black!100] (ld) {};
\draw (0,0.65) node{$0$};
\draw (1,0.65) node{$3$};
\draw (2,0.65) node{$4$};
\draw (3,0.65) node{$5$};
\draw (4,0.65) node{$6$};
\draw (5,0.65) node{$7$};
\draw (5.85,0.65) node{$8$};
\draw (6,1.15) node{$7$};
\draw (6,-0.15) node{$7$};
\end{tikzpicture}\end{center}
\begin{center}\begin{tikzpicture}[thick,scale=1.4]
\draw (0,0.5) -- (10,0.5);
\draw (7,0) -- (7,1);

\draw (0,0.5) node [draw,label=below:$\phi_{1,0}$] (l0) {};
\draw (1,0.5) node [draw,label=below:$\phi_{6,1}$] (l1) {};
\draw (2,0.5) node [draw,label=below:$\phi_{15,5}$] (l2) {};
\draw (3,0.5) node [draw,label=below:$\phi_{20,10}$] (l2) {};
\draw (4,0.5) node [draw,label=below:$\phi_{15,17}$] (l2) {};
\draw (5,0.5) node [draw,label=below:$\phi_{6,25}$] (l2) {};
\draw (6,0.5) node [draw,label=below:$\phi_{1,36}$] (l4) {};
\draw (7,0) node [draw,label=right:${E_6[\theta^2]}$] (l4) {};
\draw (7,1) node [draw,label=right:${E_6[\theta]}$] (l4) {};
\draw (8,0.5) node [draw,label=below:${D_4,\ep}$] (l4) {};
\draw (9,0.5) node [draw,label=below:${D_4,r}$] (l4) {};
\draw (10,0.5) node [draw,label=below:${D_4,1}$] (l4) {};

\draw (7,0.5) node [fill=black!100] (ld) {};
\draw (0,0.7) node{$0$};
\draw (1,0.7) node{$1$};
\draw (2,0.7) node{$2$};
\draw (3,0.7) node{$3$};
\draw (4,0.7) node{$4$};
\draw (5,0.7) node{$5$};
\draw (6,0.7) node{$6$};
\draw (7,1.2) node{$6$};
\draw (7,-0.2) node{$6$};
\draw (8,0.7) node{$6$};
\draw (9,0.7) node{$5$};
\draw (10,0.7) node{$4$};
\end{tikzpicture}\end{center}
Since the $\pi$-function increases on any path leading towards the exceptional node, Theorem \ref{thm:brauertreesperverse} holds for $E_6(q)$.

\subsection{$^2\!E_6(q)$}

The group $G=E_6(q)$ has order $q^{36}\Phi_1^4\Phi_2^6\Phi_3^2\Phi_4^2\Phi_6^3\Phi_8\Phi_{10}\Phi_{12}\Phi_{18}$; the block structure can be deduced from \cite{bmm1993}; for the blocks with cyclic defect groups we give the Brauer trees, which are given in \cite{hisslubeck1998}, with the $\pi$-function attached. (Note that, for $q\equiv 1\bmod 3$, the Brauer tree is not determined completely in \cite{hisslubeck1998}, but our $\pi$-function would yield a perverse equivalence for either possibility.)

\begin{center}\begin{tabular}{l|ccc||l|ccc}
Name & $d=3$ & $d=4$ & $d=6$ & Name & $d=3$ & $d=4$ & $d=6$
\\\hline $\phi_{1,0}$&0&0&0&$\phi_{8,3}'$&14&10&$B_2$: 4
\\$\phi_{1,24}$&24&18&12&$\phi_{8,9}''$&22&16&$B_2$: 8
\\$\phi_{1,12}'$&14&10&8&$\phi_{8,3}''$&16&-&8
\\$\phi_{1,12}''$&22&16&12&$\phi_{8,9}'$&20&-&10
\\$\phi_{2,4}'$&7&5&4&$\phi_{9,2}$&-&9&7
\\$\phi_{2,16}''$&23&17&12&$\phi_{9,6}'$&-&13&-
\\$\phi_{2,4}''$&15&-&7&$\phi_{9,6}''$&-&15&-
\\$\phi_{2,16}'$&23&-&11&$\phi_{9,10}$&-&15&11
\\$\phi_{4,1}$&11&$B_2$: 6&6&$\phi_{12,4}$&-&-&9
\\$\phi_{4,13}$&23&$B_2$: 15&12&$\phi_{16,5}$&19&14&$B_2$: 7
\\$\phi_{4,8}$&20&-&10&$^2\!A_5,1$&-&-&8
\\$\phi_{4,7}'$&17&$B_2$: 10&9&$^2\!A_5,\ep$&-&-&11
\\$\phi_{4,7}''$&21&$B_2$: 13&11&$^2\!E_6[1]$&20&15&12
\\$\phi_{6,6}'$&-&13&10&$^2\!E_6[\theta]$&19&-&10
\\$\phi_{6,6}''$&-&15&-&$^2\!E_6[\theta^2]$&19&-&10
\end{tabular}
\end{center}
For $d=1$, there are two non-principal $\ell$-blocks of $^2\!E_6(q)$, both of which only have a single unipotent character, so there is automatically a perverse equivalence. For $d=4$ the tree is similar to $E_6(q)$, except that the exceptional node has moved.
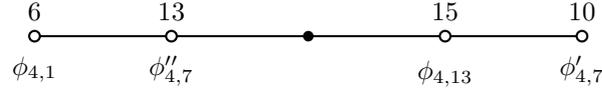
\begin{center}\begin{tikzpicture}[thick,scale=1.8]
\draw (0,0) -- (4,0);
\draw (0,0) node [draw,label=below:$\phi_{4,1}$] (l0) {};
\draw (1,0) node [draw,label=below:$\phi_{4,7}''$] (l1) {};
\draw (3,0) node [draw,label=below:$\phi_{4,13}$] (l2) {};
\draw (4,0) node [draw,label=below:$\phi_{4,7}'$] (l2) {};
\draw (2,0) node [fill=black!100] (ld) {};
\draw (0,0.18) node{$6$};
\draw (1,0.18) node{$13$};
\draw (3,0.18) node{$15$};
\draw (4,0.18) node{$10$};
\end{tikzpicture}\end{center}
For $d=6$ we get the same tree as for $E_6(q)$ and $d=3$.
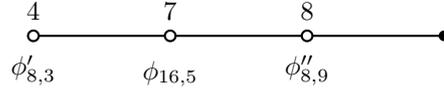
\begin{center}\begin{tikzpicture}[thick,scale=1.8]
\draw (0,0) -- (3,0);
\draw (0,0) node [draw,label=below:$\phi_{8,3}'$] (l0) {};
\draw (1,0) node [draw,label=below:$\phi_{16,5}$] (l1) {};
\draw (2,0) node [draw,label=below:$\phi_{8,9}''$] (l2) {};
\draw (3,0) node [fill=black!100] (ld) {};
\draw (0,0.18) node{$4$};
\draw (1,0.18) node{$7$};
\draw (2,0.18) node{$8$};
\end{tikzpicture}\end{center}
For $d=8$ we again get a line for the Brauer tree.
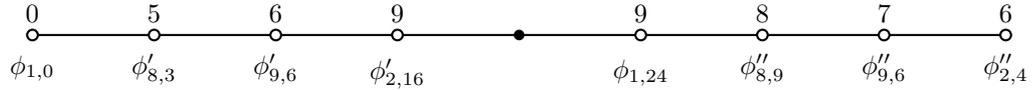
\begin{center}\begin{tikzpicture}[thick,scale=1.6]
\draw (0,0) -- (8,0);
\draw (0,0) node [draw,label=below:$\phi_{1,0}$] (l0) {};
\draw (1,0) node [draw,label=below:$\phi_{8,3}'$] (l1) {};
\draw (2,0) node [draw,label=below:$\phi_{9,6}'$] (l2) {};
\draw (3,0) node [draw,label=below:$\phi_{2,16}'$] (l2) {};
\draw (4,0) node [fill=black!100] (ld) {};
\draw (5,0) node [draw,label=below:$\phi_{1,24}$] (l2) {};
\draw (6,0) node [draw,label=below:$\phi_{8,9}''$] (l2) {};
\draw (7,0) node [draw,label=below:$\phi_{9,6}''$] (l2) {};
\draw (8,0) node [draw,label=below:$\phi_{2,4}''$] (l2) {};
\draw (0,0.18) node{$0$};
\draw (1,0.18) node{$5$};
\draw (2,0.18) node{$6$};
\draw (3,0.18) node{$9$};
\draw (5,0.18) node{$9$};
\draw (6,0.18) node{$8$};
\draw (7,0.18) node{$7$};
\draw (8,0.18) node{$6$};
\end{tikzpicture}\end{center}
In the case of $d=10$, we get two Brauer trees, as in the case of $d=5$ for $E_6(q)$.
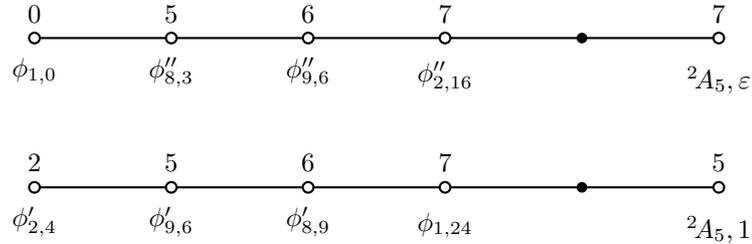
\begin{center}\begin{tikzpicture}[thick,scale=1.8]
\draw (0,1.1) -- (5,1.1);
\draw (0,1.1) node [draw,label=below:$\phi_{1,0}$] (l0) {};
\draw (1,1.1) node [draw,label=below:$\phi_{8,3}''$] (l1) {};
\draw (2,1.1) node [draw,label=below:$\phi_{9,6}''$] (l2) {};
\draw (3,1.1) node [draw,label=below:$\phi_{2,16}''$] (l2) {};
\draw (5,1.1) node [draw,label=below:${^2\!A_5,\ep}$] (l2) {};
\draw (4,1.1) node [fill=black!100] (ld) {};
\draw (0,1.28) node{$0$};
\draw (1,1.28) node{$5$};
\draw (2,1.28) node{$6$};
\draw (3,1.28) node{$7$};
\draw (5,1.28) node{$7$};

\draw (0,0) -- (5,0);
\draw (0,0) node [draw,label=below:$\phi_{2,4}'$] (l0) {};
\draw (1,0) node [draw,label=below:$\phi_{9,6}'$] (l1) {};
\draw (2,0) node [draw,label=below:$\phi_{8,9}'$] (l2) {};
\draw (3,0) node [draw,label=below:$\phi_{1,24}$] (l2) {};
\draw (5,0) node [draw,label=below:${^2\!A_5,1}$] (l2) {};
\draw (4,0) node [fill=black!100] (ld) {};
\draw (0,0.18) node{$2$};
\draw (1,0.18) node{$5$};
\draw (2,0.18) node{$6$};
\draw (3,0.18) node{$7$};
\draw (5,0.18) node{$5$};
\end{tikzpicture}\end{center}
Finally, we give the trees for $d=12$ and $d=18$; in each case there is a single Brauer tree, with a pair of non-real characters.
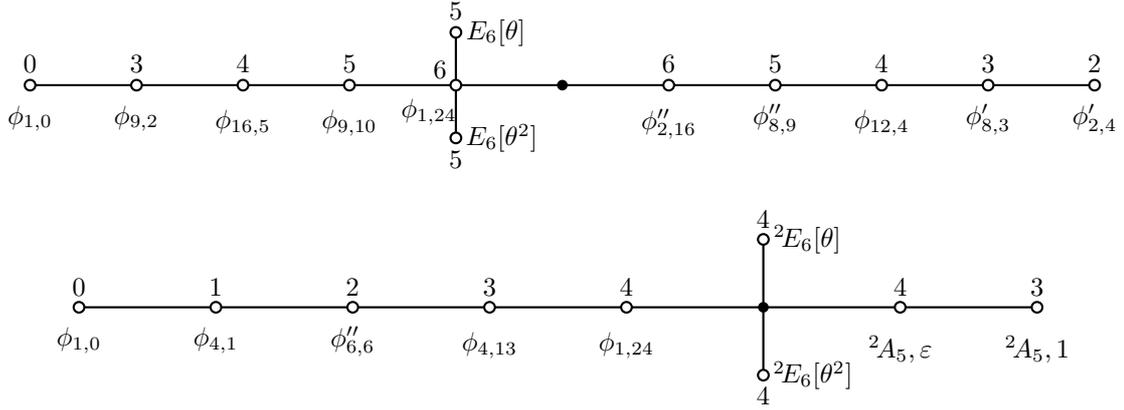
\begin{center}\begin{tikzpicture}[thick,scale=1.4]
\draw (0,0.5) -- (10,0.5);
\draw (4,0) -- (4,1);

\draw (0,0.5) node [draw,label=below:$\phi_{1,0}$] (l0) {};
\draw (1,0.5) node [draw,label=below:$\phi_{9,2}$] (l1) {};
\draw (2,0.5) node [draw,label=below:$\phi_{16,5}$] (l2) {};
\draw (3,0.5) node [draw,label=below:$\phi_{9,10}$] (l2) {};
\draw (4,0) node [draw,label=right:${E_6[\theta^2]}$] (l4) {};
\draw (4,1) node [draw,label=right:${E_6[\theta]}$] (l4) {};
\draw (4,0.5) node [draw,label=below left:$\phi_{1,24}$] (l2) {};
\draw (6,0.5) node [draw,label=below:$\phi_{2,16}''$] (l2) {};
\draw (7,0.5) node [draw,label=below:$\phi_{8,9}''$] (l4) {};
\draw (8,0.5) node [draw,label=below:$\phi_{12,4}$] (l4) {};
\draw (9,0.5) node [draw,label=below:$\phi_{8,3}'$] (l4) {};
\draw (10,0.5) node [draw,label=below:$\phi_{2,4}'$] (l4) {};

\draw (5,0.5) node [fill=black!100] (ld) {};
\draw (0,0.7) node{$0$};
\draw (1,0.7) node{$3$};
\draw (2,0.7) node{$4$};
\draw (3,0.7) node{$5$};
\draw (3.85,0.65) node{$6$};
\draw (4,1.2) node{$5$};
\draw (4,-0.2) node{$5$};
\draw (6,0.7) node{$6$};
\draw (7,0.7) node{$5$};
\draw (8,0.7) node{$4$};
\draw (9,0.7) node{$3$};
\draw (10,0.7) node{$2$};
\end{tikzpicture}\end{center}
\begin{center}\begin{tikzpicture}[thick,scale=1.8]
\draw (0,0.5) -- (7,0.5);
\draw (5,0) -- (5,1);

\draw (0,0.5) node [draw,label=below:$\phi_{1,0}$] (l0) {};
\draw (1,0.5) node [draw,label=below:$\phi_{4,1}$] (l1) {};
\draw (2,0.5) node [draw,label=below:$\phi_{6,6}''$] (l2) {};
\draw (3,0.5) node [draw,label=below:$\phi_{4,13}$] (l2) {};
\draw (4,0.5) node [draw,label=below:$\phi_{1,24}$] (l2) {};
\draw (5,0) node [draw,label=right:${^2\!E_6[\theta^2]}$] (l4) {};
\draw (5,1) node [draw,label=right:${^2\!E_6[\theta]}$] (l4) {};
\draw (5,0.5) node [fill=black!100] (ld) {};
\draw (6,0.5) node [draw,label=below:${^2\!A_5,\ep}$] (l4) {};
\draw (7,0.5) node [draw,label=below:${^2\!A_5,1}$] (l2) {};
\draw (0,0.65) node{$0$};
\draw (1,0.65) node{$1$};
\draw (2,0.65) node{$2$};
\draw (3,0.65) node{$3$};
\draw (4,0.65) node{$4$};
\draw (7,0.65) node{$3$};
\draw (6,0.65) node{$4$};
\draw (5,1.15) node{$4$};
\draw (5,-0.15) node{$4$};
\end{tikzpicture}\end{center}
Since the $\pi$-function increases on any path leading towards the exceptional node, Theorem \ref{thm:brauertreesperverse} holds for $^2\!E_6(q)$.

\section{Ree and Suzuki Groups}
\label{sec:suzree}

In order to bring the Suzuki and Ree groups into this general framework, we will have to extend our definition of $B_d(-)$ to include new polynomials, since some cyclotomic polynomials factorize over $\Q(\sqrt2)$ and $\Q(\sqrt3)$. For the groups $^2\!G_2(q^2)$ we have $\Phi_{12}=\Phi_{12}'\Phi_{12}''$, where
\[ \Phi_{12}'=q^2+\sqrt3q+1,\qquad \Phi_{12}''=q^2-\sqrt3q+1.\]
(As a warning, this labelling is the other way round from that in \cite{carterfinite}; this way round is more consistent with the theory we will give here.)  For the groups $^2\!B_2(q^2)$ and $^2\!F_4(q^2)$ we have $\Phi_8=\Phi_8'\Phi_8''$ and $\Phi_{24}=\Phi_{24}'\Phi_{24}''$, where
\[ \Phi_8'=q^2+\sqrt2q+1,\qquad \Phi_8''=q^2-\sqrt2q+1,\]
\[\Phi_{24}'=q^4+\sqrt2q^3+q^2+\sqrt2q+1,\qquad \Phi_{24}''=q^4-\sqrt2q^3+q^2-\sqrt2q+1.\]

In the previous cases, we set $\zeta$ to be $\e^{2\pi \I/d}$, which is the zero of $\Phi_d(q)$ with least argument. For $f$ not divisible by either $q-1$ or $\Phi_d$, we get that $B_d(f)$ is the sum of $\deg(f)$ and $d$ times the number of zeroes of $f$ (with multiplicity) with argument between $0$ and $2\pi/d$.

For $\alpha$ one of the $\Phi_d'$ or $\Phi_d''$, we want to define $B_\alpha(-)$, which we will write $B_{d'}(-)$ and $B_{d''}(-)$, in a similar way. For this, we need the zero of $\alpha$ of least argument, which is given in the following table.
\begin{center}\begin{tabular}{ccc}
\hline $d$&$\Phi_d'$&$\Phi_d''$
\\ \hline $8$ & $3\pi/4$ & $\pi/4$
\\ $12$ & $5\pi/6$ & $\pi/6$
\\ $24$ & $5\pi/12$ & $\pi/12$
\\ \hline\end{tabular}\end{center}
(In particular, the value for $\Phi_d''$ is the same as that of $\Phi_d$.)

With this information, we can define $B_{d'}(f)$ and $B_{d''}(f)$. For $\alpha$ one of $\Phi_d'$ or $\Phi_d''$, write the argument given in the above table as $2k\pi/d$ ($k\in\{1,3,5\}$). For $f$ not divisible by $\alpha$, $q$, or $q-1$, we define
\[ B_\alpha(f)=k\deg(f)+d\cdot|\{z\,:\,\text{$z$ is a zero of $f$ with multiplicity, and $\arg(z)\in[0,2\pi/d]$}\}|.\]
Set $B_\alpha(q)=2k$ and $B_\alpha(\Phi_1)=k+d/2$; as before, set $\pi(f)=B_d(f)/d$. This $\pi$-function extends the one produced for the other groups of Lie type, and allows us to make predictions on the cohomology of the Deligne--Lusztig varieties for the twisted groups. In most cases, these calculations have been checked by Jean Michel, and in the relevant sections we describe which calculations have been done.

For ease of use, we compile the following table of values for the $B_\alpha$-function, needed for the Suzuki and Ree groups.
\begin{center}\begin{tabular}{cccccccccc}
\hline $d$&$B_d(q)$ & $B_d(\Phi_1)$&$B_d(\Phi_2)$ & $B_d(\Phi_4)$ & $B_d(\Phi_8')$ & $B_d(\Phi_8'')$ & $B_d(\Phi_{12})$ & $B_d(\Phi_{24}')$ & $B_d(\Phi_{24}'')$
\\ \hline $8'$ & 6 & 7 & 3 & 14 & - & 14 & 20 & 20 & 28
\\ $24'$ & 10 & 17 & 5 & 10 & 10 & 34 & 44 & - & 44
\\ \hline\end{tabular}\end{center}

\begin{center}\begin{tabular}{cccccc}
\hline $d$&$B_d(q)$ & $B_d(\Phi_1)$&$B_d(\Phi_2)$ & $B_d(\Phi_4)$ & $B_d(\Phi_{12}'')$
\\ $12'$ & 10 & 11 & 5 & 22 & 14
\\ \hline\end{tabular}\end{center}

We now provide information on each group in turn, proving that the $\pi$-function as defined here produces similar answers as for the other series of groups.

\subsection{$^2\!B_2(q^2)$}

The Suzuki groups $^2\!B_2(q^2)$ have order $q^4\Phi_1\Phi_2\Phi_8$. Complete information about the Brauer trees was given in \cite{burkhardt1979}, but we use the notation for the unipotent characters from \cite{carterfinite}.

As $q$ is an odd power of $\sqrt2$, $\ell$ divides one of $q^2-1=\Phi_1\Phi_2$, $\Phi_{8}'$ or $\Phi_{8}''$. For $q^2-1$, it is not clear exactly what to use: we either set $d=2$, so that $B_2(\ep)=12$ and $\pi(\ep)=4$, or we consider $q^2-1=\Phi_1(q^2)$, so that $\ep$ has degree $(q^2)^2$, and we get $B_1((q^2)^2)=4$ and $\pi(\ep)=4$ again. Either way, we get the following tree.
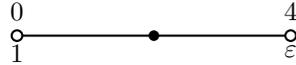
\begin{center}\begin{tikzpicture}[thick,scale=1.8]
\draw (0,0) -- (2,0);
\draw (0,0) node [draw,label=below:$1$] (l0) {};
\draw (2,0) node [draw,label=below:$\ep$] (l2) {};
\draw (1,0) node [fill=black!100] (ld) {};
\draw (0,0.18) node{$0$};
\draw (2,0.18) node{$4$};
\end{tikzpicture}\end{center}
For $\ell\mid\Phi_8'$ and $\Phi_8''$ we get the trees below.
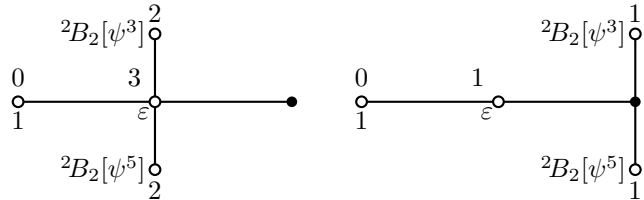
\begin{center}\begin{tikzpicture}[thick,scale=1.8]
\draw (0,0.5) -- (2,0.5);
\draw (1,0) -- (1,1);
\draw (0,0.5) node [draw,label=below:$1$] (l0) {};
\draw (1,0.5) node [draw,label=below left:$\ep$] (l2) {};
\draw (2,0.5) node [fill=black!100] (ld) {};
\draw (0,0.68) node{$0$};
\draw (0.85,0.68) node{$3$};

\draw (1,1) node [draw,label=left:${^2\!B_2[\psi^3]}$] (l0) {};
\draw (1,0) node [draw,label=left:${^2\!B_2[\psi^5]}$] (l0) {};

\draw (1,1.15) node{$2$};
\draw (1,-0.15) node{$2$};
\end{tikzpicture}\qquad\begin{tikzpicture}[thick,scale=1.8]
\draw (0,0.5) -- (2,0.5);
\draw (2,0) -- (2,1);
\draw (0,0.5) node [draw,label=below:$1$] (l0) {};
\draw (1,0.5) node [draw,label=below left:$\ep$] (l2) {};
\draw (2,0.5) node [fill=black!100] (ld) {};
\draw (0,0.68) node{$0$};
\draw (0.85,0.68) node{$1$};

\draw (2,1) node [draw,label=left:${^2\!B_2[\psi^3]}$] (l0) {};
\draw (2,0) node [draw,label=left:${^2\!B_2[\psi^5]}$] (l0) {};

\draw (2,1.15) node{$1$};
\draw (2,-0.15) node{$1$};
\end{tikzpicture}\end{center}
Each of these cases has been verified, in calculations by Jean Michel.

\subsection{$^2\!G_2(q^2)$}
The groups $^2\!G_2(q^2)$ have order $q^6\Phi_1\Phi_2\Phi_4\Phi_{12}$. The Brauer trees were determined in \cite{hiss1991}, and since no cyclotomic polynomial divides $|G|$ to more than the first power, this is all the information that we need. Because we have labelled $\Phi_{12}'$ and $\Phi_{12}''$  in a different way to \cite{carterfinite}, we give a table of the unipotent characters, together with a labelling of the cuspidal characters.
\begin{center}\begin{tabular}{ll}
Name & Degree
\\ \hline $1$ & $1$
\\ $\ep$ & $q^6$
\\ ${}^2\!G_2[\xi^5]$ & $q\Phi_1\Phi_2\Phi_4/\sqrt3$
\\ ${}^2\!G_2[\xi^7]$ & $q\Phi_1\Phi_2\Phi_4/\sqrt3$
\\ ${}^2\!G_2^\mathrm{I}[\I]$ & $q\Phi_1\Phi_2\Phi_{12}'/2\sqrt3$
\\ ${}^2\!G_2^\mathrm{I}[-\I]$ & $q\Phi_1\Phi_2\Phi_{12}'/2\sqrt3$
\\ ${}^2\!G_2^\mathrm{II}[\I]$ & $q\Phi_1\Phi_2\Phi_{12}''/2\sqrt3$
\\ ${}^2\!G_2^\mathrm{II}[-\I]$ & $q\Phi_1\Phi_2\Phi_{12}''/2\sqrt3$
\end{tabular}\end{center}
As $q$ is an odd power of $\sqrt3$, we have that $\ell$ divides one of $(q^2-1)$, $\Phi_4$, $\Phi_{12}'$ or $\Phi_{12}''$. For $\ell\mid(q^2-1)$, we use the same idea as for $^2\!B_2(q^2)$, and so $\pi(\ep)=6$.
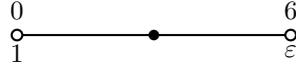
\begin{center}\begin{tikzpicture}[thick,scale=1.8]
\draw (0,0) -- (2,0);
\draw (0,0) node [draw,label=below:$1$] (l0) {};
\draw (2,0) node [draw,label=below:$\ep$] (l2) {};
\draw (1,0) node [fill=black!100] (ld) {};
\draw (0,0.18) node{$0$};
\draw (2,0.18) node{$6$};
\end{tikzpicture}\end{center}
For $\ell\mid\Phi_4(q)$, we get $B_4(\Phi_{12}')=4$ and $B_4(\Phi_{12}'')=8$, using the formula given earlier in this section, and get the following Brauer tree.
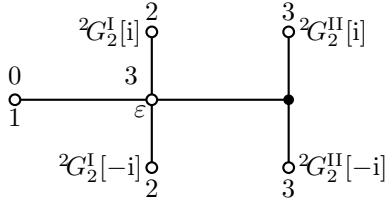
\begin{center}\begin{tikzpicture}[thick,scale=1.8]
\draw (0,0.5) -- (2,0.5);
\draw (1,0) -- (1,1);
\draw (2,0) -- (2,1);
\draw (0,0.5) node [draw,label=below:$1$] (l0) {};
\draw (1,0.5) node [draw,label=below left:$\ep$] (l2) {};
\draw (2,0.5) node [fill=black!100] (ld) {};
\draw (0,0.68) node{$0$};
\draw (0.85,0.68) node{$3$};

\draw (1,1) node [draw,label=left:${{}^2\!G_2^\mathrm{I}[\I]}$] (l0) {};
\draw (1,0) node [draw,label=left:${{}^2\!G_2^\mathrm{I}[-\I]}$] (l0) {};

\draw (2,1) node [draw,label=right:${{}^2\!G_2^\mathrm{II}[\I]}$] (l0) {};
\draw (2,0) node [draw,label=right:${{}^2\!G_2^\mathrm{II}[-\I]}$] (l0) {};

\draw (1,1.15) node{$2$};
\draw (1,-0.15) node{$2$};
\draw (2,1.15) node{$3$};
\draw (2,-0.15) node{$3$};
\end{tikzpicture}\end{center}
When $\ell\mid\Phi_{12}'$ and $\ell\mid\Phi_{12}''$, we get the following trees. (The Brauer tree for $\Phi_{12}''$ is from \cite{hiss1991}, but the planar embedding is determined in \cite{dudas2010un2}.)
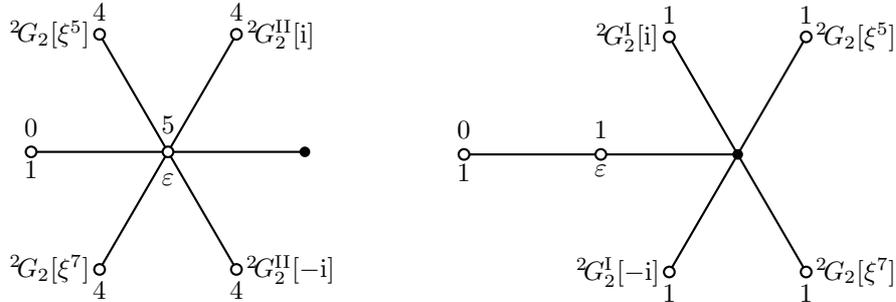
\begin{center}\begin{tikzpicture}[thick,scale=1.8]
\draw (0,1) -- (2,1);
\draw (0.5,0.1339) -- (1.5,1.866);
\draw (0.5,1.866) -- (1.5,0.1339);
\draw (0,1) node [draw,label=below:$1$] (l0) {};
\draw (1,1) node [draw] (l2) {};
\draw (2,1) node [fill=black!100] (ld) {};
\draw (0,1.18) node{$0$};
\draw (1,1.2) node{$5$};

\draw (1,0.8) node{$\ep$};

\draw (0.5,1.866) node [draw,label=left:${{}^2\!G_2[\xi^5]}$] (l0) {};
\draw (1.5,1.866) node [draw,label=right:${{}^2\!G_2^\mathrm{II}[\I]}$] (l0) {};
\draw (0.5,0.1339) node [draw,label=left:${{}^2\!G_2[\xi^7]}$] (l0) {};
\draw (1.5,0.1339) node [draw,label=right:${{}^2\!G_2^\mathrm{II}[-\I]}$] (l0) {};

\draw (0.5,2.03) node{$4$};
\draw (1.5,2.03) node{$4$};
\draw (0.5,-0.02) node{$4$};
\draw (1.5,-0.02) node{$4$};
\end{tikzpicture}
\qquad\qquad\begin{tikzpicture}[thick,scale=1.8]
\draw (-1,1) -- (1,1);
\draw (0.5,0.1339) -- (1.5,1.866);
\draw (0.5,1.866) -- (1.5,0.1339);
\draw (-1,1) node [draw,label=below:$1$] (l0) {};

\draw (0,1) node [draw,label=below:$\ep$] (l0) {};
\draw (1,1) node [fill=black!100] (ld) {};
\draw (-1,1.18) node{$0$};
\draw (0,1.18) node{$1$};

\draw (1.5,1.866) node [draw,label=right:${{}^2\!G_2[\xi^5]}$] (l0) {};
\draw (0.5,1.866) node [draw,label=left:${{}^2\!G_2^\mathrm{I}[\I]}$] (l0) {};
\draw (1.5,0.1339) node [draw,label=right:${{}^2\!G_2[\xi^7]}$] (l0) {};
\draw (0.5,0.1339) node [draw,label=left:${{}^2\!G_2^\mathrm{I}[-\I]}$] (l0) {};

\draw (0.5,2.03) node{$1$};
\draw (1.5,2.03) node{$1$};
\draw (0.5,-0.02) node{$1$};
\draw (1.5,-0.02) node{$1$};
\end{tikzpicture}\end{center}

In calculations by Jean Michel, this has been confirmed in all cases.

\subsection{$^2\!F_4(q^2)$}

The groups $^2\!F_4(q^2)$ have order $q^{24}\Phi_1^2\Phi_2^2\Phi_4^2\Phi_8^2\Phi_{12}\Phi_{24}$. The Brauer trees appeared in \cite{hiss1991}, and considerable information on the decomposition matrices given in \cite{himstedt2011}. In this case, if $\ell\nmid q$ then $\ell$ divides one of $q^2-1$, $\Phi_4$, $\Phi_8'$, $\Phi_8''$, $\Phi_{12}$, $\Phi_{24}'$ or $\Phi_{24}''$. Because there are misprints in the table of degrees in \cite{carterfinite}, we give the degrees here.
\begin{center}\begin{tabular}{ll|cccc}
Name & Degree & $d=4$ & $d=8'$ & $d=8''$
\\ \hline $\phi_{1,0}=\chi_1=1$ & $1$ & 0 & 0 & 0
\\ $\phi_{1,4}''=\chi_4=\ep'$ & $q^2\Phi_{12}\Phi_{24}$ & 7 & 10 & 4
\\ $\phi_{1,4}'=\chi_{18}=\ep''$ & $q^{10}\Phi_{12}\Phi_{24}$ & 11 & 16& 6
\\ $\phi_{1,8}=\chi_{21}=\ep$ & $q^{24}$ & 12 & 18 & 6
\\ $\phi_{2,3}=\chi_5=\rho_2'$ & $q^4\Phi_4^2\Phi_8''^2\Phi_{12}\Phi_{24}'/4$ & - & 15 & -
\\ $\phi_{2,1}=\chi_6=\rho_2''$ & $q^4\Phi_4^2\Phi_8'^2\Phi_{12}\Phi_{24}''/4$ & - & - & 5
\\ $\phi_{2,2}=\chi_7=\rho_2$ & $q^4\Phi_8^2\Phi_{24}/2$ & 10 & - & -
\\ $\chi_2={}^2\!B_2[\psi^3],1$ & $q\Phi_1\Phi_2\Phi_4^2\Phi_{12}/\sqrt2$ & - & 8 & 3
\\ $\chi_3={}^2\!B_2[\psi^5],1$ & $q\Phi_1\Phi_2\Phi_4^2\Phi_{12}/\sqrt2$ & - & 8 & 3
\\ $\chi_{19}={}^2\!B_2[\psi^3],\ep$ & $q^{13}\Phi_1\Phi_2\Phi_4^2\Phi_{12}/\sqrt2$ & - & 17 & 6
\\ $\chi_{20}={}^2\!B_2[\psi^5],\ep$ & $q^{13}\Phi_1\Phi_2\Phi_4^2\Phi_{12}/\sqrt2$ & - & 17 & 6
\\ ${}^2F_4^{\mathrm{I}}[-1]=\chi_{10}$ & $q^4\Phi_1^2\Phi_2^2\Phi_4^2\Phi_{24}/6$ & - & 15 & 5
\\ ${}^2\!F_4^{\mathrm{II}}[-1]=\chi_9$& $q^4\Phi_1^2\Phi_2^2\Phi_8''^2\Phi_{12}\Phi_{24}''/12$ & 11 & 15 & -
\\ ${}^2\!F_4^{\mathrm{III}}[-1]=\chi_8$ & $q^4\Phi_1^2\Phi_2^2\Phi_8'^2\Phi_{12}\Phi_{24}'/12$ & 9 & - & 5
\\ ${}^2\!F_4^{\mathrm{IV}}[-1]=\chi_{17}$&$q^4\Phi_1^2\Phi_2^2\Phi_{12}\Phi_{24}/3$ & 10 & 14 & 6
\\ ${}^2\!F_4^{\mathrm{I}}[\I]=\chi_{13}$ & $q^4\Phi_1^2\Phi_2^2\Phi_4^2\Phi_{12}\Phi_{24}'/4$ & - & 14 & 5
\\ ${}^2\!F_4^{\mathrm{II}}[\I]=\chi_{11}$ & $q^4\Phi_1^2\Phi_2^2\Phi_4^2\Phi_{12}\Phi_{24}''/4$ & - & 15 & 6
\\ ${}^2\!F_4^{\mathrm{I}}[-\I]=\chi_{14}$ & $q^4\Phi_1^2\Phi_2^2\Phi_4^2\Phi_{12}\Phi_{24}'/4$ & - & 14 & 5
\\ ${}^2\!F_4^{\mathrm{II}}[-\I]=\chi_{12}$ & $q^4\Phi_1^2\Phi_2^2\Phi_4^2\Phi_{12}\Phi_{24}''/4$ & - & 15 & 6
\\ ${}^2\!F_4[-\theta]=\chi_{15}$ & $q^4\Phi_1^2\Phi_2^2\Phi_4^2\Phi_8^2/3$ & - & - & -
\\ ${}^2\!F_4[-\theta^2]=\chi_{16}$ & $q^4\Phi_1^2\Phi_2^2\Phi_4^2\Phi_8^2/3$ & - & - & -
\end{tabular}\end{center}

For $\ell\mid(q^2-1)$ we get the following two trees, easily deducible although not described in \cite{hiss1991}.
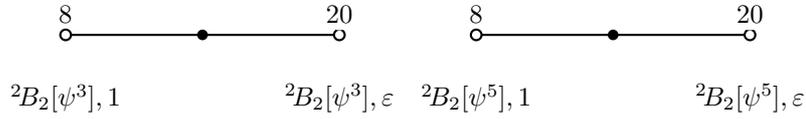
\begin{center}\begin{tikzpicture}[thick,scale=1.8]
\draw (0,0) -- (2,0);
\draw (3,0) -- (5,0);

\draw (0,0) node [draw,label=below:${^2\!B_2[\psi^3],1}$] (l0) {};
\draw (2,0) node [draw,label=below:${^2\!B_2[\psi^3],\ep}$] (l1) {};
\draw (3,0) node [draw,label=below:${^2\!B_2[\psi^5],1}$] (l2) {};
\draw (5,0) node [draw,label=below:${^2\!B_2[\psi^5],\ep}$] (l2) {};

\draw (1,0) node [fill=black!100] (ld) {};
\draw (4,0) node [fill=black!100] (ld) {};

\draw (0,0.15) node{$8$};
\draw (2,0.15) node{$20$};
\draw (3,0.15) node{$8$};
\draw (5,0.15) node{$20$};
\end{tikzpicture}\end{center}
An easy calculation yields $B_4(\Phi_{24}')=B_4(\Phi_{24}'')=8$, using the rule given in the introduction to this section. This yields the integers $\pi(-)$ above for $d=4$, producing a triangular decomposition matrix for the principal $\ell$-block in both $d=4$ and $d=8'$, according to \cite{himstedt2011}, although the matrix is sparse so it is not surprising. For $d=8''$ this requires some currently unknown parameters to be $0$.

If $\ell\mid\Phi_{12}$ then we get the following labelled Brauer tree.
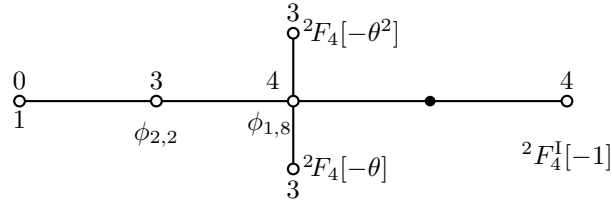
\begin{center}\begin{tikzpicture}[thick,scale=1.8]
\draw (0,0.5) -- (4,0.5);
\draw (2,0) -- (2,1);

\draw (0,0.5) node [draw,label=below:$1$] (l0) {};
\draw (1,0.5) node [draw,label=below:$\phi_{2,2}$] (l1) {};
\draw (2,0.5) node [draw,label=below left:$\phi_{1,8}$] (l2) {};
\draw (4,0.5) node [draw,label=below:${{}^2F_4^{\mathrm{I}}[-1]}$] (l2) {};

\draw (2,0) node [draw,label=right:${{}^2\!F_4[-\theta]}$] (l4) {};
\draw (2,1) node [draw,label=right:${{}^2\!F_4[-\theta^2]}$] (l4) {};

\draw (3,0.5) node [fill=black!100] (ld) {};

\draw (0,0.65) node{$0$};
\draw (1,0.65) node{$3$};
\draw (4,0.65) node{$4$};
\draw (1.85,0.65) node{$4$};

\draw (2,1.15) node{$3$};
\draw (2,-0.15) node{$3$};
\end{tikzpicture}\end{center}
If $\ell\mid\Phi_{24}'$, then the $\pi$-function is the following.
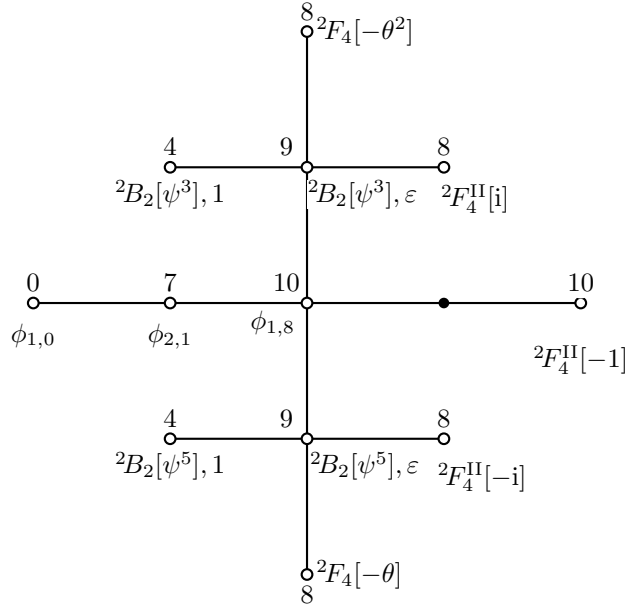
\begin{center}\begin{tikzpicture}[thick,scale=1.8]
\draw (0,2) -- (4,2);
\draw (1,1) -- (3,1);
\draw (1,3) -- (3,3);

\draw (2,0) -- (2,4);

\draw (0,2) node [draw,label=below:$\phi_{1,0}$] (l0) {};
\draw (1,2) node [draw,label=below:$\phi_{2,1}$] (l1) {};
\draw (2,2) node [draw] {};
\draw (3,2) node [fill=black!100] (ld) {};
\draw (4,2) node [draw,label=below:${{}^2\!F_4^{\mathrm{II}}[-1]}$] (l2) {};

\draw (1,1) node [draw] {};
\draw (2,1) node [draw] (l4) {};
\draw (3,1) node [draw,label=below right:${{}^2\!F_4^{\mathrm{II}}[-\I]}$] (l4) {};

\draw (2,0) node [draw,label=right:${{}^2\!F_4[-\theta]}$] (l4) {};

\draw (1,3) node [draw] {};

\draw (2,3) node [draw] (l4) {};
\draw (3,3) node [draw,label=below right:${{}^2\!F_4^{\mathrm{II}}[\I]}$] (l4) {};

\draw (2,4) node [draw,label=right:${{}^2\!F_4[-\theta^2]}$] (l4) {};

\draw (1,2.80) node[rectangle] {$^2\!B_2[\psi^3],1$};
\draw (1,0.80) node[rectangle] {$^2\!B_2[\psi^5],1$};

\draw (2.405,2.80) node[rectangle] {$^2\!B_2[\psi^3],\ep$};
\draw (2.422,0.80) node[rectangle] {$^2\!B_2[\psi^5],\ep$};

\draw (1.75,1.85) node[rectangle] {$\phi_{1,8}$};

\draw (0,2.15) node{$0$};
\draw (1,2.15) node{$7$};
\draw (1.85,2.15) node{$10$};
\draw (4,2.15) node{$10$};

\draw (3,3.15) node{$8$};
\draw (3,1.15) node{$8$};
\draw (1,3.15) node{$4$};
\draw (1,1.15) node{$4$};
\draw (2,4.15) node{$8$};
\draw (2,-0.15) node{$8$};

\draw (1.85,3.15) node{$9$};
\draw (1.85,1.15) node{$9$};

\end{tikzpicture}\end{center}

Finally, if $\ell\mid\Phi_{24}''$, we get the standard picture for the Coxeter case. (Note that the planar embedding is not known, and that this is just the guess from a conjecture of Hiss--L\"ubeck--Malle \cite{hlm1995}.
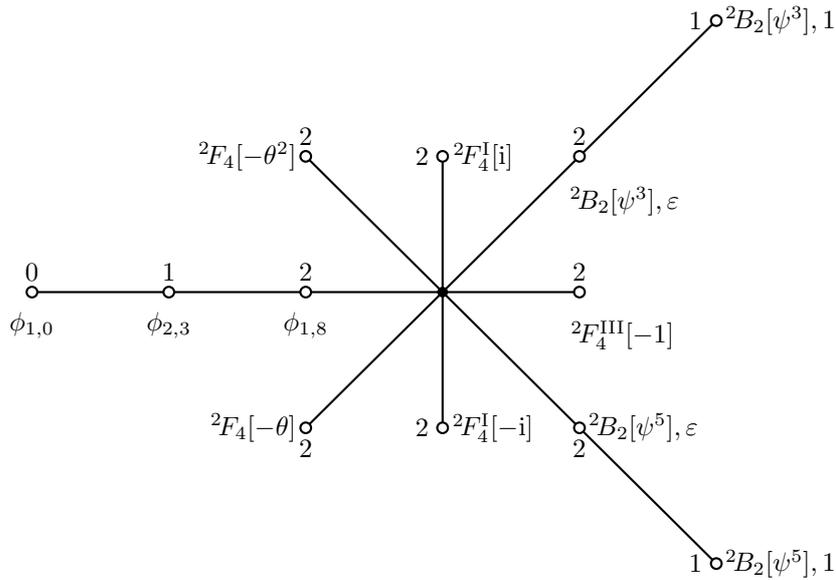
\begin{center}\begin{tikzpicture}[thick,scale=1.8]
\draw (0,2) -- (4,2);
\draw (3,1) -- (3,3);
\draw (2,1) -- (5,4);
\draw (2,3) -- (5,0);

\draw (0,2) node [draw,label=below:$\phi_{1,0}$] (l0) {};
\draw (1,2) node [draw,label=below:$\phi_{2,3}$] (l1) {};
\draw (2,2) node [draw,label=below:$\phi_{1,8}$] (l1) {};
\draw (3,2) node [fill=black!100] (ld) {};
\draw (4,2) node [draw,label=below right:${{}^2\!F_4^{\mathrm{III}}[-1]}$] (l2) {};

\draw (5,0) node [draw,label=right:${{}^2\!B_2[\psi^5],1}$] (l4) {};
\draw (4,1) node [draw,label=right:${{}^2\!B_2[\psi^5],\ep}$] (l4) {};
\draw (4,3) node [draw,label=below right:${{}^2\!B_2[\psi^3],\ep}$] (l4) {};
\draw (5,4) node [draw,label=right:${{}^2\!B_2[\psi^3],1}$] (l4) {};
\draw (2,1) node [draw,label=left:${{}^2\!F_4[-\theta]}$] (l4) {};
\draw (3,1) node [draw,label=right:${{}^2\!F_4^{\mathrm{I}}[-\I]}$] (l4) {};
\draw (3,3) node [draw,label=right:${{}^2\!F_4^{\mathrm{I}}[\I]}$] (l4) {};
\draw (2,3) node [draw,label=left:${{}^2\!F_4[-\theta^2]}$] (l4) {};
%

%
%
%
%
%
\draw (0,2.15) node{$0$};
\draw (1,2.15) node{$1$};
\draw (2,2.15) node{$2$};
\draw (4,2.15) node{$2$};

\draw (4.85,0) node{$1$};
\draw (2.85,1) node{$2$};
\draw (2.85,3) node{$2$};
\draw (4.85,4) node{$1$};

\draw (2,0.85) node{$2$};
\draw (2,3.15) node{$2$};
\draw (4,0.85) node{$2$};
\draw (4,3.15) node{$2$};
%

\end{tikzpicture}\end{center}

\bigskip\bigskip

\noindent\textbf{Acknowledgements}: I would like to thank Gunter Malle for suggesting that I consider the $d$-cuspidal pair in a preliminary version of Conjecture \ref{conj:DLcohom}. I thank Olivier Dudas for giving me access to various unpublished calculations on the cohomology of Deligne--Lusztig varieties, which I used to verify Conjecture \ref{conj:DLcohom} in Section \ref{sec:previouswork}, for reading and commenting on this manuscript, and very helpful conversations. Jean Michel also performed several calculations with Deligne--Lusztig varieties that backed up Conjecture \ref{conj:DLcohom}, and more importantly helped me to straighten out the problems I was having with the `very twisted' Ree and Suzuki groups. I would finally like to thank Rapha\"el Rouquier for letting me bounce ideas off him, and for introducing me to the problem of determining the `geometric perversity function'.

\bibliography{references}

\end{document}